\documentclass{ejm}
\usepackage{amsmath}
\usepackage{amssymb}
\usepackage{graphicx}
\usepackage[lined,boxed]{algorithm2e}
\usepackage{gensymb}
\usepackage{fancybox}
\usepackage{listings}
\usepackage{color}
\usepackage{enumerate}
\usepackage{marginnote}
\usepackage{breqn}
\usepackage{mathrsfs}
\usepackage{cite}

\newcommand{\red}[1]{\ifmmode\mathbf{\textcolor{red}{#1}}\else \textbf{\textcolor{red}{#1}}\fi}

\newcommand{\blue}[1]{\ifmmode\mathbf{\textcolor{blue}{#1}}\else \textbf{\textcolor{blue}{#1}}\fi}
\newcommand{\imagi}{\mathrm{i}}

\newcommand{\gray}[1]{\ifmmode\mathbf{\textcolor{gray}{#1}}\else \textbf{\textcolor{gray}{#1}}\fi}

\newcommand{\lbl}[1]{\label{#1}}

\def\K{K}

\def\eps{\epsilon}
\def\hmin{h_{\min}}
\def\peps{\mathcal{P}_{\max}}
\def\A{\mathcal{A}}
\def\p0{\bar{\mathcal{P}}}

\def\hpeak{h_{\mathrm{peak}}}
\def\fld{\mathrm{fold}}

\def\HA{H_A} % second-order steady state
\def\HB{H_B} % fourth-order steady state
\def\HLC{H_C} % limit cycle soln
\def\Hc{\bar{H}_+} %original H_c, larger uniform steady state
\def\Hm{\bar{H}_-} %original H_m, smaller uniform steady state

\def\lB{{\bar{\ell}_{-}}^{\,\gamma}}%original \ell_B, bifurcation point for domain size L where the primary fourth-order steady state bifurcates from \Hm uniform steady state

\def\lA{{\bar{\ell}_{+}}^{\,p}}%original \ell_A, bifurcation point for domain size L where the primary second-order steady state bifurcates from \Hc uniform steady state

\def\lBc{{\bar{\ell}_{+}}^{\,\gamma}}%original \ell_B^c, bifurcation point for domain size L where the primary fourth-order steady state bifurcates from \Hc uniform steady state

\def\MM{\mathcal{M}}
\def\EE{\mathcal{E}}
\def\ts{\textstyle}
\def\Mc{M_{\mathrm{c}}}
\def\Mnc{M_{\mathrm{nc}}}
\usepackage{hyperref}
\hypersetup{
    colorlinks,
    citecolor=black,
    filecolor=black,
    linkcolor=black,
    urlcolor=black}

\begin{document}
%%%%%%%%%%%%%%%%%%%%%%%%%%%%%%%%%%%%%
%  2. the title and author, etc.
%%%%%%%%%%%%%%%%%%%%%%%%%%%%%%%%%%%%%

%\begin{frontmatter}

\title{Steady states and dynamics of a thin film-type equation with non-conserved mass}
%\shorttitle{Steady states and dynamics of a non-conserved mass thin film equation}
\shorttitle{A thin film equation with non-conserved mass}
\author[H. Ji and T. P. Witelski]{%
    H.\ns J\ls I$\,^1$\ns
    \and
    T.\ns P.\ns W\ls I\ls T\ls E\ls L\ls S\ls K\ls I $\,^2$
}
\affiliation{%
    $^1\,$Department of Mathematics, University of California, Los Angeles, Los Angeles, CA 90095, USA\\
      email\textup{\nocorr: \texttt{hangjie@math.ucla.edu}}\\
     $^2\,$Department of Mathematics, Duke University, Durham, NC 27708-0320, USA\\
      email\textup{\nocorr: \texttt{witelski@math.duke.edu}}
}

\maketitle
\begin{abstract}
We study the steady states and dynamics of a thin film-type equation with non-conserved mass in one dimension. The evolution equation is a nonlinear fourth-order degenerate parabolic PDE motivated by a model of volatile viscous fluid films allowing for condensation or evaporation.
We show that by changing the sign of the non-conserved flux and breaking from a gradient flow structure, the problem can exhibit novel behaviors including having two distinct classes of coexisting steady state solutions. Detailed analysis of the bifurcation structure for these steady states and their stability reveal several possibilities for the dynamics.
For some parameter regimes, solutions can lead to finite-time rupture singularities. Interestingly, we also show that a finite amplitude limit cycle can occur as a singular perturbation in the nearly-conserved limit.
\end{abstract}


\begin{section}{Introduction}
\par
Phase field models are widely used in many branches of continuum mechanics to describe pattern formation and evolving interfaces with respect to an order parameter $\phi({\bf x},t)$ that could represent a volume fraction of one phase in a binary mixture \cite{grant1993spinodal,ANC1984,fife2000}, or the density of healthy tissue in a biological system \cite{garcke2017well}, for example. 
\par
Such models can be formulated in terms of a free energy, which in simple cases will have contributions from a potential energy of homogeneous phases, $U(\phi)$, and an interfacial energy from the formation of gradients,
$$
\EE[\phi] = \iint_D U(\phi) + {\ts {1\over 2}}|\nabla \phi|^2\,dA.
$$
A chemical potential can then be defined from the variational derivative of the energy,
$$
\mu[\phi]\equiv  {\delta \EE\over \delta \phi}=U'(\phi) -\nabla^2 \phi.
$$
As described by Thiele in \cite{thiele2010thin,thiele2016}, a general dissipative evolution equation for $\phi$ can then be written as 
\begin{equation}
\frac{\partial \phi}{\partial t} = \nabla\cdot \left({\Mc (\phi)\nabla \mu[\phi]}\right) - \Mnc (\phi)\mu[\phi],
\lbl{generalModel}
\end{equation} 
where $\Mc (\phi)$ and $\Mnc (\phi)$ are non-negative mobility functions for the conserved and non-conserved parts of the dynamics, respectively. These coefficients define a spatial flux following Fick's law, ${\bf J} = -\Mc \nabla \mu$, and a mass loss rate,
$$
\MM=\iint_D \phi\,dA\qquad {d\MM\over dt} = -\iint_D \Mnc \mu\,dA.
$$
This model monotonically dissipates the energy as a Lyapnunov functional,
\begin{equation}
\frac{d\EE}{dt} = -\iint_D \Mc |\nabla \mu|^2 + \Mnc \mu^2 \,dA\le 0,
\lbl{generalEdiss}
\end{equation}
which shows that all equilibrium states must correspond to constant chemical potentials, and further if mass is not conserved ($\Mnc\neq 0$), they must have $\mu\equiv 0$.
We note that the $U(\phi)$ term in the energy may include both an intrinsic potential energy  $\tilde{U}(\phi)$ and a linear term, $U(\phi) = \tilde{U}(\phi)-\mu_0 \phi$. This leads to $\mu[\phi] = \tilde{\mu}[\phi]-\mu_0$, where $\tilde{\mu}$ represents an intrinsic chemical potential, so that for equilibrium states in the non-conserved system the condition $\mu = 0$ corresponds to a constant intrinsic potential with $\tilde{\mu} = \mu_0$.
\par 
One model fitting into this framework is
the mixed Allen-Cahn/Cahn-Hilliard equation (AC/CH) \cite{miranville2017cahn, israel2013well,israel2013wella, karali2012existence} that describes microscopic pattern formation mechanisms like surface diffusion, adsorption and desorption \cite{karali2007role}. 
The basic form of this model has constants for both mobilities $M_c$ and $M_{nc}$ ($\Mc \equiv D >0$ and $\Mnc\equiv 1$) and the potential energy functional is a symmetric double well, $U(\phi) = (1-\phi^2)^2$, written in one spatial dimension as
$$
{\partial \phi\over \partial t} = D \frac{\partial^2 \mu}{\partial x^2} - \mu, \qquad
\mu =  {U'(\phi)}-\frac{\partial^2 \phi}{\partial x^2}.
$$
The classic Cahn-Hilliard equation \cite{ANC1984,ANC1988} for phase separation of binary mixtures is recovered by eliminating the non-conserved flux,
$\Mnc \equiv 0$.
Many variants of the mixed AC/CH model have been used to study more complex interfacial dynamics, one such example is the coupled system describing a model of tumor growth investigated in \cite{garcke2017well}.  Studies have also addressed mathematical properties of solutions when the mobility coefficient is a degenerate function of $\phi$ in the Cahn-Hilliard equation \cite{dai2016} and for the AC/CH model \cite{zhang2016existence}. There is a vast body of literature on existence, regularity and attracting states for the Cahn-Hilliard equation \cite{miranville2017cahn}. Some results on steady states in convective Cahn-Hilliard models subject to external fields have also been obtained \cite{munch2008}. Fewer results have been obtained specifically for the AC/CH equation \cite{israel2013well,karali2010,karali2012existence,zhang2016existence}, but many extensions with systems of equations having phase-field models coupled to other physical effects have been used in applications \cite{chen2002}. The Cahn-Hilliard equation with other forms of non-conserved terms (not proportional to the chemical potential) have been used to describe linear adsorption/desorption (the Cahn-Hilliard-Oono equation) and other applications \cite{miranville2017cahn}, but will not have an energy dissipation result like \eqref{generalEdiss}.
\par
While the derivation of lubrication models for coating flows of thin viscous films on solid substrates follows from an asymptotic reduction of the Navier-Stokes equations for low Reynolds numbers \cite{oron1997long,craster2009dynamics,myers,OckendonOckendon}, the resulting thin film equations also fit the form (\ref{generalModel}) \cite{Mitlin,thiele2010thin,Thiele2014,thiele2016,Thiele2018}, with the order parameter representing the film height, $\phi\to h(x,t)\ge 0$. Likewise, the role of the chemical potential is taken on by the hydrodynamic pressure $\mu\to p$, which will have contributions from surface tension and a hydrostatic potential function.  
\par
One class of models derived for evaporating thin films of volatile fluids on heated substrates \cite{lyushnin2002fingering,ajaev2005spreading, ajaev2005evolution, ajaev2001steady,moosman1980evaporating}
is consistent with the form of equation \eqref{generalModel}, also see references listed in \cite{Thiele2014}.
For Newtonian fluids on solid substrates with no-slip boundary conditions, the mobility functions for the conserved and non-conserved terms are $\Mc(h) = h^3$ and $\Mnc(h) = \beta/(h+\K)$, where the constants $\beta> 0$ and $\K>0$
are set by material properties of the fluid and a thermodynamic-kinetic condition at the fluid-vapor interface \cite{ajaev2005evolution}. Alternative forms for the mobility $M_c(h)$ have been discussed for slip models in \cite{munch2005lubrication}.
The $M_c(h)$ function makes thin film models comparable to degenerate-mobility Cahn-Hilliard equations \cite{garcke1996,dai2016,zhang2016existence}.
The resulting model is
\begin{equation}
{\partial h\over \partial t} = 
{\partial\over \partial x}
\left(h^3 {\partial p\over \partial x}\right)-\frac{\beta}{h+\K}\, p , \qquad p = \Pi(h) - 
{\partial^2 h\over \partial x^2},
\lbl{evapModel}
\end{equation}
where $\Pi(h)=U'(h)$  is a generalized disjoining pressure that incorporates both spatial wetting properties of the substrate and the thermodynamics driving evaporation, $\Pi(h) = \tilde{\Pi}(h)-\p0$ \cite{ji2018instability}.
The contribution from the standard disjoining pressure function $\tilde{\Pi}(h)$ gives the pressure due to the hydrophobic or hydrophilic properties of the substrate, as typically described in relation to contact angles of droplets \cite{glasner2003pf,de1985wetting,teletzke1988wetting}. For partially-wetting liquids, the form $\tilde{\Pi}(h) = h^{-3}-h^{-4}$ has been frequently used and provides a positive lower bound for the film thickness \cite{bertozzi2001dewetting,glasner2003coarsening,ji2018instability}. 
The pressure offset  $\p0$ then encodes the difference between the fluid temperature and the surrounding vapor, determining whether evaporation or condensation is favored \cite{ji2018instability}. 
For $\beta\ge 0$
since this model matches \eqref{generalModel} it will have an energy dissipation equation of the form \eqref{generalEdiss}.
\par
In the special case $\beta=0$, the model \eqref{evapModel} reduces to the mass-conserving thin film equation where no evaporation or condensation occurs, 
\begin{equation}
{\partial h\over \partial t} = 
{\partial\over \partial x}
\left(h^3 {\partial\over \partial x}\left[\Pi(h) - 
{\partial^2 h\over \partial x^2}\right]\right),
\lbl{mainnoevap}
\end{equation}
which has been studied extensively in mathematical fluid dynamics \cite{munch2005dewetting,craster2009dynamics,oron1997long,myers}
and PDE analysis  \cite{laugesen2000properties,laugesen2002energy,novick2010thin,taranets2014unstable,bertozzi1994lubrication}.
For $\beta=0$, \cite{bertozzi2001dewetting} gives a proof showing that given appropriate initial conditions, there is global existence of smooth positive solutions for all times. This proof can be extended to \eqref{evapModel} with $\beta>0$ \cite{ji2018instability,ji2017thin}.
However these proofs do not hold for $\beta<0$ and hence global existence is not guaranteed for that case.
\par
There have also been studies using other physically-based models of evaporating fluid films 
\cite{burelbach1988nonlinear, oron2001dynamics,oron1999dewetting, bestehorn2006regular,todorova2012relation,Thiele2014}
that have used forms of the evaporative flux that are not consistent with \eqref{generalModel}.
For instance, the linear evaporative flux used in \cite{bestehorn2006regular} can yield a pattern-forming instability.
In \cite{todorova2012relation} a base model for an evaporating film of the form \eqref{generalModel} was used, but the addition of an inhomogeneous term representing an influx of fluid from a porous substrate to maintain steady droplet solutions,  breaking the variational form.
Some PDE analysis has also addressed unstable thin film equations with linear or nonlinear absorption (or ``proliferation'') terms \cite{VSS,LP,Taranets2003,harwin2009}. For these equations the evolution of the energy does not follow \eqref{generalEdiss} and different approaches to the analysis of solutions must be used for each model \cite{engelnkemper2019continuation}. 
\par
In \cite{ji2017finite} it was shown that if the non-conserved term in \eqref{evapModel} is modified to have the opposite sign, then the PDE yields finite-time rupture singularities driven by the non-conserved loss term. Here we will more fully explore the other consequences that occur in a model that breaks from the gradient flow framework of \eqref{generalModel}.
\par 
In particular, we 
study the nonlinear partial differential equation in one dimension on a periodic domain $0 \le x \le L$,
\begin{subequations}\lbl{main}
\begin{align}
{\partial h\over \partial t} &= {\partial \over \partial x}\left( h^3
{\partial p \over \partial x}\right)+ \gamma p,
\lbl{mainpde}\\
\intertext{with}
p &= \Pi(h) - \frac{\partial^2 h}{\partial x^2},
\lbl{pressure} \\
\intertext{and}
\Pi(h) &= \frac{1}{h^3}- \frac{1}{h^4}-\p0
\lbl{pi},
\end{align}
\end{subequations}
with $\p0$ being a pressure-offset constant and $\gamma>0$ scaling the strength of the non-conserved effects.
\par
This system retains the mobility $\Mc = h^3$ and the generalized disjoining pressure $\Pi(h)$ from the volatile thin film equation \eqref{evapModel}, but uses a negative constant for the mobility $\Mnc = -\gamma$; apart from the sign-change, this is analogous to \eqref{evapModel} with $\K \gg h$.
For a thorough discussion of mass-conserving and mass non-conserving models in different contexts with variational or non-variational contributions, we refer readers to \cite{engelnkemper2019continuation}; our model with $\gamma > 0$ belongs to their mass non-conserving non-variational case.
\par
Note that \eqref{main} can also be written as
\begin{equation}
    {\partial h\over \partial t} = {\partial\over \partial x}\left(h^3
    {\partial\over \partial x} \left[\Pi(h)+ {\gamma\over 2h^2}-{\partial^2 h\over \partial x^2}\right]\right)
    +\gamma \Pi(h),
\end{equation}
namely as a Cahn-Hilliard equation with a modified potential and an added source term. It can be shown that to recast
\eqref{generalModel} in the form
$${\partial \phi\over \partial t}=\nabla\cdot(\Mc{}(\phi) \nabla \hat{\mu})+\hat{S}(\phi)\quad
\mbox{with $\hat{\mu}=\hat{U}'-\nabla^2 \phi$},$$
the mobility $\Mnc{}$ must be a constant. To define a gradient flow with energy dissipation \eqref{generalEdiss}, the source and modified potential must be related by $\hat{U}''=-\hat{S}'/\Mnc{}+\Mnc{}/\Mc{}$ for some $\Mnc{}>0$. This is a less convenient form than \eqref{generalModel}, but again serves to show that added lower order terms can change significantly the properties of the mass-conserving model. Rather than attempting to describe the influence of broad classes of source terms, we formulated \eqref{main} to take advantage of the conserved and non-conserved fluxes both being expressed simply in terms of a single pressure functional, but showing the consequences of not being a gradient flow.
\par
Starting from positive, finite-mass initial data $h_0(x)>0$, we will show that the interactions between the conserved and non-conserved terms in \eqref{main} yield novel behaviors that are not possible in \eqref{evapModel}. In particular we will see that there exist two distinct classes of nontrivial steady states, and that finite-time rupture and limit cycles can occur.
\par  
The structure of this paper is as follows. In Sec.~\ref{sec:properties} general properties and the steady state structure of the PDE \eqref{main} are outlined. In Sec.~\ref{sec:uniform2ndss} spatially uniform states and their instabilities with respect to spatial perturbations are investigated. Properties of nonuniform second-order and fourth-order states and the corresponding bifurcations will be discussed in Sec.~\ref{sec:2ndsteadystate} and Sec.~\ref{sec:4thss}. Stability of these states will be analyzed in Sec.~\ref{sec:stability}, followed by discussions of limit cycle dynamics in Sec.~\ref{sec:limitcycle} and the formation of finite time singularities in Sec.~\ref{sec:singularity}. Numerical simulations showing dynamic transitions among different steady states are presented in Sec.~\ref{sec:transition}. Concluding notes and a discussion of the remaining open questions are presented in Sec.~\ref{sec:conclusion}.

\end{section}

\section{Properties of PDE \eqref{main}}
\lbl{sec:properties}

The energy functional and its rate of change for equation \eqref{main} are given by
\begin{equation}
\EE = \int_0^L\frac{1}{2}
\left({\partial h\over \partial x}\right)^2 + U(h)\,dx,\qquad 
\frac{d\EE}{dt} = -\int_0^L h^3\left(\frac{\partial p}{\partial x}\right)^2\,dx +\gamma
\int_0^L {p^2} \,dx,
\lbl{dissipation}
\end{equation}
where 
$p=\delta \EE/\delta h$ is given by \eqref{pressure} and
the potential $U=\int \Pi(h) dh$ with $\Pi(h)$ defined in \eqref{pi}.
For $\gamma>0$, \eqref{dissipation}$_2$ shows that the contribution from the second integral can overcome the dissipation from the first integral, and hence the energy can evolve non-monotonically.
\par
The total mass of the solution and its rate of change are
\begin{equation}
\MM=\int_0^L h\,dx, \qquad {d\MM\over dt} = \gamma\int_0^L \Pi(h)\,dx.
\lbl{Dmass}
\end{equation}
This indicates that if $\Pi(h)\le 0$ (or $\Pi(h)\ge 0$) everywhere in the domain, then $\gamma >0$ yields a monotone decreasing (or increasing) mass over time. When $\Pi(h)$ changes sign, the evolution of the  mass depends on the form of the solution $h$ and the system parameter $\p0$, similar to the volatile thin film model \eqref{evapModel} studied in \cite{ji2018instability}.
\par
Now we focus on the positive steady states $H(x)>0$ of the model \eqref{main}.
Unlike most mass-conserving thin film equations, the interplay between the conserved and non-conserved terms in the equation \eqref{main} leads to more interesting sets of equilibria.
By setting the time derivative in \eqref{main} to be zero, we get the ODE system for all steady state solutions of the model as
\begin{equation}
\lbl{4thsteadystate}
0 = {d \over d x}\left( h^3
{d p \over d x}\right) +
{\gamma p}, \qquad
p = \Pi(h) - {d^2 h\over d x^2},
\end{equation}
which we study subject to periodic boundary conditions.
This system is equivalent to a fourth-order nonlinear ODE for $h$. If we set the pressure $p\equiv 0$ then \eqref{4thsteadystate} reduces to the second-order differential equation
\begin{equation}
\frac{d^2h}{dx^2}-\Pi(h)=0.
\lbl{2ndsteadystate}
\end{equation}
By setting the derivative term in \eqref{2ndsteadystate} to be zero, we can further reduce the second-order equation to an algebraic equation
\begin{equation}
\Pi(h)=0
\lbl{pi0}
\end{equation}
for spatially-uniform steady states. 
In particular, we will refer to these steady states as
\begin{itemize}
\item \textit{Uniform steady states}, $\bar{H}$, satisfying \eqref{pi0}.
\item \textit{Second-order nonuniform steady states}, $\HA(x)$, satisfying \eqref{2ndsteadystate} but not \eqref{pi0}.
\item \textit{Fourth-order nonuniform steady states}, $\HB(x)$, satisfying \eqref{4thsteadystate} but not \eqref{2ndsteadystate}.
\end{itemize}

There have been systematic studies of second-order steady states of both volatile and nonvolatile thin film equations \cite{perazzo2017analytical, bertozzi2001dewetting, thiele2010thin}. Investigations on fourth-order steady states are given in \cite{todorova2012relation, engelnkemper2019continuation}.
However, the present study is the first, to our knowledge, to classify equilibrium solutions of a thin film-type equation into both second- and fourth-order states.
In our framework, the second-order states $\HA$ and the uniform states $\bar{H}$ are associated with a constant zero pressure $p$, and the fourth-order states $\HB$ have nontrivial pressures.

\par 
The system parameter $\gamma$ is important for the existence and multiplicity of steady states of the model. For $\gamma \le 0$, only the uniform and second-order nonuniform steady states can exist \cite{ji2018instability}. That is, since the energy \eqref{dissipation} is a Lyapunov functional in this case, the steady states are given by extrema of the energy. Since the contributions from each integral in \eqref{dissipation}$_2$ are non-negative, they must equal to zero independently when the extrema of the energy is attained. Namely by setting ${d\EE}/{dt}=0$ in equation \eqref{dissipation}$_2$, from the first integral one obtains 
${\partial p}/{\partial x} = 0,$
which indicates that $p$ is a constant over the domain. Equivalently, for $\gamma=0$ the governing equation \eqref{main} reduces to a conserved equation, and the equilibria are given by solutions with a constant pressure, $p(x,t)\equiv P$. For $\gamma<0$ the second integral in \eqref{dissipation} then leads to $p=0.$
Therefore based on the form of the dynamic pressure $p$ in \eqref{pressure}, the equilibria satisfy the second-order nonlinear ODE \eqref{2ndsteadystate}. 
\par 
For $\gamma> 0$, we will show that the fourth-order equilibria satisfying the fourth-order ODE \eqref{4thsteadystate} co-exist with the second-order steady states of \eqref{main} in sufficiently large domains. This is a key  difference between the steady state structure of our problem and that of classic conserved thin film models. The second-order steady states are solutions of the mass-conserving model and are independent of $\gamma$ while the fourth-order steady states depend on balances of the conserved and non-conserved fluxes.
\par 
We start with a brief discussion of the spatially uniform and second-order steady states. Since these solutions satisfy the second-order ODE \eqref{2ndsteadystate}, 
their properties are similar to those of the steady states of the model \eqref{mainpde} for $\gamma < 0$ that have been studied in \cite{ji2018instability}.

\section{Spatially uniform steady states and their stability}
\lbl{sec:uniform2ndss}

\begin{figure}
\centering
\includegraphics[width=6.3cm, height=4.2cm]{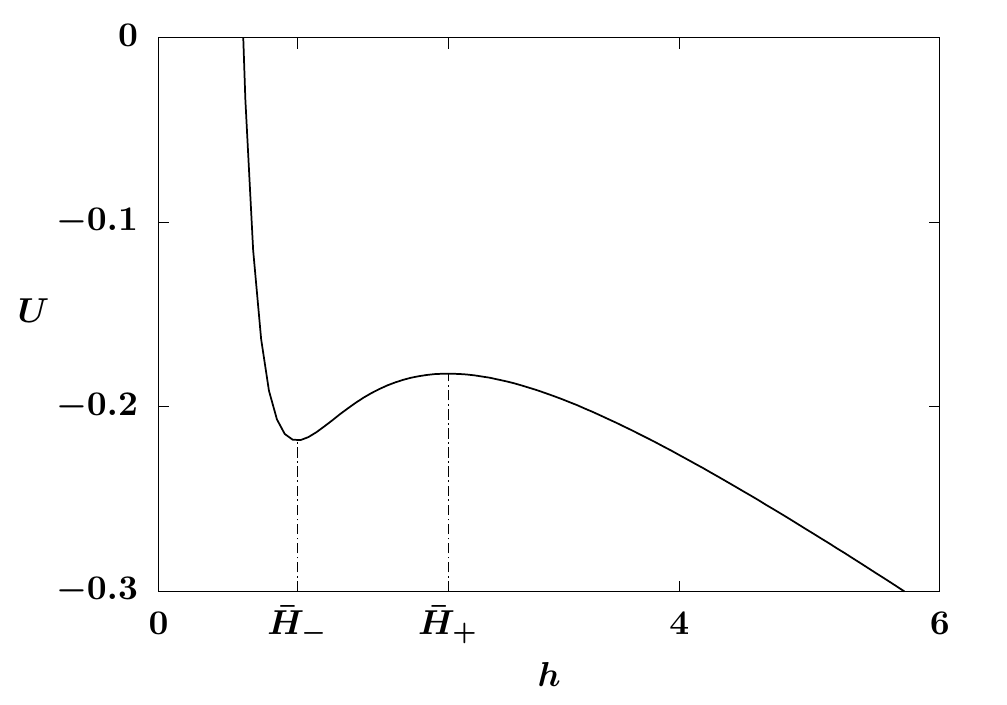}
\includegraphics[width=6.3cm]{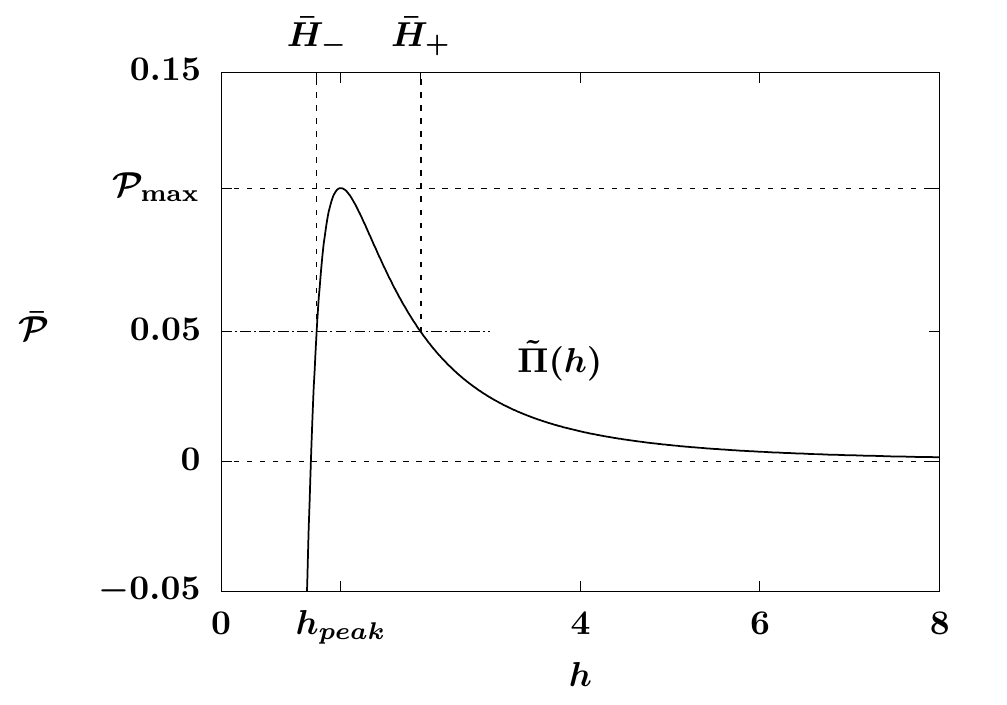}
\caption{(Left) The potential $U(h)$ in \eqref{Potential} with  $\p0 = 0.05$. (Right) The base disjoining pressure $\tilde{\Pi}(h) = h^{-3}-h^{-4}$  with the dashed lines giving the range $0\le \p0 < \peps$ where multiple uniform steady states exist, and $\Hm \approx 1.06, \Hc\approx 2.22$ for $\p0=0.05$.} 
\lbl{fig:uniformsteadystate}
\end{figure}
\par 
The spatially uniform steady states of equation \eqref{main} are determined by the critical points of the potential $U(h)$ used in \eqref{dissipation}.
 For clarity, we explicitly write out the potential $U(h)$ as the integral of $\Pi(h)$ from \eqref{pi},
\begin{equation}
U(h) = \int \Pi(h)\,dh = -\frac{1}{2h^2}+\frac{1}{3h^3}-\p0h.
\lbl{Potential}
\end{equation}
Fig.~\ref{fig:uniformsteadystate}~(left) shows a plot of $U(h)$ with $\p0 = 0.05$ which has a local minimum at $\Hm$ and a local maximum at $\Hc$, both of which satisfy \eqref{pi0}.
\par 
This potential $U(h)$ has a unique inflection point at $h_{\text{peak}} =4/3$ (where $\Pi(h)$ has its global maximum), and on the semi-infinite range $h>\hpeak$, $U''(h)<0$ (see Fig.~\ref{fig:uniformsteadystate}~(left)). This can be regarded as a degenerate case of the bi-stable double well potential used in Cahn-Hilliard models, which have a finite range of $h$ satisfying $U''(h)<0$, called the \textit{spinodal} range.
\par
The number of spatially uniform steady states of equation \eqref{main} depends on the value of $\p0$.
Fig.~\ref{fig:uniformsteadystate}~(right) depicts a plot of $\tilde{\Pi}(h) = h^{-3}-h^{-4}$  which has its unique maximum at $\hpeak$, yielding
\begin{equation}
\Pi(h_{\text{peak}}) = \peps-\p0, \quad \peps = {\ts\frac{27}{256}}>0.
\end{equation}
If $\p0<0$, there is one spatially uniform steady state $\Hm < h_{\text{peak}}$. 
For 
$0<\p0 < \peps$, there are two uniform states $\Hm< h_{\text{peak}}< \Hc$; we will call this the critical range of the pressures $\p0$. As $\p0$ approaches $\peps$, the two states coincide at the double root $h = h_{\text{peak}}$. 
For $\p0> \peps$, there are no uniform steady states.
\par 
These ranges of $\p0$ are important for the analysis of the steady states and the dynamics of the model \eqref{main}.
For $\p0\le \peps$, the uniform steady state $\Hm$ for $\p0 \to 0$ is given by
\begin{equation}
\Hm = 1+\p0+4\p0^2+O(\p0^{3}).
\lbl{hmin*}\\
\end{equation}
For $0<\p0 \le \peps$ the other uniform steady state $\Hc$ can be written as
\begin{align}
\Hc =& \p0^{-1/3}-{\ts \frac{1}{3}-\frac{2}{9}\p0^{1/3}-\frac{20}{81}}\p0^{2/3}+O(\p0).
\lbl{hc*}
\end{align}
While $\Hc$ has a leading order dependence on $\p0$, $\Hm$ has a weaker dependence on $\p0$ with a saddle-node bifurcation occurring at $\p0 = \peps$. 
\par 
For $0<\p0 \le \peps$, the coexisting uniform steady states $\Hc$ and $\Hm$ give rise to novel and interesting solution structures and dynamics. Therefore,  for the rest of this paper we focus on the critical range.
\par 
Now we consider the linear stability of a uniform steady state $\bar{H}$ with respect to an infinitesimal Fourier mode disturbance, $h(x,t) = \bar{H} + \delta e^{i2k\pi x/L}e^{\lambda t}$,
where $k=0,1,2,3,\cdots$ and $\lambda$ describes the growth rate for perturbations starting from the initial amplitude $\delta \ll 1$. Expanding \eqref{main} about the uniform steady state $h = \bar{H}$ then gives the $O(\delta)$ equation
\begin{equation}
\lambda = 
\left[-M_c(\bar{H})\left(\frac{2k\pi}{L}\right)^2+\gamma\right]
\left[\Pi'(\bar{H})+\left(\frac{2k\pi}{L}\right)^2\right],
\lbl{dispersion}
\end{equation}
where the first factor corresponds to an operator that includes the mobility functions from both conserved and non-conserved parts of \eqref{mainpde}, and the second factor comes from the linearized pressure operator set by \eqref{pressure} \cite{ji2018instability}. 
The dispersion relation \eqref{dispersion} shows that the stability of the uniform steady states depends on both the parameter $\gamma$ and the domain size. 
% In the long-wave limit, $L\to \infty$, equation \eqref{dispersion} reduces to $\lambda = \gamma \Pi'(\bar{H})$ for any fixed $k$, and
% since $\Pi'(\Hm) > 0$ and $\Pi'(\Hc) < 0$, $\Hm$ is unstable and $\Hc$ is stable.

Fig.~\ref{fig:dispersion} gives a plot of the dispersion relation \eqref{dispersion}. For $\bar{H} = \Hm$, we have $\Pi'(\Hm) > 0$ and the state $\Hm$ is long-wave unstable with respect to perturbations with $0 \le 2k\pi/L<[\gamma/M_c(\Hm)]^{1/2}$. For $\bar{H} = \Hc$, since $\Pi'(\Hc) < 0$, the state $\Hc$ is stable in the long-wave limit $L\to \infty$. The form of the curve depends on the parameter $\gamma$, and for small $\gamma$ (as in Fig.~\ref{fig:dispersion}) the state $\Hc$ is unstable for $[\gamma/M_c(\Hc)]^{1/2}<2k\pi/L < [-\Pi'(\Hc)]^{1/2}$.
For $k=0$, \eqref{dispersion} reduces to $\lambda=\gamma \Pi'(\bar{H})$ corresponding to the spatially-uniform 
rate of change of mass per unit length in \eqref{Dmass}.
We note that while the dispersion relation $\lambda(k)$ can be long-wave unstable, the influence of the fourth-order term always yields strong damping for large $k$. For small $k$ the form of the curves are similar to the type-II and type-III curves in \cite{cross2009pattern}. These behaviors are comparable with dispersion relations for other volatile thin film models with gradient and non-gradient dynamics forms \cite{thiele2010thin,bestehorn2006regular}.

\begin{figure}
\centering
\includegraphics[width=8cm]{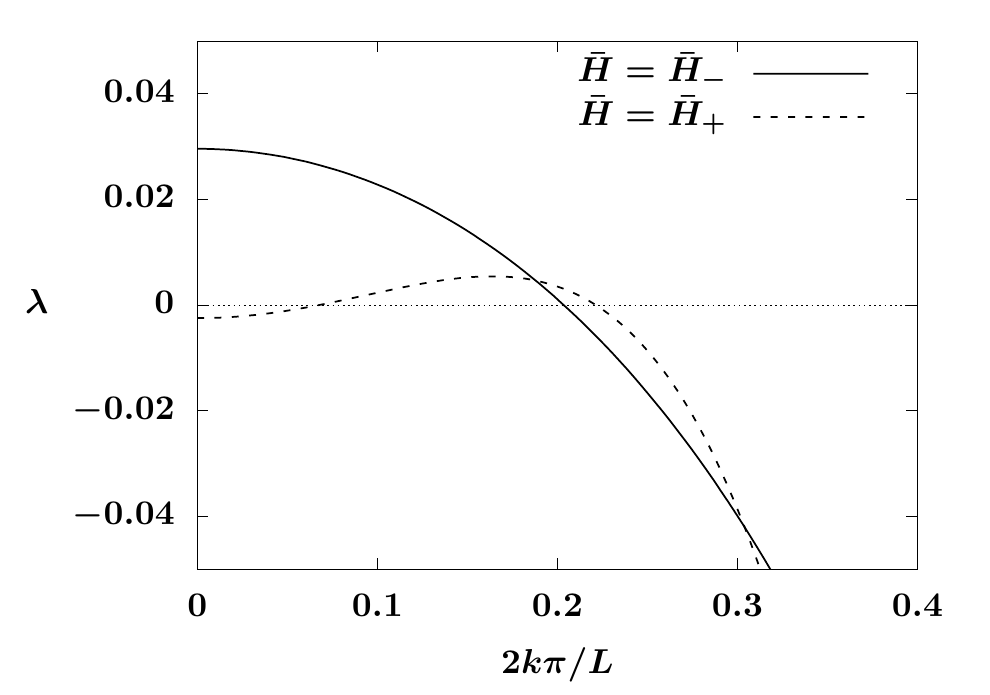}
\caption{A plot of the dispersion relation \eqref{dispersion} for $\bar{H} = \Hm$ and $\bar{H}=\Hc$ with $\p0 = \gamma = 0.05$.} 
\lbl{fig:dispersion}
\end{figure}

\par 
The growth rate $\lambda$ changes sign at roots of the two factors in \eqref{dispersion}. In terms of the domain size $L$, these roots correspond to critical lengths at which the uniform state $\bar{H}$ changes stability. 
In particular, for $M_c(H) = H^3$ the first factor gives roots satisfying $\bar{H}^3(2k\pi/L)^2=\gamma$.
We express these roots as $L=k\bar{\ell}^\gamma$ for $k=1,2,3,\cdots$ in terms of primary critical lengths, given by 
\begin{equation}
\lBc = \frac{2\pi \Hc^{3/2}}{\gamma^{1/2}}\qquad \mbox{for $\bar{H}=\Hc$,}
\lbl{LBC}
\end{equation}
and
\begin{equation}
\lB = \frac{2\pi \Hm^{3/2}}{\gamma^{1/2}} \qquad \mbox{for $\bar{H}=\Hm$.}
\lbl{LB}
\end{equation}
Both $\lB$ and $\lBc$ strongly depend on the parameter $\gamma$ and involve the interaction of both the conserved and non-conserved fluxes.
\par 
The roots of the second factor in \eqref{dispersion} satisfy $(2k\pi/L)^2 = -\Pi'(\bar{H})$. Since $\Pi'(\Hm) > 0$, the $\Hm$ state yields no real roots. But for $\bar{H}=\Hc$, $\Pi'(\Hc) < 0$ and we obtain another family of critical lengths $L = k\lA$ based on the linearized pressure, where the primary critical length is
\begin{equation}
\lA = \frac{2\pi}{\sqrt{-\Pi'(\Hc)}}.
\lbl{LA}
\end{equation}

These critical lengths are also useful in characterizing the dependence of the stability of the uniform states $\Hm$ and $\Hc$ on the domain size $L$. In the limit $L \to 0$, both states are linearly stable to spatial perturbations. The state $\Hm$ is stable for $0<L < \lB$ and unstable for $L > \lB$, while the stability of the state $\Hc$ depends on the values of $\lA$ and $\lBc$.

\par 
At these critical lengths, $k\lB$, $k\lBc$, and k$\lA$, nonuniform steady state solutions bifurcate from the uniform states. These will be the branches of second-order and fourth-order steady states. Moreover, these critical lengths give the only possible bifurcation points from the uniform states. Other apparent crossings between branches of nonuniform and uniform states at points not in this set are due to the projection of the full bifurcation structure onto the two-dimensional diagrams that we will be using, showing the average of the solution vs.\ the domain length. For an example of this in the next section, see Fig.~\ref{fig:2ndsteadystate}~(left).

\section{Second-order nonuniform steady states}
\lbl{sec:2ndsteadystate}

Here we briefly describe the nonuniform second-order steady states.
They are like steady states of broad families of conserved thin film equations that have been extensively studied in \cite{laugesen2000properties,laugesen2002energy} and elsewhere.
The second-order ODE \eqref{2ndsteadystate} can be written as a phase plane system
\begin{equation}
{dh\over dx} =s \qquad {ds\over dx}= \Pi(h).
\lbl{2ndss}
\end{equation}
Equilibrium points of this system correspond to the spatially uniform states $\Hm$, $\Hc$ satisfying \eqref{pi0}.  Linearizing the system around these equilibrium points using $h = \bar{H}+\delta e^{\sigma x}$ yields $\sigma = \pm \sqrt{\Pi'(\Hm)}$ for $\bar{H}=\Hm$, and $\sigma = \pm \imagi \sqrt{-\Pi'(\Hc)}$ for $\bar{H}=\Hc$. Hence $h=\Hm$ is a hyperbolic saddle point, and $h = \Hc$ is a center point. There is a homoclinic orbit through $\Hm$ with $\Hm\le h\le H_{\max}$ where $U(H_{\max})=U(\Hm)$ and 
\begin{equation}
H_{\max}=(6\p0)^{-1} + 1 + O(\p0)\qquad \mbox{for $\p0\to 0$}
\lbl{HmaxEqn}
\end{equation}
that represents the maximal size of a single droplet state on the whole real axis. In the phase plane this orbit encloses the continuous family of periodic solutions centered around $\Hc$.
\par 
For any fixed value of $\p0$ in the critical range $0<\p0 \le \peps$, selecting a minimum height $\hmin$ in the range $\Hm \le \hmin \le \Hc$ determines a periodic steady state solution $\HA(x)$ which has a maximum height $h_{\max}$ that satisfies $U(h_{\max}) = U(\hmin)$. 
Since we can write the first integral of \eqref{2ndss} as
$s^2 = 2(U(h)-U(\hmin))$,
the period $\ell(h)$ of the steady state is given by
\begin{equation}
\ell(h_{\min}) = \int_0^{\ell(h_{\min})}~dx = 2\int_{h_{\min}}^{h_{\max}}\frac{1}{\sqrt{2U(h)-2U(h_{\min})}}\,dh.
\end{equation}
In the limit case when $\hmin = \Hm$, we have the solitary droplet solution with the length $\ell \to \infty$.

\begin{figure}
\centering
\includegraphics[width=6.6cm]{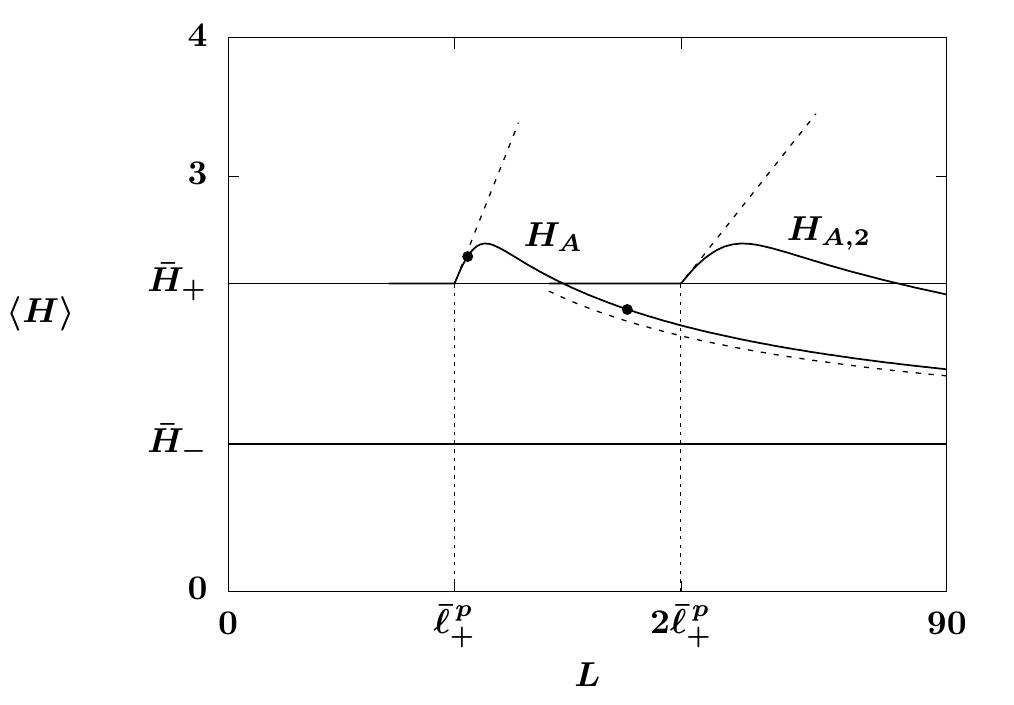}
\includegraphics[width=6.6cm]{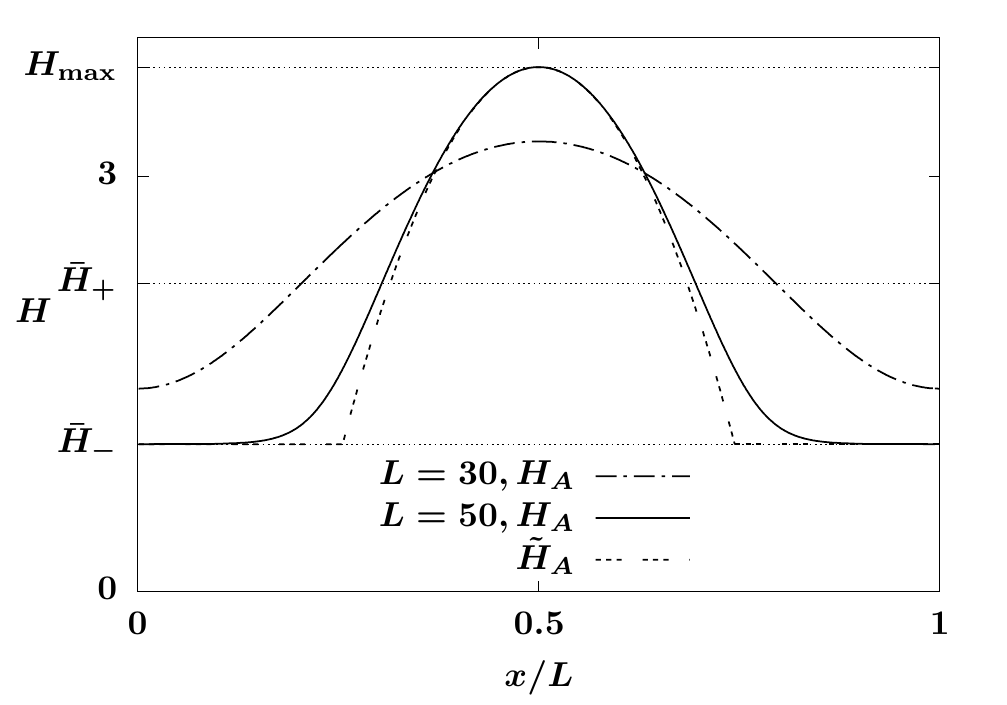}
\caption{(Left) Bifurcation diagram for $\left\langle H \right\rangle$ for spatially-uniform and second-order states $\HA$ parametrized by the domain size $L$. The dashed curve below $\left\langle H_A \right\rangle$ and the dashed lines tangent to $\left\langle H_{A,k} \right\rangle$ at $L = k\lA$ correspond to local estimates described by \eqref{LALinf} and \eqref{bifurcationLAk}, respectively. These results are valid for any $\gamma \neq 0$ and are plotted here for $\p0 = 0.05$ which yields $\lA\approx 28.4$. (Right) Plots of rescaled states $\HA$ corresponding to the two marked states (dots) in Fig.~\ref{fig:2ndsteadystate}~(left)
satisfying \eqref{2ndsteadystate}, and a rescaled plot of the truncated parabolic profile $\tilde{H}_A$ defined in \eqref{eq:2ndApprox}.}
\lbl{fig:2ndsteadystate}
\end{figure}
\par 
Small amplitude second-order steady states bifurcate from the spatially-uniform steady state $\Hc$. Note that the second-order states are only related to the pressure \eqref{pressure} and do not depend on the mobility, and their minimum period of oscillations is given by $L = \lA$ in \eqref{LA}.
For $L \ge \lA$, we define $h = \HA(x)$ as the principal second-order steady state with a single maximum on the domain $0 \le x \le L$. With the average of the solution defined by 
\begin{equation}
\langle h \rangle = \frac{1}{L}\int_0^L h(x)\,dx,
\lbl{hAve}
\end{equation}
Fig.~\ref{fig:2ndsteadystate}~(left) shows the average of the solution, $\langle \HA\rangle$,
for $\p0 = 0.05$ parametrized by $L$ bifurcating supercritically from $\Hc$ at $L = \lA$. Corresponding to the dots marked in Fig.~\ref{fig:2ndsteadystate}~(left), typical profiles of steady states $H_{A}(x)$ with periods $L=30$ and $L = 50$ are presented in Fig.~\ref{fig:2ndsteadystate}~(right) using rescaled spatial variables. 
A bifurcation diagram similar to Fig.~\ref{fig:2ndsteadystate}~(left) for steady states of a thin volatile film model can be found in \cite{thiele2010thin}.
\par 
To get an estimate of $\langle \HA\rangle$ for the limit $L \to \infty$, we use the result that
 if $h_{\max}\gg \Hm$ the core structure of $\HA(x)$ can be approximated by a truncated parabolic profile \cite{bertozzi2001dewetting,glasner2003coarsening}
centered at $x=L/2$,
\begin{equation}
H_A(x) \ge \tilde{H}_A(x) =\max\left(h_{\max}+\frac{\Pi(h_{\max})}{2}\left(x-{\textstyle {1\over 2}} L\right)^2, \Hm \right).
\lbl{eq:2ndApprox}
\end{equation}
Here $h_{\max}$ is the maximum of $\HA$ attained at $x = L/2$ and satisfies $\Pi(h_{\max}) = d^2{H}_A/dx^2|_{x=L/2}$ from \eqref{2ndsteadystate}. 
For $L \to \infty$, $H_A$ approaches the homoclinic solution, and we can use $h_{\max}\approx H_{\max}$. 
As a consequence, a lower bound for the $\langle \HA\rangle$ branch as $L \to 
\infty$ is given by
\begin{equation}
   \langle \HA\rangle > \frac{1}{L}\int_0^L \tilde{H}_A(x)\, dx   = \Hm + 
      \left(\frac{4\sqrt{2} (H_{\max}-\Hm)^{3/2}}{3 \sqrt{|\Pi(H_{\max})|}}\right) {1\over L}.
      \lbl{LALinf}
\end{equation}
For $\p0 = 0.05$, this corresponds to the approximation $\langle \HA\rangle \approx 1.06 + 44.4/L$, which is plotted in Fig.~\ref{fig:2ndsteadystate}~(left) in a dashed curve beneath the $\langle \HA\rangle$ branch. A plot of the truncated parabolic profile $\tilde{H}_A$ is also included in Fig.~\ref{fig:2ndsteadystate}~(right) which shows a good agreement with the core structure of $H_A(x)$ for $L = 50$.
\par
For the domain size $L \ge 2\lA$, multiple periodic nonuniform second-order steady states can coexist in the system, as families of second-order solutions $H_{A,k}(x)$ with $k$ periods bifurcate from $\Hc$ at $L = k\lA$ for $k = 2, 3, \cdots$
We will use $\p0 = 0.05$ as a typical value for the parameter $\p0$ for the rest of the paper.

We will focus on the primary bifurcation point $L = \lA$.
Using the rescaling $x = LX$, $\HA(x) = H(X)$ and perturbing the domain size $L$ around the bifurcation point $\lA$ by a small parameter $\epsilon\ll 1$, $L = \lA+\epsilon$, we rewrite the ODE \eqref{2ndsteadystate} as
\begin{equation}
\frac{d^2 H}{d X^2} - (\lA+\epsilon)^2 \Pi(H) = 0, \qquad 0 \le X \le 1.
\lbl{2ndODErescale}
\end{equation}
Following the local bifurcation analysis in \cite{witelski2000dynamics,bertozzi2001dewetting}, we expand the solution in the neighborhood of $\lA$ as
\begin{equation}
H(X) = \Hc + \delta H_1(X)+\delta^2 H_2(X) + \delta^3 H_3(X) + O(\delta^4),
\lbl{hAexpansion}
\end{equation}
where $\delta \ll 1$ is a small perturbation parameter, and we need to determine the relation between $\delta$ and $\epsilon$.
Substituting \eqref{hAexpansion} in \eqref{2ndODErescale} yields a harmonic oscillator equation at $O(\delta)$, 
\begin{equation}
\hat{\mathcal{L}}H_1=0\qquad\mbox{where}\quad
\hat{\mathcal{L}}H_1 \equiv \frac{d^2 H_1}{d X^2}+\left(- {\left(\lA\right)}^2 \Pi'(\Hc)\right) H_1,
\end{equation}
which for $\lA = {2\pi}/\sqrt{-\Pi'(\Hc)}$ from \eqref{LA} has the solution 
$H_1(X) = \A\cos(2\pi X),$ where $\A$ is the amplitude.
The $O(\delta^2)$ equation
\begin{equation}
\hat{\mathcal{L}}H_2 = \frac{1}{2}{\left(\lA\right)}^2 \Pi''(\Hc) H_1^2
\lbl{HA2nd}
\end{equation}
gives the solution
$$
H_2(X) = -\frac{\A^2\Pi''(\Hc)}{12\Pi'(\Hc)}\left(3-\cos(4\pi X)\right).
$$
In obtaining this form, we have assumed that $\epsilon \ll O(\delta)$ in order to exclude a resonant term from \eqref{HA2nd} that would have made it impossible to find a periodic solution.
\par
Assuming that $\epsilon = O(\delta^2)$, at $O(\delta^3)$ the equation for $H_3(X)$ is $\hat{\mathcal{L}}H_3 = \mathcal{R}_3(\Hc, H_1, H_2)$, where $\mathcal{R}_3$ is the right-hand side inhomogeneous term. From the solvability condition for the $O(\delta^3)$ equation, $\int_0^1 \mathcal{R}_3H_1(x)~dx = 0$, and the condition that the amplitude $\A$ of the $O(\delta)$ perturbation is real, we determine that the parameter $\epsilon$ is positive, $\epsilon >0$. 
\par 
This indicates that a supercritical bifurcation occurs at $L = \lA$, and the amplitude $\A$ of the leading-order perturbation is
\begin{equation}
\A^2 = \left|\frac{24\left[-\Pi'(\Hc)\right]^{5/2}}{\displaystyle{5\pi \Pi''(\Hc)^2-3\pi \Pi'''(\Hc)\Pi'(\Hc)}}\right|.
\end{equation}
The local structure of $\langle \HA \rangle$ near the bifurcation point $L = \lA$ is then given by
\begin{equation}
\langle \HA \rangle = \int_0^1 \HA(X)\,dX \sim \Hc + \delta^2\int_0^1 H_2(X)\,dX  = \Hc  -\frac{\A^2\Pi''(\Hc)}{4\Pi'(\Hc)}(L-\lA).
\lbl{bifurcationLA}
\end{equation}
Equivalently, the approximation $\HA(X) \sim \Hc + \A\cos(2\pi X)(L - \lA)^{1/2}$ holds in the limit $L \searrow \lA$; 
plotting the extrema of the solution vs.\ $L$ would then show the local structure as a supercritical pitchfork bifurcation.
Large-amplitude solutions of the strongly nonlinear equation \eqref{2ndODErescale} can also be constructed by matched asymptotic expansions \cite{bertozzi2001dewetting}.
\par 
Without loss of generality, the linear approximation \eqref{bifurcationLA} can be extended to
\begin{equation} 
\langle H_{A,k} \rangle \sim  \Hc  -\frac{\A^2\Pi''(\Hc)}{4k\Pi'(\Hc)}(L-k\lA)\qquad \mbox{for $L\to k\lA$}
\lbl{bifurcationLAk}
\end{equation}
for the branches of second-order states $H_{A,k}$ with multiple periods bifurcating from $\Hc$.
These approximations in the neighborhood of bifurcation points $L = \lA$ and $L = 2\lA$ are plotted against numerical solutions in Fig.~\ref{fig:2ndsteadystate}~(left). 
This bifurcation diagram will be revisited later for the regime where fourth-order nonuniform steady states also exist in the system.
\par 
The slopes of these linear approximations are positive since $\Pi''(\Hc) >0$ and $\Pi'(\Hc) < 0$ for $\p0 = 0.05$. 
The disjoining pressure has an inflection point at $h_{\mathrm{I}}=5/3$ where $\tilde{\Pi}''(h_{\mathrm{I}})=0$, with $\mathcal{P}_{\mathrm{I}}=54/625=0.0864$.
For $\p0$ in the range $\mathcal{P}_{\mathrm{I}}<\p0 < \peps$, the uniform state $\Hc$ satisfies $\Pi''(\Hc)<0$ (see Fig.~\ref{fig:uniformsteadystate}~(right)). This choice of $\p0$ will lead to 
the approximation \eqref{bifurcationLA} with a negative slope and the average $\langle H_{A,k} \rangle \le \Hc$ for all second-order nonuniform solutions.

The results about the local structure of the second-order steady states are valid for any $\gamma \neq 0$. However, the stability of these states depend on the value of $\gamma$, which will be studied in Sec.~\ref{sec:stability}.
As a reminder, the appearance of each $H_{A,k}$ branch crossing $\Hc$ at some point to the right of $k\lA$ for $\p0$ in the range $0<\p0 <\mathcal{P}_{\mathrm{I}}$ is a consequence of the fact that Fig.~\ref{fig:2ndsteadystate}~(left) is a two-dimensional projection of the full structure, as described in Sec.~\ref{sec:uniform2ndss}.
These branches must pass by each other in parameter space since there are no bifurcation points of $\Hc$ at these values of $L$.

\section{Fourth-order nonuniform steady states}
\lbl{sec:4thss}

Next we investigate the family of fourth-order nonuniform steady states that satisfy the fourth-order ODE \eqref{4thsteadystate} but not the second-order ODE \eqref{2ndsteadystate}. 
Equation \eqref{4thsteadystate} can be written as a fourth-order autonomous system
\begin{equation}
{dh\over dx} =s \qquad {ds\over dx}= \Pi(h)-p \qquad {dp\over dx}= {q\over h^3}
\qquad {dq\over dx} = - \gamma p
\lbl{4thss}
\end{equation}
with equilibrium points of the system corresponding to spatially uniform states, with
$\Pi(h)=0$ (and $s=p=q=0$). 
As in the second-order case \eqref{2ndss} these uniform states are $h=\Hm$ and $h=\Hc$. But now, linearization of the system around these equilibrium points with $h = \bar{H}+\delta e^{\sigma x}$ leads to
$\sigma = \pm\sqrt{\Pi'(\Hm)}$, $\pm \imagi\sqrt{\gamma/\Hm^3}$ for $\bar{H} = \Hm$, and 
$\sigma = \pm\imagi\sqrt{-\Pi'(\Hc)}$, $\pm \imagi\sqrt{\gamma/\Hc^3}$ for $\bar{H} = \Hc$.
Therefore $h = \Hm$ is a saddle-focus and $h = \Hc$ is a focus-focus type equilibrium point \cite{kuznetsov}.
\par
We note that following our choice of terminology from Section~\ref{sec:properties}, for the fourth-order steady states, we will consider only solutions with nontrivial pressures, thus excluding the second-order and uniform steady states.

\begin{figure}
\centering
\includegraphics[width=6.6cm]{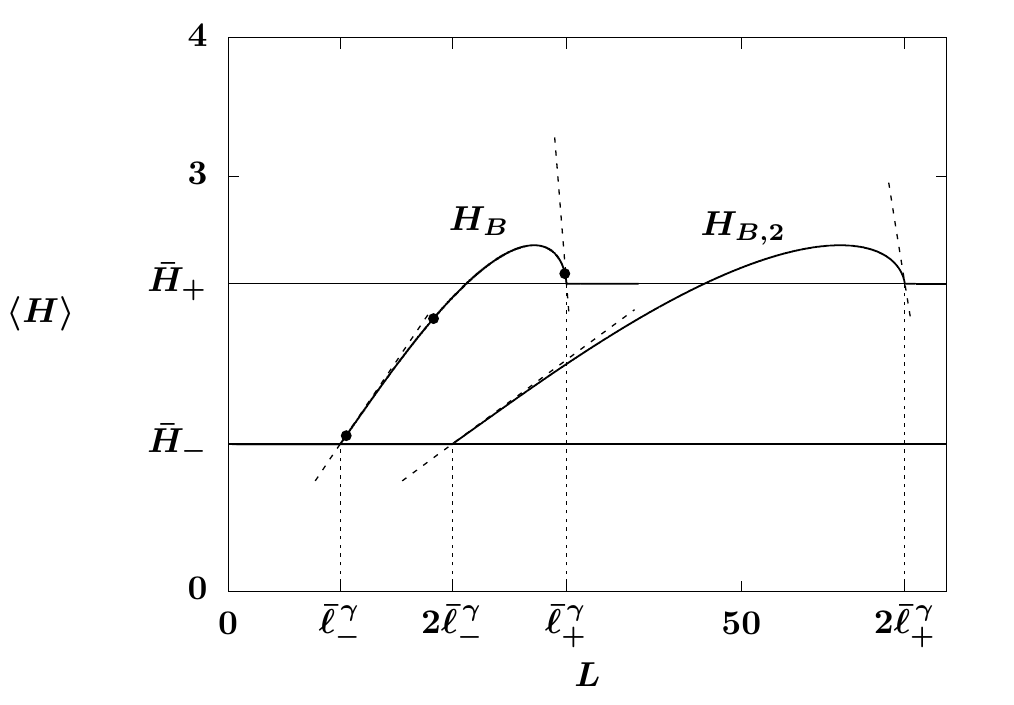}
\includegraphics[width=6.6cm]{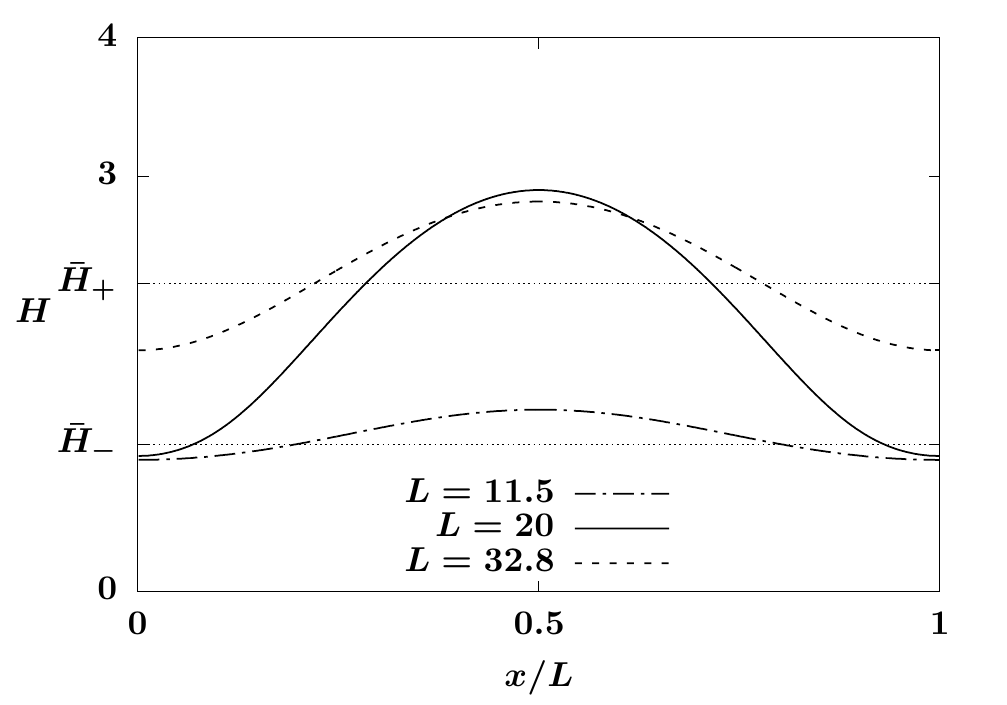}
\caption{(Left) Bifurcation diagrams of $\left\langle H \right\rangle$ for fourth-order steady states $H_{B,k}(x)$ parametrized by the domain size $L$, compared against predictions \eqref{4thbif_lb} and \eqref{4thbif_lbc} for the local bifurcation structure (dashed lines). (Right) Plots of typical rescaled states $\HB(x)$ corresponding to the marked dots in Fig.~\ref{fig:4thss}~(left). The system parameters are given by $\gamma = 0.4$, $\p0 = 0.05$, yielding $\lB \approx 10.91$ and $\lBc \approx 32.96.$}
\lbl{fig:4thss}
\end{figure}

\subsection{Primary bifurcations from uniform states}
\lbl{sec:primaryBif}
Small amplitude fourth-order nonuniform steady states bifurcate from both spatially-uniform states $\Hm$ and $\Hc$, and their minimum periods of oscillations are given by $L = \lB$ in \eqref{LB} and $L=\lBc$ in \eqref{LBC}, respectively.  
Similar to the second-order case, we denote $H_B(x)$ as the principal branch of fourth-order states with period $L$, and denote the family of fourth-order states with $k$ periods as $H_{B,k}(x)$, $k=2, 3, \cdots$,
\par 
For $\gamma=0.4$, Fig.~\ref{fig:4thss}~(left) depicts the first two branches of fourth-order states $H_{B}(x)$ and $H_{B,2}(x)$ in a bifurcation diagram parametrized by $\left\langle H \right\rangle$. These solutions start from the spatially-uniform state $\Hm$ at bifurcation points $L = k\lB$, and end at the uniform state $\Hc$ at $L = k \lBc$.
Typical plots of these states on the principal branch $\HB$ with periods $L = 11.5$, $L= 20$, and $L = 32.8$ are shown in Fig.~\ref{fig:4thss}~(right) which correspond to the marked dots in Fig.~\ref{fig:4thss}~(left). Note that the $\HA$ nonuniform states  are not included in this figure.

To determine the local bifurcation structure near a critical period $L = \ell_c$, we apply the approach
used in Sec.~\ref{sec:2ndsteadystate} again. Applying the rescaling $x \to LX$ to \eqref{4thsteadystate} with the perturbation $L = \ell_c + \epsilon$ for $\epsilon \ll 1$ leads to the rescaled fourth-order ODE on $0 \le X \le 1$
\begin{equation}
\frac{d}{dX}\left[M_c(H)\frac{d}{dX}\left((\ell_c+\epsilon)^2\Pi(H)-\frac{d^2 H}{dX^2}\right)\right] + \gamma(\ell_c+\epsilon)^2\left[(\ell_c+\epsilon)^2\Pi(H)-\frac{d^2 H}{dX^2}\right] = 0.
\lbl{4thODErescale}
\end{equation}
In the neighborhood of $\ell_c$ we expand the steady state $H_B(X)$ as in \eqref{hAexpansion}, $H_B(X) = \bar{H} + \delta H_1(X)+\delta^2 H_2(X) + \delta^3 H_3(X) + O(\delta^4)$ where $\delta \ll 1$. The $O(\delta)$ equation of the expansion of \eqref{4thODErescale} leads to the critical period $\ell_c = 2\pi\bar{H}^{3/2}/\sqrt{\gamma}$ as in \eqref{LBC} and \eqref{LB}. \par 
Similar to the construction of the second-order solutions $\HA(x)$, the form of the fourth-order solutions for $\delta \to 0$ is given by
\begin{equation}
H_B(X) = \bar{H} +\delta \A \cos(2\pi X) + O(\delta^2),
\end{equation}
where $\A$ is the amplitude of the leading-order perturbation to $\bar{H}$. We omit the details of calculating the $H_2(X)$ solution and move directly to the Fredholm solvability condition for the $O(\delta^3)$ equation which determines that 
\begin{equation}
    \A^2 = \frac{\epsilon}{\delta^2}\mathcal{W},
\end{equation}
where
\begin{equation}
\mathcal{W} = -\,{\frac {24\sqrt {{\gamma}}{{ M_c(\bar{H})}}^{3/2}\, \left( {\Pi'(\bar{H})}\,{ M_c(\bar{H})}+4\,{\gamma}
 \right){\Pi'(\bar{H})} }{ \pi \left(  \mathcal{W}_1 {{
\Pi'(\bar{H})}}^{2}+ \mathcal{W}_2 {\Pi'(\bar{H})}+24\,{ M_c(\bar{H})}\,{ M_c'(\bar{H})}\,{ \Pi''(\bar{H})\,{ \gamma}}
 \right)}},
\end{equation}
where
\begin{equation}
    \mathcal{W}_1 = -3\,{{
M_c(\bar{H})}}^{2}{M_c''(\bar{H})}+4\,{M_c(\bar{H})}\,{{M_c'(\bar{H})}}^{2},\nonumber
\end{equation}
\begin{equation}
    \mathcal{W}_2 = 3\,{{M_c(\bar{H})}}^{2}{M_c'(\bar{H})}\,{ \Pi''(\bar{H})}-12\,{M_c(\bar{H})}\,{ M_c''(\bar{H})}\,{ \gamma}+28\,{{ M_c'(\bar{H})}}^{2}{\gamma}.\nonumber
\end{equation}
Since $\A^2 \ge 0$, the sign of $\mathcal{W}$ determines the relation between the small parameters $\epsilon$ and $\delta$. In particular, for $M_c(\bar{H}) = \bar{H}^3$ and $\bar{H} = \Hm$, we have $\mathcal{W} > 0$ for all $\gamma > 0$, which yields $\epsilon = \delta^2$. 
Therefore, the average of the solution $\langle \HB\rangle$ near the supercritical bifurcation point $\ell_c = \lB$ is approximated by 
\begin{equation}
\langle \HB\rangle \sim \Hm - \frac{\A^2\Pi''(\Hm)}{4\Pi'(\Hm)}(L-\lB) \quad \text{for}\, L \to \lB \text{ with } L \ge \lB.
\lbl{4thbif_lb}
\end{equation}
For $\bar{H} = \Hc$, there are two cases depending on the value of the parameter $\gamma$ relative to functions of $\Hc$:
\begin{itemize}
\item 
For values of $\gamma$ yielding $\mathcal{W}<0$ then we take $\eps= -\delta\le 0$ to indicate a subcritical bifurcation with
\begin{subequations}
\begin{equation}
\langle \HB\rangle \sim \Hc - \frac{\A^2\Pi''(\Hc)}{4\Pi'(\Hc)}(\lBc-L) \quad \text{for}\, L \to \lBc \text{ with } L \le \lBc,
\lbl{4thbif_lbc}
\end{equation}
\item For values of $\gamma$ yielding $\mathcal{W}>0$ then $\eps=\delta\ge 0$ for a supercritical bifurcation with
\begin{equation}
\langle \HB\rangle \sim \Hc - \frac{\A^2\Pi''(\Hc)}{4\Pi'(\Hc)}(L-\lBc) \quad \text{for}\, L \to \lBc \text{ with } L \ge \lBc.
\lbl{4thbif_lbc_small_gamma}
\end{equation}
\end{subequations}
\end{itemize}
For $\gamma = 0.4$ with $\mathcal{W}<0$, a comparison between $\langle \HB\rangle$ and their linear approximations near these bifurcation points are shown in Fig.~\ref{fig:4thss}~(left). 
We also note that the critical value $\gamma = -\Pi'(\Hc)\Hc^3/4$ is a root of the equation $\lBc=2\lA$ (and a zero of $\mathcal{W}(\gamma)$) where $\lA$ and $\lBc$ are defined in \eqref{LA} and \eqref{LBC}.  Therefore we expect the local structure of the $\HB$ branch to change qualitatively as the critical period $\lBc$ exceeds $2\lA$.

\begin{figure}
\centering
\includegraphics[width=6.5cm]{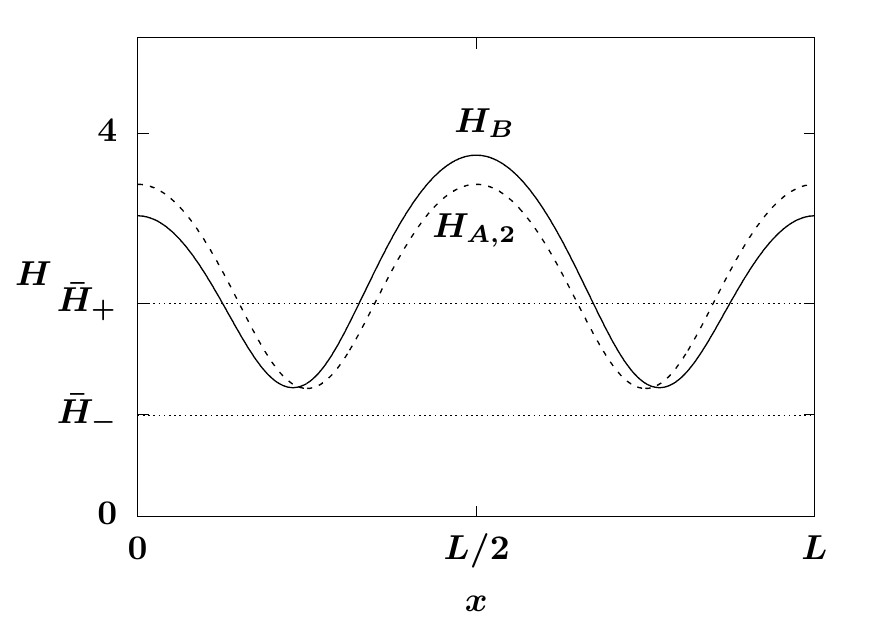}
\includegraphics[width=6.5cm]{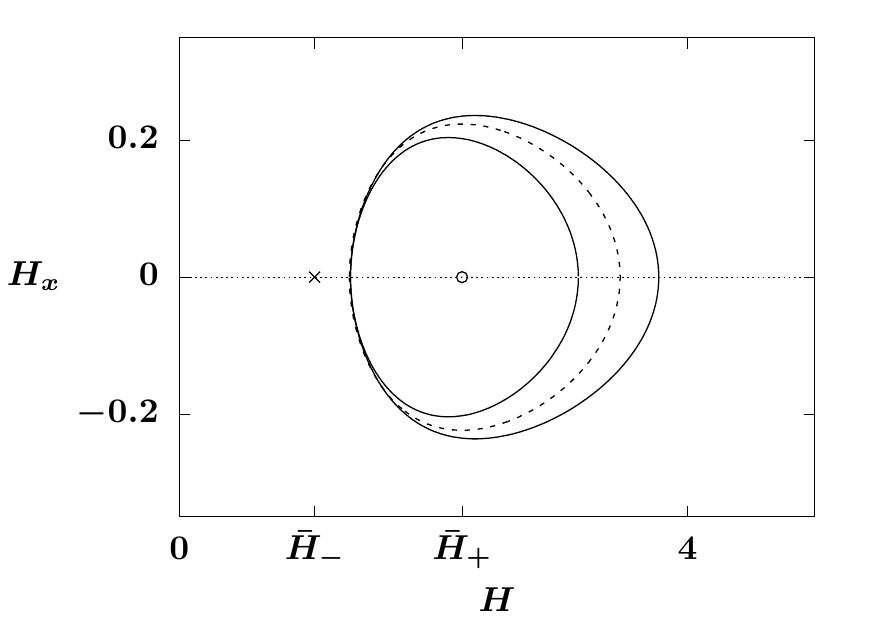}
\caption{Profiles and $(H,H_x)$ phase portrait of the fourth-order steady state $\HB(x)$ and the second-order steady state $H_{A,2}(x)$ showing that $\HB(x)$ bifurcates from $H_{A,2}(x)$ at $L = 61.7$ (near  $\ell_{AB}^{\,(2)}$), with $\gamma= 0.05$ and $\p0=0.05$.}
\lbl{fig:4thss_p0l0}
\end{figure}

\subsection{Secondary bifurcations from nonuniform second-order states}
\lbl{sec:secondaryBif}
The family of fourth-order states $H_{B,k}(x)$ may also undergo secondary bifurcations from the nonuniform second-order states $H_{A,k}(x)$. As an example, Fig.~\ref{fig:4thss_p0l0} shows the profiles and the $(H, H_x)$ phase portraits of the second-order and fourth-order steady states, $H_{A,2}(x)$ and $H_{B}(x)$, that coexist at $L = 61.7$.  The profiles of the two steady states are close to each other,
 with the main peak of $H_{B}(x)$ being higher than that of $H_{A,2}(x)$, and the secondary peak of $H_{B}(x)$ is lower than that of $H_{A,2}(x)$, as in a period-doubling bifurcation.
In this section, we use Floquet theory \cite{iakubovich1975linear,teschl2012ordinary,chicone2006ordinary} to identify these secondary bifurcations.

We consider a fixed domain size $L$ near a critical period $\ell_c$
at which a fourth-order state $\HB$ bifurcates from a branch of second-order states $\HA$. From \eqref{4thsteadystate}, the fourth-order state $\HB$ corresponds to a non-trivial pressure function $P_B(x) = \Pi(\HB)-d^2\HB/dx^2$. In contrast, a second-order state $\HA$ corresponds to a trivial pressure $p\equiv 0$. 
We expand the fourth-order state and its corresponding pressure $(\HB,P_B)$ around the principal second-order state $\HA(x)$ and its corresponding zero pressure using
\begin{equation}
    H_B(x) = \HA(x)+\delta \tilde{H}(x)+O(\delta^2), \qquad
    P_B(x) = \delta\tilde{P}(x) + O(\delta^2), \nonumber
\end{equation}
where the small parameter $\delta \ll 1$ is set by the distance $|L-\ell_c|$.
% , and the exponent $\alpha > 0$.
Substituting these expansions into \eqref{4thsteadystate} and using \eqref{2ndsteadystate} we obtain the $O(\delta)$ linearized problem 
\begin{equation}
\frac{d}{dx}\left(\HA^3 \frac{d\tilde{P}}{dx}\right)+\gamma \tilde{P}=0,\qquad
\tilde{P} = \Pi'(\HA)\tilde{H}-\frac{d^2\tilde{H}}{dx^2}.
\lbl{floquet1}
\end{equation}
A secondary bifurcation of $\HB$ from $\HA$ will be identified if \eqref{floquet1} has a non-trivial $L$-periodic solution $(\tilde{H},\tilde{P})$, $\tilde{H}\not\equiv 0$ and $\tilde{P}\not\equiv 0$.
% Here we have assumed that the exponent $\alpha = 1$ to obtain the balanced $\tilde{H}$ terms and the $\tilde{P}$ term in the second equation in \eqref{floquet1}, otherwise \eqref{floquet1} will lead to a translational mode $\tilde{H} = d\HA/dx$ with $\tilde{P} \equiv 0$ which does not produce a non-trivial $P_B(x)$. 
\par 
For a given $L$-periodic $H_A$ state, the problem \eqref{floquet1} is a fourth-order differential equation with periodic coefficients set by $\HA$. 
% It is useful to rewrite \eqref{floquet1} as an autonomous system
% \begin{equation}
% \frac{d\tilde{H}}{dx} = W, \quad 
% \frac{d\tilde{W}}{dx} = \Pi'(\HA)\tilde{H} - \tilde{P},\quad
% \frac{d\tilde{P}}{dx} =\frac{\tilde{Q}}{\HA^3},\quad
% \frac{d\tilde{Q}}{dx} = -\gamma\tilde{P}.
% \lbl{eqn:floquet4}
% \end{equation}
While a periodic solution to \eqref{floquet1} needs to satisfy the full system, we can get necessary (but not sufficient) conditions for the existence of such periodic solutions by first looking at the $2\times 2$ system for $\mathbf{Y}(x) = [\tilde{P}(x),\tilde{Q}(x)]^{T}$ from the first equation in \eqref{floquet1},
\begin{equation}
\frac{d\mathbf{Y}}{dx} = \mathbf{A}_2(x)\mathbf{Y},\qquad 
\mathbf{A}_2(x) = 
\begin{bmatrix}
    0       & \HA(x)^{-3} \\
    -\gamma   &    0 \\
\end{bmatrix}.
\lbl{floquet3}
\end{equation}
We numerically solve this system as an initial value problem, and obtain the principal fundamental matrix solution $\mathbf{\Phi}(x)$ which satisfies
\begin{equation}
\frac{d\mathbf{\Phi}}{dx} = \mathbf{A}_2(x)\mathbf{\Phi}, \qquad \mathbf{\Phi}(0) = \mathbf{I}_2.
\end{equation}
Here we use a predetermined $\HA$ with its peak located at $x = L/2$ (see Fig.~\ref{fig:2ndsteadystate}~(right)), and $\mathbf{I}_2$ denotes the $2\times 2$ identity matrix. 
Evaluating the matrix solution $\mathbf{\Phi}(x)$ at $x = L$, we obtain the monodromy matrix $\mathbf{B}_2 = \mathbf{\Phi}(L)$ whose eigenvalues are characteristic multipliers of  equation \eqref{floquet3} \cite[Section 3.6]{teschl2012ordinary}. 
Based on Floquet theory \cite[Chapter II section 1]{iakubovich1975linear}, if the matrix $\mathbf{B}_2$ has an eigenvalue $\rho$ that satisfies $\rho^k = 1$ for a positive integer $k$, then there exists a $kL$-periodic solution $\mathbf{Y}$ of ODE \eqref{floquet3}. 
In particular, when $\mathbf{B}$ has an eigenvalue $\rho = 1$, it corresponds to an $L$-periodic nonuniform solution $\mathbf{Y}$ of \eqref{floquet3}.

\begin{figure}
\centering
\includegraphics[width=6cm,height = 4.5cm]{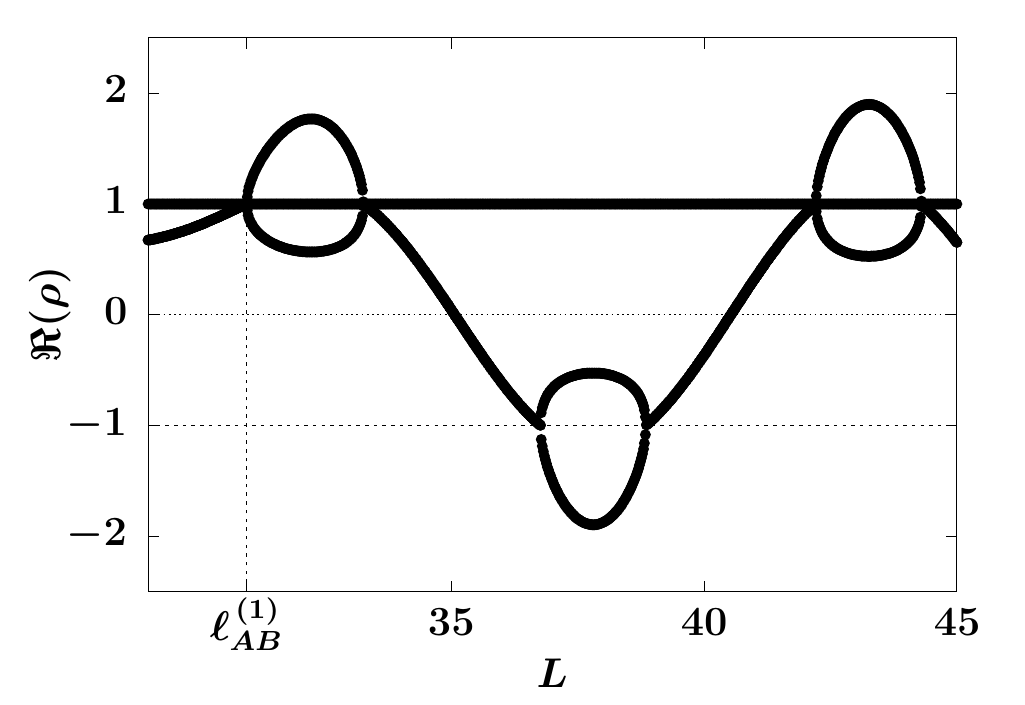}
\includegraphics[width=7.2cm,height = 4.5cm]{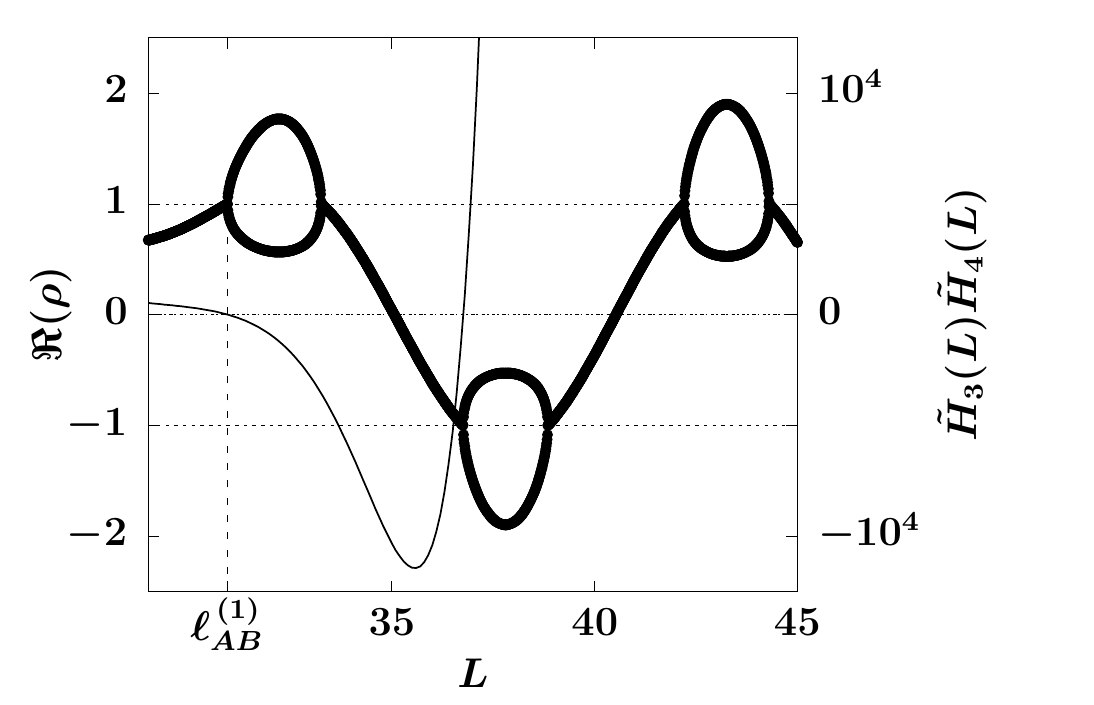}
 \caption{The real part of the characteristic multipliers $\rho$ of (left) the fourth-order system \eqref{floquet4} and (right) the second-order ODE system  \eqref{floquet3} plotted against the domain size $L$ for $\gamma = 0.4$. The smallest critical period $L = \ell_{AB}^{\,(1)}\approx 30.95$ associated with $\rho = 1$ is a secondary bifurcation of $\HB$ from $\HA$ since this is also a zero of $\tilde{H}_3(L)\tilde{H}_4(L)$ (thin curve).}
%  and (Right) Perturbation solutions $\tilde{P}$ associated with $\ell_c = \ell_{AB}^{\,(k)}$ for $k = 1,2,3$.
\lbl{fig:floquet}
\end{figure}

Fig.~\ref{fig:floquet}~(right) shows the dependence of characteristic multipliers $\rho$ on the domain size $L$ for $\gamma = 0.4$.
The ``egg-shaped'' portions of the plot correspond to pairs of real multipliers, and the remaining portions represent complex conjugates pairs of multipliers $\rho, \bar{\rho}$ satisfying $|\rho| = 1$. As $L$ increases the structure of the set of multiplier values appears to converge to a repeating pattern. Our numerical study also shows that at the critical domain sizes corresponding to a characteristic multiplier $\rho = 1$, the monodromy matrix $\mathbf{B}_2$ takes the form
\begin{equation}
\mathbf{B}_2 = 
\begin{bmatrix}
    1       & 0 \\
    a   &    1 \\
\end{bmatrix}\qquad \text{ or } \qquad
\mathbf{B}_2 = 
\begin{bmatrix}
    1       & b \\
    0   &    1 
\end{bmatrix},
\lbl{eqn:B2}
\end{equation}
where the constants $a, b \neq 0$, indicating that the multiplier $\rho = 1$ is a double eigenvalue of geometric multiplicity one.

While the Floquet analysis for the second-order system \eqref{floquet3} identifies multiple critical domain sizes $L$ with $\rho = 1$, we will show that only the first critical $L$ corresponds to a secondary bifurcation of $\HA$. To exclude the other spurious critical $L$, we apply Floquet analysis to the full system in \eqref{floquet1}, and consider the fourth-order differential equation system for $\mathbf{Z}(x) = [\tilde{H}(x), \tilde{W}(x), \tilde{P}(x), \tilde{Q}(x)]^T$,
\begin{equation}
\frac{d\mathbf{Z}}{dx} = \mathbf{A}_4(x)\mathbf{Z},\qquad
\mathbf{A}_4(x) = 
\begin{bmatrix}
 0 & 1 & 0 & 0 \\
 \Pi'(\HA) & 0 & -1 & 0\\
 0 & 0 & 0 & \HA(x)^{-3} \\
 0 & 0 &  -\gamma   &  0
\end{bmatrix}.
\lbl{floquet4}
\end{equation}
Again the principal fundamental solution matrix $\mathbf{\Psi}$ of the system \eqref{floquet4} satisfies
\begin{equation}
\frac{d\mathbf{\Psi}}{dx} = \mathbf{A}_4(x)\mathbf{\Psi}, \qquad \mathbf{\Psi}(0) = \mathbf{I}_4,
\end{equation}
where the solution matrix
$\mathbf{\Psi} = \left[\psi_1, \psi_2, \psi_3,\psi_4\right]$, 
with the $i$-th column vector 
$\psi_i(x) = [\tilde{H}_i(x), \tilde{W}_i(x), \tilde{P}_i(x), \tilde{Q}_i(x)]^T$, 
and
$\mathbf{I}_4$ denotes the $4\times 4$ identity matrix.
Since the differential equations for $\tilde{P}$ and $\tilde{Q}$ do not depend on $\tilde{H}$ or $\tilde{W}$, 
the monodromy matrix $\mathbf{B}_4 = \mathbf{\Psi}(L)$ of equation \eqref{floquet4} can be written as a block matrix of four $2\times 2$ matrices,
\begin{equation}
\mathbf{B}_4 = 
\begin{bmatrix}
    \mathbf{C}_2       &  \mathbf{D}_2\\
    \mathbf{0}   &    \mathbf{B}_2 \\
\end{bmatrix},
\end{equation}
where
\begin{equation}
\mathbf{C}_2 = 
\begin{bmatrix}
    1       &  0\\
    c   &    1
\end{bmatrix},\quad
\mathbf{D}_2 = 
\begin{bmatrix}
    \tilde{H}_3(L)       &  \tilde{H}_4(L) \\
    \tilde{W}_3(L)   &    \tilde{W}_4(L) 
\end{bmatrix},\quad
\mathbf{B}_2 = 
\begin{bmatrix}
    \tilde{P}_3(L)       &  \tilde{P}_4(L) \\
    \tilde{Q}_3(L)   &    \tilde{Q}_4(L) 
\end{bmatrix},
\lbl{eqn:B4}
\end{equation}
where the constant $c\neq 0$.
At the critical domain sizes identified in Fig.~\ref{fig:floquet}~(right), the matrix $\mathbf{B}_2$ is given by \eqref{eqn:B2}.
We numerically calculate the eigenvalues of the monodromy matrix $\mathbf{B}_4$ and plot the corresponding characteristic multipliers $\rho$ in Fig.~\ref{fig:floquet}~(left).
In addition to the characteristic multipliers that have appeared in Fig.~\ref{fig:floquet}~(right) for the second-order system \eqref{floquet3}, there exists another double multiplier $\rho = 1$ of geometric multiplicity $1$ for all $L$. This additional $\rho = 1$ corresponds to the translational mode $\tilde{H} = d\HA/dx$ that lies on the branch of second-order states $\HA$ and does not yield any non-trivial $\tilde{P}$ solutions. 
\par
Since ${\psi}_1(x)$ and ${\psi}_2(x)$ only contain trivial $\tilde{P}$ solutions, in order to construct a periodic solution to \eqref{floquet4} with a non-trivial $\tilde{P}$, one needs to use a linear combination of ${\psi}_3(x)$ and ${\psi}_4(x)$ which satisfy ${\psi}_3(0) = [0,0,1,0]^T$ and ${\psi}_4(0) = [0,0,0,1]^T$, $\mathbf{Z}(x) = c_3{\psi}_3(x) + c_4{\psi}_4(x)$. The periodicity of $\mathbf{Z}$ requires $\mathbf{Z}(0) = \mathbf{Z}(L)$, which leads to
\begin{subequations}
\lbl{eqn:floquetcombination}
\begin{equation}
    c_3\tilde{H}_3(L) + c_4\tilde{H}_4(L) = 0,
\end{equation}
\begin{equation}
    c_3\tilde{P}_3(L) + c_4\tilde{P}_4(L) = c_3,
\end{equation}
\begin{equation}
    c_3\tilde{Q}_3(L) + c_4\tilde{Q}_4(L) = c_4.
\end{equation}
\end{subequations}
Note that at the critical periods, $\mathbf{B}_2$ takes the form \eqref{eqn:B2}, therefore we have $\tilde{P}_3(L) = \tilde{Q}_4(L) = 1$, and \eqref{eqn:floquetcombination} leads to $c_4\tilde{P}_4(L) = c_3\tilde{Q}_3(L) = 0$. Moreover, since $(\tilde{P}_4(L),\tilde{Q}_3(L)) = (0, a) $ or $(b,0)$ with $a,b\neq 0$, the coefficients $c_3$ and $c_4$ satisfy $c_3 = 0$ or $c_4 = 0$. This indicates that the periodic solution $\tilde{H}$ is given by a multiple of $\tilde{H}_3(x)$ or $\tilde{H_4}(x)$. 

Therefore, it suffices to check whether $\tilde{H_3}(L)  = 0$ or $\tilde{H_4}(L) = 0$ to determine if a candidate critical period of the system \eqref{floquet3} indeed corresponds to a periodic solution $\tilde{H}$ of \eqref{floquet4}  with a non-trivial $\tilde{P}$.
The plot of $\tilde{H}_3(L)\tilde{H}_4(L)$ as a function of $L$ in Fig.~\ref{fig:floquet}~(right) shows that only the smallest critical period $L = \ell_{AB}^{\,(1)}$ associated with $\rho=1$ corresponds to such a periodic solution.
This domain size yields a symmetry-preserving perturbation solution $\tilde{H}$ with respect to the reflectional symmetry of the second-order state $H_A(x)$ about $x = L/2$, and corresponds to a secondary bifurcation from the principal second-order state $H_A$ to the principal fourth-order state $H_B$. 

For other characteristic multipliers satisfying $\rho^k = 1$ with higher values of $k$, more perturbation solutions $(\tilde{H},\tilde{P})$ with period $kL$ exist in the system. They can lead to more interesting symmetry-preserving and symmetry-breaking secondary bifurcations, similar to the bifurcation of $H_B$ from $H_{A,2}$ observed in Fig.~\ref{fig:4thss_p0l0}. We will not attempt an exhaustive investigation of these bifurcations in this paper.

\subsection{Overall structure of the solution branches}
\begin{figure}
\centering
\includegraphics[width=6.6cm]{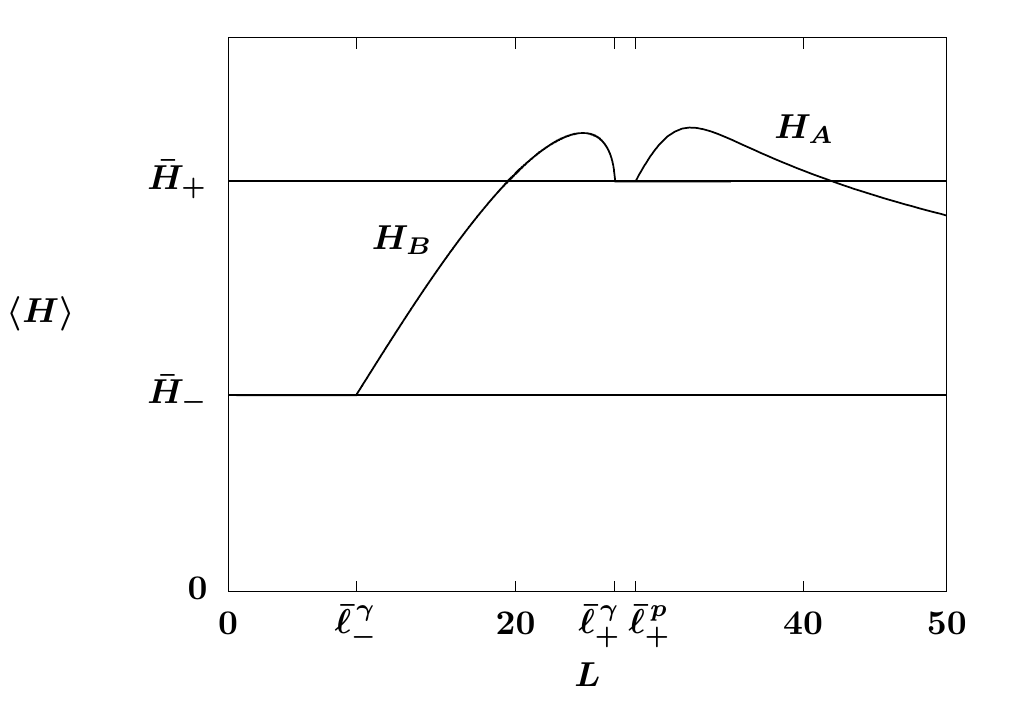}
\includegraphics[width=6.6cm]{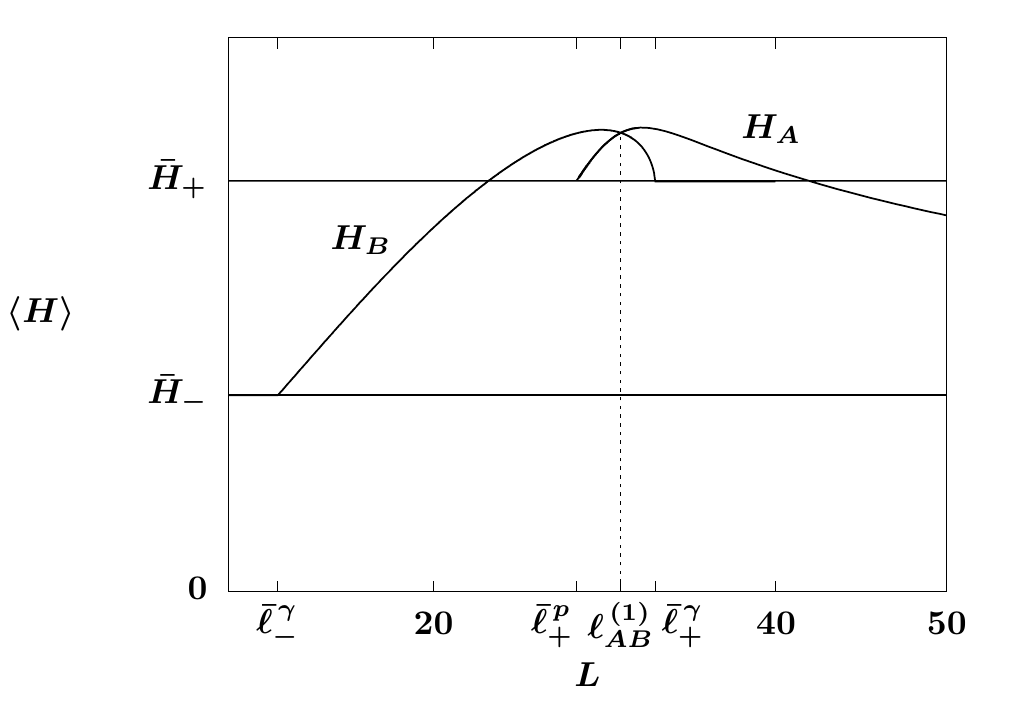}
\caption{Bifurcation diagrams for $\left\langle H \right\rangle$ for coexisting principal steady states $\HA$ and $\HB$ parametrized by the domain size $L$ with (left) $\gamma = 0.6$, $\lBc < \lA$ and  (right) $\gamma = 0.4$, $\lA  < \lBc < 2\lA$.}
\lbl{fig:4thss_bifurcation}
\end{figure}

So far we have worked out some local properties of the $H_{A,k}$ and $H_{B,k}$ branches. We will now consider some overall aspects of how these branches interact with each other.
\par 
Recall that 
the $H_{A,k}$ branches always start from a supercritical bifurcation at $L = k\lA$ from the uniform state $\Hc$ and have $\langle H_A\rangle$ asymptotically approaching $\Hm$  as $L\to \infty$ (see Fig.~\ref{fig:2ndsteadystate}~(left)).
These bifurcation points $L = k\lA$ are set by $\p0$ and are independent of $\gamma$. In contrast, the bifurcation points $L = k\bar{\ell}^{\,\gamma}_{\pm}$ at which the $H_{B,k}$ branches bifurcate from $\bar{H}_{\pm}$ depend strongly on $\gamma$. Therefore, different configurations of the $H_{A,k}$ and $H_{B,k}$ branches can be expected for different values of $\gamma$.

For simplicity we focus on the principal branch of fourth-order states $\HB(x)$ (we will briefly mention the generalization for $H_{B,k}$ with $k=2,3,\cdots$ later)
and discuss the following two basic cases:
\begin{itemize}
    \item  \underline{Case I}: Primary-Primary branches -- branches starting and ending from primary bifurcations with constant states $\bar{H}_{\pm}$ at $\bar{\ell}^{\,\gamma}_{\pm}$.
    \item \underline{Case II}: Primary-Secondary branches -- branches starting at a primary bifurcation from $\bar{\ell}^{\,\gamma}$, and ending at a secondary bifurcation with a $H_{A,k}$ branch.
\end{itemize}

\subsubsection{\underline{Case I}: Primary-Primary branches}
We start with the simplest case of an $\HB$ branch that starts and ends from primary bifurcations, and does not intersect with any $\HA$ branches. To ensure that there are no intersections, $\gamma$ must satisfy $\lBc < \lA$. Having $\gamma$ in this range was the first case discussed in Sec.~\ref{sec:primaryBif}. 
An example of such an $\HB$ branch is shown in Fig.~\ref{fig:4thss_bifurcation}~(left) with $\gamma = 0.6$, where the $H_B$ branch starts from a supercritical bifurcation at $L=\lB$, and ends at a subcritical bifurcation at $L=\lBc$. The $\HB(x)$ solutions on this branch all have a single local maxima.
\par 
A primary-primary $\HB$ branch may also undergo a secondary bifurcation as a transcritical crossing with the principal branch $\HA$.
This regime corresponds to the value of $\gamma$ satisfying $\lA < \lBc < 2\lA$, so that $\HB$ intersects the principal second-order branch $\HA$, but not other second-order branches. 
%%%%%
Using the local properties of the $\HB$ branch discussed in section \ref{sec:primaryBif}, having 
$\gamma$ satisfy the upper bound $\lBc = 2\lA$ 
corresponds to the primary bifurcation at $L = \lBc$ changing from a subcritical bifurcation to a supercritical when as $\lBc$ exceeds $2\lA$.
%%%%%%
Fig.~\ref{fig:4thss_bifurcation}~(right) shows an example of this case with $\gamma = 0.4$, where the $\HB$ branch starts from $\Hm$, ends at $\Hc$, and crosses the $\HA$ branch at a transcritical secondary bifurcation point $L = \ell_{AB}^{\,(1)}$. This secondary bifurcation at $L = \ell_{AB}^{\,(1)}$ has been identified in section \ref{sec:secondaryBif}.
The fourth-order states on this $\HB$ branch also have a unique maximum (see Fig.~\ref{fig:4thss}~(right)) and coincide with the second-order state at $L = \ell_{AB}^{\,(1)}$.

\subsubsection{\underline{Case II}: Primary-Secondary Branches}

\begin{figure}
\centering
\hbox{\hspace{0.2in}\includegraphics[width=12.3cm,height = 6cm]{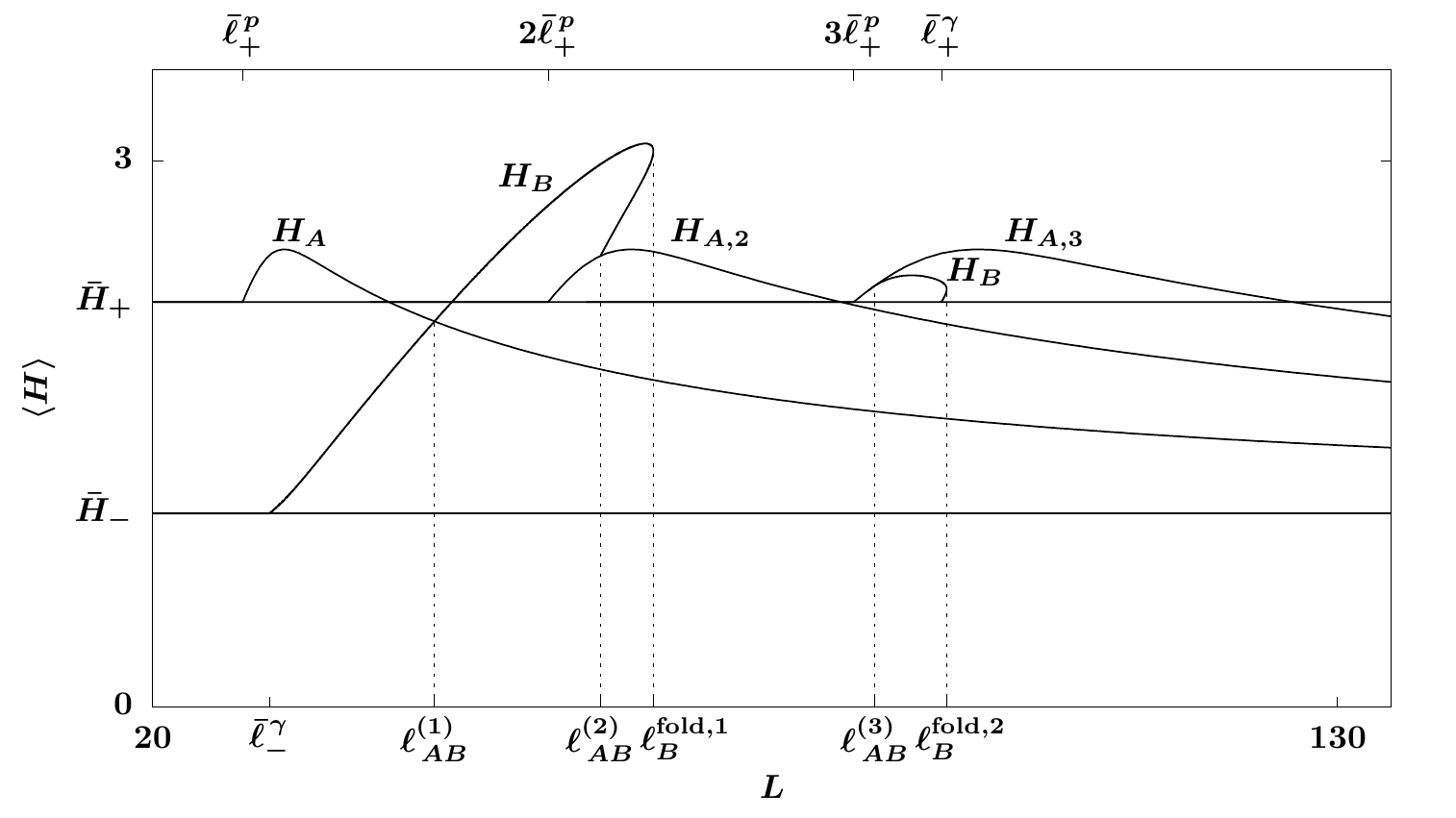}}
\vspace{-0.1in}
\hbox{\hspace{-0.05in}\includegraphics[width=12.9cm,height = 6cm]{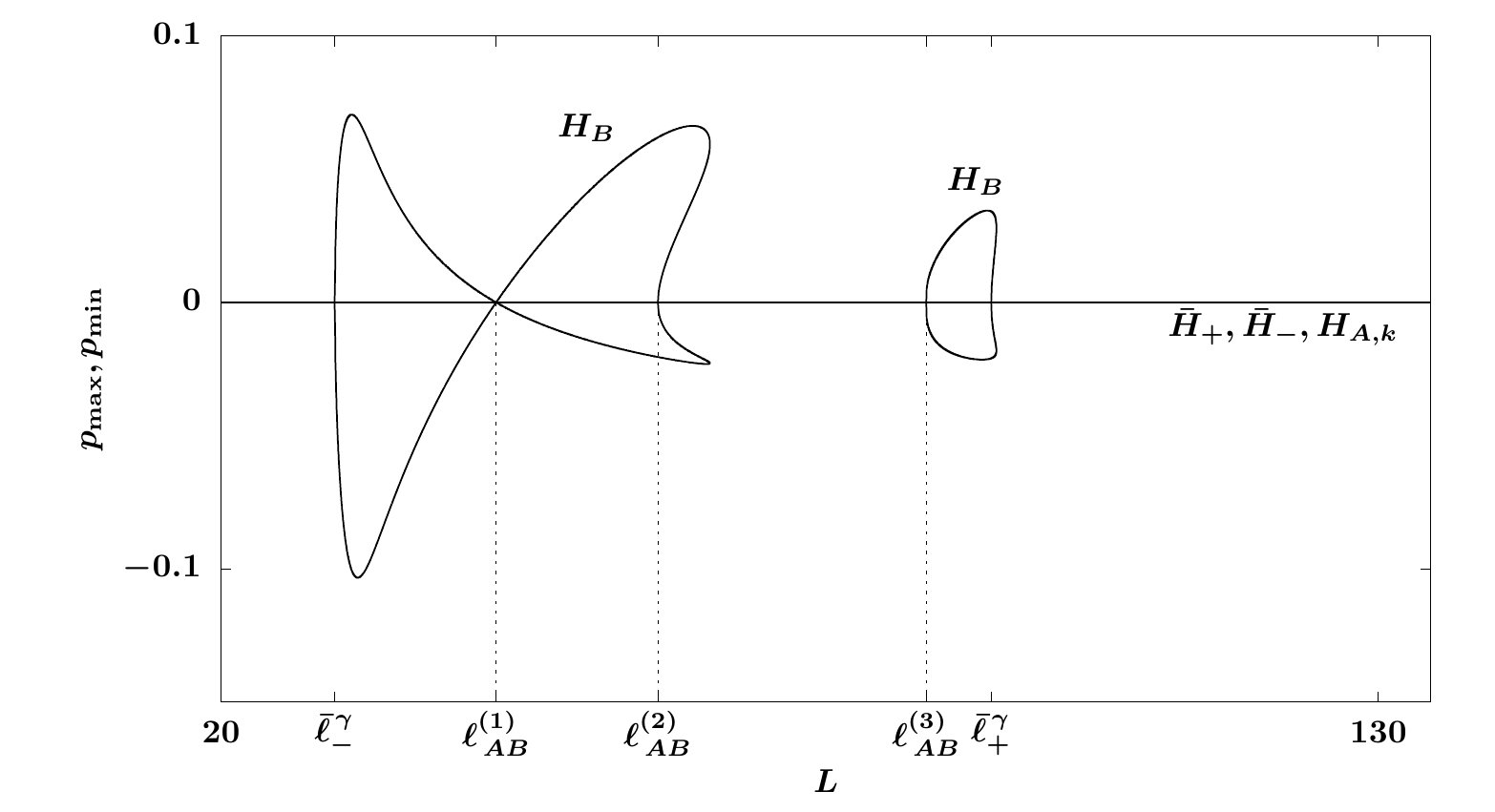}}
\vspace{-0.1in}
\hbox{\hspace{0.25in}\includegraphics[width=13.1cm,height = 6cm]{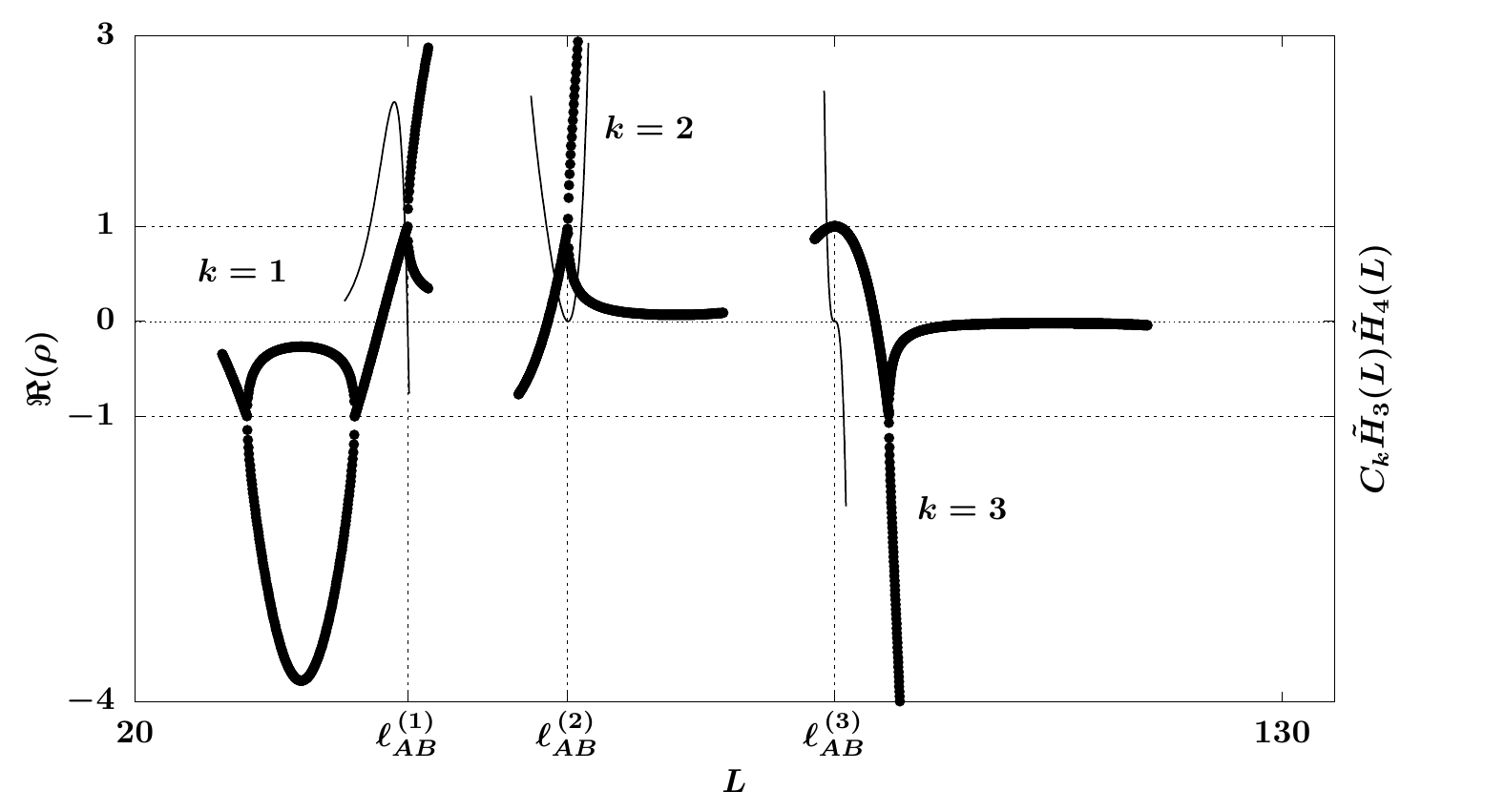}}
\vspace{-0.05in}
\caption{(Top and Middle) Bifurcation diagrams for the average of the solution $\langle H \rangle$ (top) and extremal values of the pressure $p$ (middle) parametrized by the domain size $L$ with  $\gamma=0.05$.
(Bottom) The real part of the characteristic multipliers $\rho$ (dots) of the second-order system \eqref{floquet3} and the rescaled product $C\tilde{H}_3(L)\tilde{H}_4(L)$ (thin curves) from \eqref{eqn:B4}, plotted against $L$ in three sections for the base states $\HA$, $H_{A,2}$ and $H_{A,3}$, respectively. 
Three critical periods $L = \ell_{AB}^{\,(k)}$ for secondary bifurcations of $\HB$ from $H_{A,k}$ are identified by $\rho = 1$ and $\tilde{H}_3(L)\tilde{H}_4(L) = 0$.
The critical lengths are given by $\lA \approx 28.4$, $\lB \approx 30.8$, $\lBc=93.2$, $\ell_{AB}^{\,(1)}\approx 46.2$, $\ell_{AB}^{\,(2)}\approx 61.5$. $\ell_{AB}^{\,(3)}\approx 87.1$, 
$\ell_B^{\fld,1}\approx 66.5$, $\ell_B^{\fld,2}\approx 93.7$.
The scaling parameters for $\tilde{H}_3(L)\tilde{H}_4(L)$ are $C_1 = 6\times 10^{-7}$, $C_2 = 0.01$, and $C_3 = 0.8$.
Higher order branches of $\HB$ ($H_{B,2}$ and $H_{B,3}$) are not shown in this figure.}
\lbl{fig:bifurcation}
\end{figure}

More interesting primary-secondary $\HB$ branches emerge for small $\gamma$ values satisfying $\lBc > 2\lA$. One such example for $\gamma = 0.05$ is shown in Fig.~\ref{fig:bifurcation}~(top).
In this case, the overall $\HB$ branch takes the form of two disconnected primary-secondary branches, where the first branch starts from the uniform state $\Hm$ and ends at the $H_{A,2}$ branch, and the second branch connects the uniform state $\Hc$ and the $H_{A,3}$ branch. 
\par
In particular, the first branch starts at the primary bifurcation at $L =\lB$, ends at the secondary bifurcation at $L=\ell_{AB}^{\,(2)}$, and crosses the principal second-order branch $\HA$ at a transcritical bifurcation point $L=\ell_{AB}^{\,(1)}$ similar to the primary-primary branch shown in Fig.~\ref{fig:4thss_bifurcation}~(right).
Note that the fourth-order states $\HB$ on this branch connect the $\HA$ states with a single peak and the $H_{A,2}$ states with two peaks. That is, the properties of $\HB$ change qualitatively as the branch is followed from $L=\lB$ to $L=\ell_{AB}^{\,(2)}$.
Fig.~\ref{fig:4thss_LBfold}~(left) depicts typical profiles of $\HB(x)$ on a closeup view of the bifurcation diagram from Fig.~\ref{fig:bifurcation}~(top). It shows that a secondary peak develops in the $\HB$ profiles moving along this branch
towards $L=\ell_{AB}^{\,(2)}$, and the two peaks become identical when $L=\ell_{AB}^{\,(2)}$ is reached (also see Fig.~\ref{fig:4thss_p0l0}). This observation indicates that the $\HB$ solutions from
the pitchfork bifurcation at $L=\ell_{AB}^{\,(2)}$ break the discrete $L/2$ translation symmetry of the $H_{A,2}$ solutions. The point where the secondary peak first develops is above $\ell_{AB}^{\,(1)}$ and is distinct from that point; it has not been labeled in the figures.
The fold point $L = \ell_B^{\fld,1}$ on this branch does not change the qualitative structure of the profiles, but does influence the stability of the solutions as a saddle-node bifurcation. 

\begin{figure}
\centering
\includegraphics[width=6.6cm]{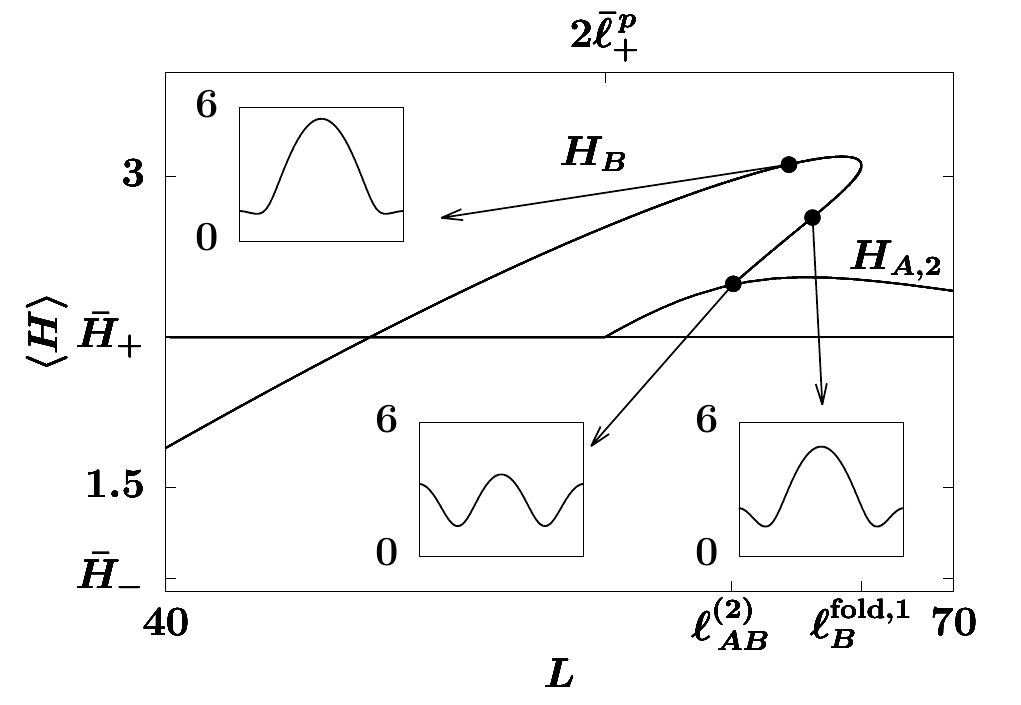}
\includegraphics[width=6.6cm]{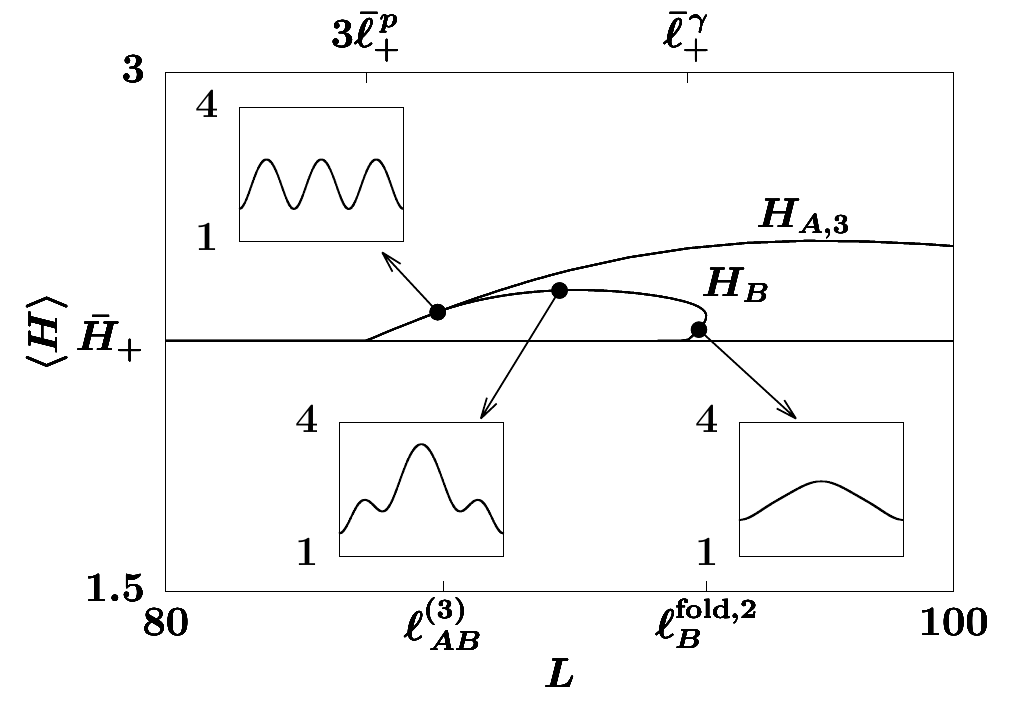}
\caption{Closeup views of the bifurcation diagram of $\left\langle H \right\rangle$ in Fig.~\ref{fig:bifurcation}~(top) parametrized by the domain size $L$, with typical profiles of fourth-order steady states on the $\HB$ branches.}
\lbl{fig:4thss_LBfold}
\end{figure}

The second $\HB$ branch connects the supercritical primary bifurcation at $L = \lBc$ and the secondary bifurcation at $L=\ell_{AB}^{\,(3)}$. No transcritical secondary bifurcations are observed on this branch. 
Typical profiles of $\HB(x)$ on this branch are plotted in Fig.~\ref{fig:4thss_LBfold}~(right). At the pitchfork bifurcation $L = \ell_{AB}^{\,(3)}$, where $\HB$ bifurcates from $H_{A,3}$, the discrete $L/3$ translation symmetry of $H_{A,3}$ is broken. This bifurcation also breaks most of the reflection symmetries about the extrema of $H_{A,3}$. The $\HB$ solutions have one large primary peak and two smaller secondary peaks that shrink as $L$ approaches $\lBc$, and transition to single-peak solutions as the branch connects to $\Hc$ at the pitchfork bifurcation $L = \lBc$. Similar to the first $\HB$ branch, the fold point $L =\ell_B^{\fld,2} > \lBc$ on this branch impacts the stability of the solutions, but does not change their symmetry. The transition point for solutions changing from three peaks to one peak has not been labeled but was found to occur after the fold point.

\par
To better characterize the coexisting states in this regime, we also plot the extremal values of the pressures $p$ in Fig.~\ref{fig:bifurcation}~(middle) corresponding to the steady states from Fig.~\ref{fig:bifurcation}~(top). The $p\equiv 0$ trivial branch corresponds to the spatially uniform steady states $\Hm$, $\Hc$, and the second-order steady states $H_{A,k}(x)$. The branches with $p\not\equiv 0$ exist for $\lB < L< \ell_B^{\fld,1}$ and $\ell_{AB}^{\,(3)}<L<\ell_B^{\fld,2}$,
and represent the two branches of fourth-order steady states $\HB$ that undergo primary bifurcations from $\bar{\ell}^{\,\gamma}_{\pm}$ and secondary bifurcations from $\ell_{AB}^{\,(k)}$ where $k = 1,2,3$. This is an alternative way of detecting transcritical bifurcation points involving $\HB$ that does not require doing the full linear stability analysis.
\par
Using the approach developed in Sec.~\ref{sec:secondaryBif}, we also identify these secondary bifurcations $\ell_{AB}^{\,(k)}$. Fig.~\ref{fig:bifurcation}~(bottom) shows the real part of the characteristic multipliers $\rho$ for the second-order system \eqref{floquet3} for $\gamma = 0.05$. The three sections correspond to a partial structure of the multipliers obtained for the states $\HA$, $H_{A,2}$, and $H_{A,3}$, respectively. The full structures for these states are similar to the one presented in Fig.~\ref{fig:floquet}~(right) for $\gamma = 0.4$.
Then for $k = 1,2,3$, we solve the fourth-order system \eqref{floquet4} using $H_{A,k}$ as the base state, and obtain the corresponding $\tilde{H}_3(L)\tilde{H}_4(L)$ as functions of $L$. Recall that an admissible critical domain size $L = \ell_{AB}^{\,(k)}$ for a secondary bifurcation from $H_{A,k}$ of period $\ell_{AB}^{\,(k)}/k$ needs to satisfy the following two conditions for the base state $H_{A,k}$: 
\begin{itemize}
    \item The fourth-order system \eqref{floquet4} gives $\tilde{H}_3(L)\tilde{H}_4(L) = 0$;
    \item There exists a characteristic multiplier $\rho = 1$ for the second-order system \eqref{floquet3}.
\end{itemize}
The rescaled $\tilde{H}_3(L)\tilde{H}_4(L)$ curves for $k = 1,2,3$ in Fig.~\ref{fig:bifurcation}~(bottom) show that their values reach zero, along with the corresponding characteristic multiplier $\rho = 1$, at the bifurcation points $L = \ell_{AB}^{\,(k)}$ identified in the top and middle panels.
It is also worth mentioning that the bifurcation points $\ell_{AB}^{\,(k)}$ with $k > 1$ can also be folded back into the characteristic multiplier plots of the principal state $\HA$ instead of $H_{A,k}$, since a multiplier satisfying $\rho^k = 1$ at $L = \ell_{AB}^{\,(k)}/k$ implies the existence of a periodic solution for $L = \ell_{AB}^{\,(k)}$.
\par 
There are very likely more complicated cases involving higher order bifurcations, yielding period doubling of $H_B$ solutions, and an extensive investigation of all the interesting structures is beyond the scope of this study.
The rest of the paper will focus on the stability and dynamical behavior of the principal branches of the solutions.

\section{Stability analysis of nonuniform steady states}
\lbl{sec:stability}
To examine the stability of the steady states, we consider an $L$-periodic positive steady state $H(x)$ over the domain $0\le x \le L$ and perturb it by setting $h(x,t) \sim H(x)+\delta \Psi(x)e^{\lambda t}$, where $\delta \ll 1$ and $\Psi(x)$ is also $L$-periodic.
We linearize equation \eqref{main} around the base state $H(x)$, a solution of \eqref{4thsteadystate}, and obtain the $O(\delta)$ equation
\begin{equation}
\lambda \Psi = \mathscr{L} \Psi,
\lbl{4thstability}
\end{equation}
where the linear operator $\mathscr{L}$ is given by
\begin{equation}\lbl{op4}
\mathscr{L}\Psi\equiv  \left[ {\gamma} + \frac{d}{d x}\left(H^3\frac{d}{dx}\right)\right] \bigg(\Pi'(H)\Psi-\frac{d^2 \Psi}{d x^2}\bigg)
 +\frac{d}{d x}\bigg[3H^2 
 \left(\Pi'(H) {dH\over dx} -{d^3H\over dx^3}\right)
 \Psi \bigg].
\end{equation}
If any eigenvalue $\lambda$ of problem \eqref{4thstability} has a positive real part, then the steady state $H(x)$ is unstable. 
Since the second-order steady states  satisfy $p(H_A)\equiv 0$, the linear operator $\mathscr{L}$ in \eqref{op4} reduces to a simpler operator $\mathscr{L}_A$ for $H=\HA(x)$
\begin{equation}
\lbl{op2}
\mathscr{L}_A\Psi\equiv \left[ {\gamma} + \frac{d}{d x}\left(H^3\frac{d}{dx}\right)\right] \bigg(\Pi'(H)\Psi-\frac{d^2 \Psi}{d x^2}\bigg).
\end{equation} 
For the case $\gamma \le 0$ the instability of the second-order nonuniform  steady states has been studied in \cite{ji2018instability}. Similar stability analysis for second-order steady states in related equations has been performed in \cite{LP, engelnkemper2019continuation}.
Here we will extend this study to the case $\gamma > 0$ which reveals critical domain sizes that give rise to limit cycles and other interesting dynamics. 
For simplicity, for the rest of the paper we only focus on the case where only a single-period solution fits in the domain, but these results can be easily extended to states with multiple periods. For instance, symmetry breaking, droplet merging, and leveling dynamics arising from a second-order state with two periods have been numerically investigated in \cite{ji2018instability}.

\begin{figure}
\centering
\includegraphics[width=6.5cm]{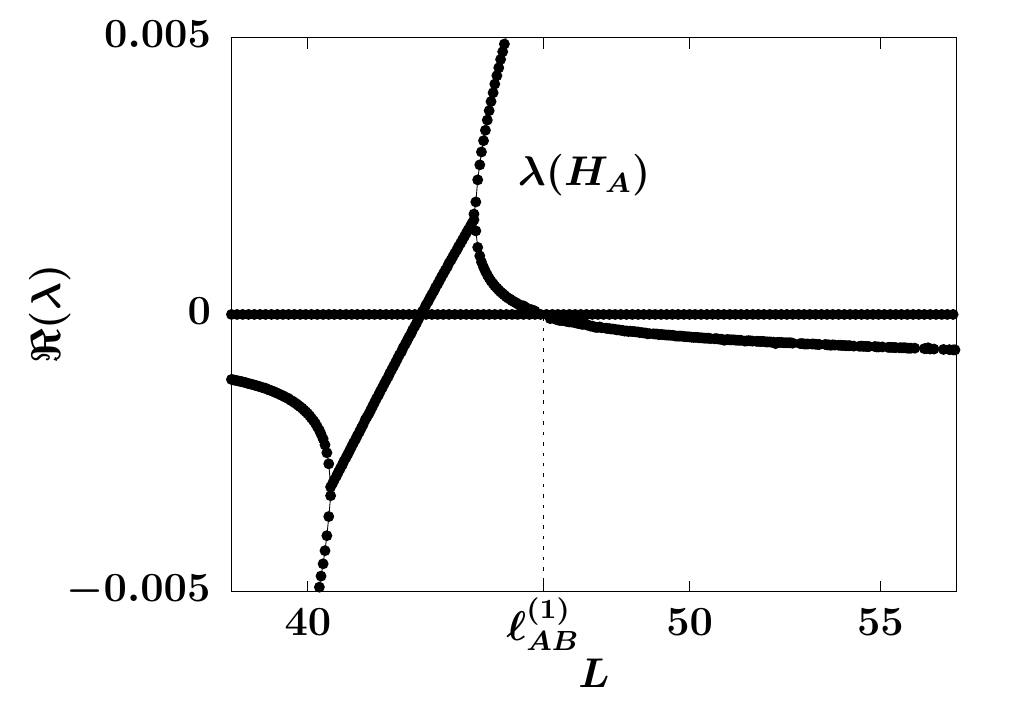}
\includegraphics[width=6.5cm]{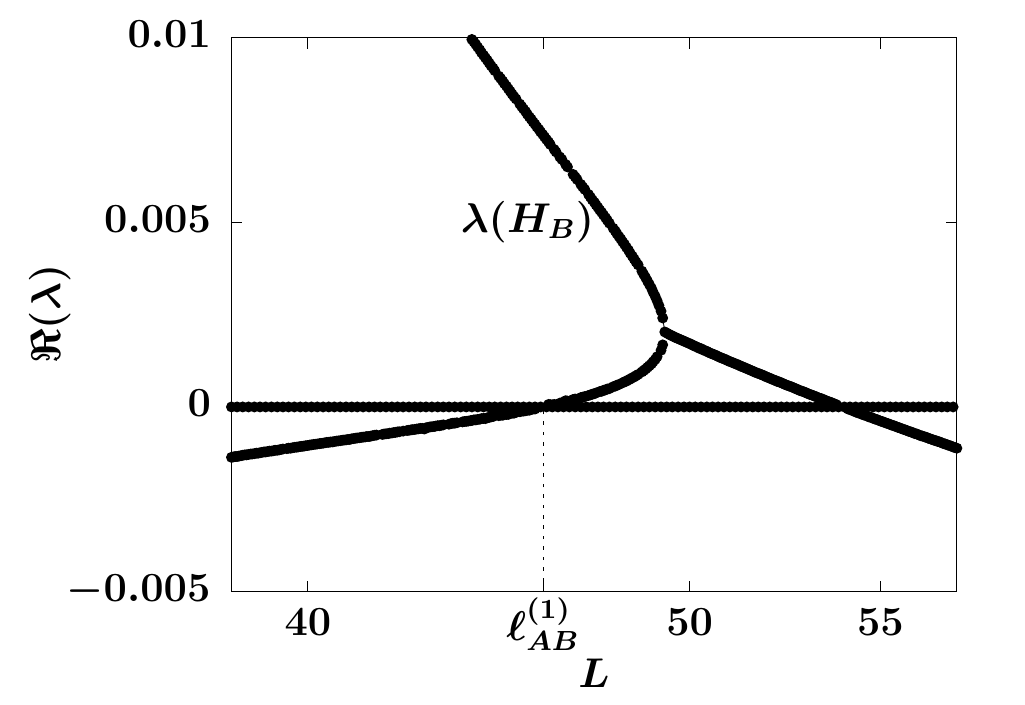}
\caption{The dependence of eigenvalues of $\HA(x)$ (left)  and  $\HB(x)$ (right) parametrized by the domain size $L$ with fixed $\gamma = 0.05$,  showing that $\ell_{AB}^{(1)}$ labelled in Fig.~\ref{fig:bifurcation} is a  transcritical bifurcation point at which both $\HA$ and $\HB$ have an eigenvalue crossing zero.}
\lbl{fig:2nd4thev}
\end{figure}

First we numerically solve the eigenproblem \eqref{4thstability} for both nonuniform steady states $\HA$ and $\HB$ for a range of domain sizes $L$. For all $L$, a translational mode $\lambda^T=0$ with $\Psi=H'(x)$ is present on the periodic domain.
The dependence of the dominant eigenvalues on $L$ is plotted in Fig.~\ref{fig:2nd4thev}, showing that each state has a dominant real eigenvalue crossing zero at the bifurcation point $\ell_{AB}^{\,(1)}$. This indicates that $\ell_{AB}^{(1)}$ labelled in Fig.~\ref{fig:bifurcation} is a transcritical bifurcation point.
In addition, both plots show pairs of real eigenvalues transitioning to become complex conjugate pairs that yield secondary Hopf bifurcations that have not been labelled in Fig.~\ref{fig:bifurcation}. Further, the fold points shown in Fig.~\ref{fig:bifurcation} indicate saddle-node bifurcations, where a real eigenvalue crosses zero.
The dependence of the dominant eigenvalues on the system parameter $\gamma$ is illustrated in Fig.~\ref{fig:beta_ev_LlL*} for a second-order state $\HA$ with a fixed domain size $L$.
\par 

These numerical results indicate that the stability of the steady states depends strongly on both $L$ and $\gamma$. Therefore, stability in plots of the solution branches parametrized by the domain size $L$ can change significantly when $\gamma$ is changed. For example, Fig.~\ref{fig:beta_ev_LlL*} shows that the stability of the $\HA$ solution with a fixed domain size undergoes qualitative changes as $\gamma$ varies, despite the fact that $\HA$ itself is independent of $\gamma$.

\par
Next we study the stability of $\HA(x)$ in the limit of weak non-conserved flux for $\gamma \to 0$. Since the stability of these states is closely related to the stability of equilibria of the mass-conserving equation \eqref{main} with $\gamma = 0$, we present below a brief overview of the $\gamma = 0$ case.

For $\gamma = 0$, the linear operator $\mathscr{L}_A$ in \eqref{op2} is self-adjoint with respect to a weighted $H^{-1}$ norm \cite{witelski2000dynamics, bernoff1998axisymmetric}. Hence from spectral theory for $\gamma = 0$ the state $\HA$ on a compact domain is associated with a spectrum that is real and discrete. Moreover, it has been shown in \cite{ji2018instability} that for $\gamma = 0$ the state $\HA$ changes stability at a critical $L = \ell^*$, where $\ell^*$ depends on the parameter $\p0$. That is, for $ \lA <L < \ell^*$ the state $\HA$ is unstable, and for $L > \ell^*$ it is stable.

\lbl{sec:2ndStability}
\begin{figure}
\centering
\includegraphics[width=10cm]{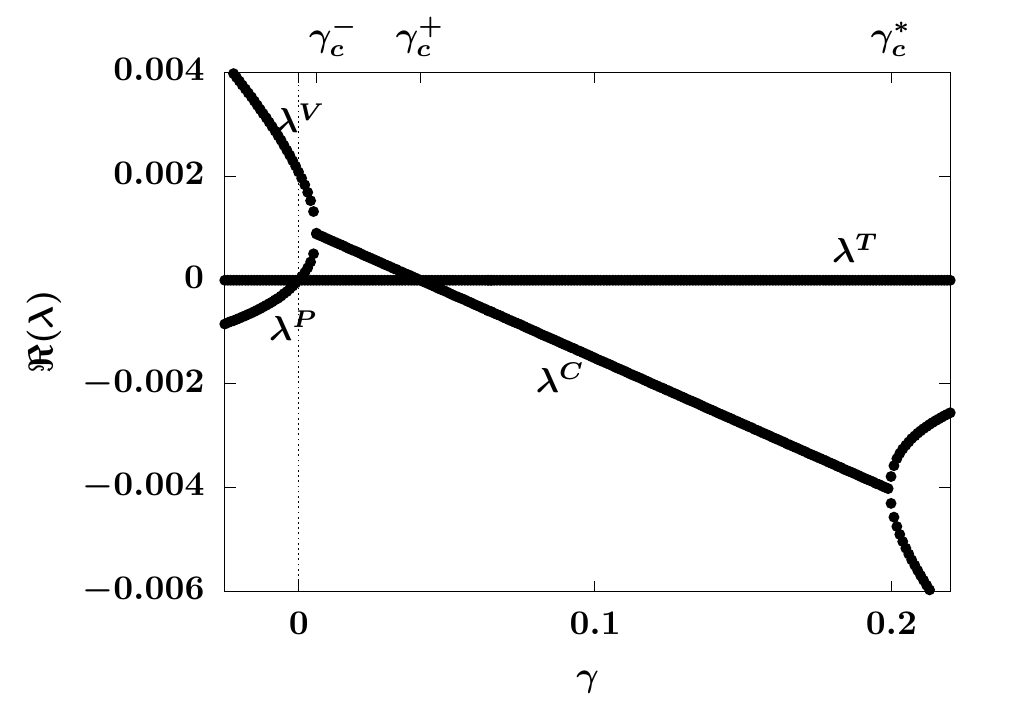}
\caption{The dependence of eigenvalues of $\HA(x)$ on the parameter $\gamma$ for $-0.025\le \gamma \le 0.22$ with the fixed period $L=29.2$  with $L < \ell^*$,
where $\ell^*\approx 30.05$.}\lbl{fig:beta_ev_LlL*}
\end{figure}

In the non-conserved case with small $\gamma > 0$, the stability of the state $\HA$ depends on both $L$ and $\gamma$.
For a fixed period $L = 29.2$ in the range $\lA<L<\ell^*$,
we numerically solve the eigenproblem \eqref{4thstability} with $\mathscr{L} = \mathscr{L}_A$ for the state $\HA(x)$.
The dominant eigenvalues of $\HA$ parametrized by the parameter $\gamma$ are plotted in Fig.~\ref{fig:beta_ev_LlL*}. The figure shows that in addition to the translational eigenvalue $\lambda^T=0$, for small positive $\gamma > 0$, $\HA(x)$ has two unstable eigenvalues $\lambda^P, \lambda^V > 0$. In the limit $\gamma \to 0^+$, $\lambda^P$ approaches zero, and $\lambda^V$ approaches a positive value given by the unstable eigenvalue of $\HA$ in the $\gamma = 0$ case.
At $\gamma_c^- \approx 0.0059$,  the two eigenvalues
merge together and become a pair of complex conjugate eigenvalues $\lambda^C, \bar{\lambda}^C$ with positive real parts.
At the critical value $\gamma_c^+\approx 0.041$, the complex conjugate pair cross the imaginary axis and
a Hopf bifurcation occurs.
The real part of the complex conjugate pair is positive for $\gamma_c^-\le \gamma \le \gamma_c^{+}$, and is negative for $\gamma_c^+\le \gamma \le \gamma_c^{*}$, where $\gamma_c^* \approx 0.2$. Therefore this suggests a limit cycle bifurcating from the steady state $\HA(x)$ at the Hopf bifurcation point $\gamma = \gamma_c^+$.
\par
In summary, on the range $\lA<L < \ell^*$ Hopf bifurcations occur at $\gamma=\gamma_c^+(L)$ with $\Re(\lambda^C(\gamma)) = 0$. To express this in terms of critical domain sizes, we define two critical periods $\ell_{AC}^{\,(1)}$ and $\ell_{AC}^{\,(2)}$ as functions of $\gamma$ such that $\ell_{AC}^{\,(1,2)}(\gamma_c^+) = L$. The two curves $\ell_{AC}^{\,(1)}(\gamma)$ and $\ell_{AC}^{\,(2)}(\gamma)$ merge at $\gamma^* \approx 0.041$. The region bounded by $\ell_{AC}^{\,(1,2)}(\gamma)$ on $0 <\gamma < \gamma^*$ is shown in the bifurcation diagram parametrized by $(\gamma, L)$ in Fig.~\ref{fig:bifurcation_limitcycle}.

\begin{figure}
\centering
\includegraphics[width=10cm]{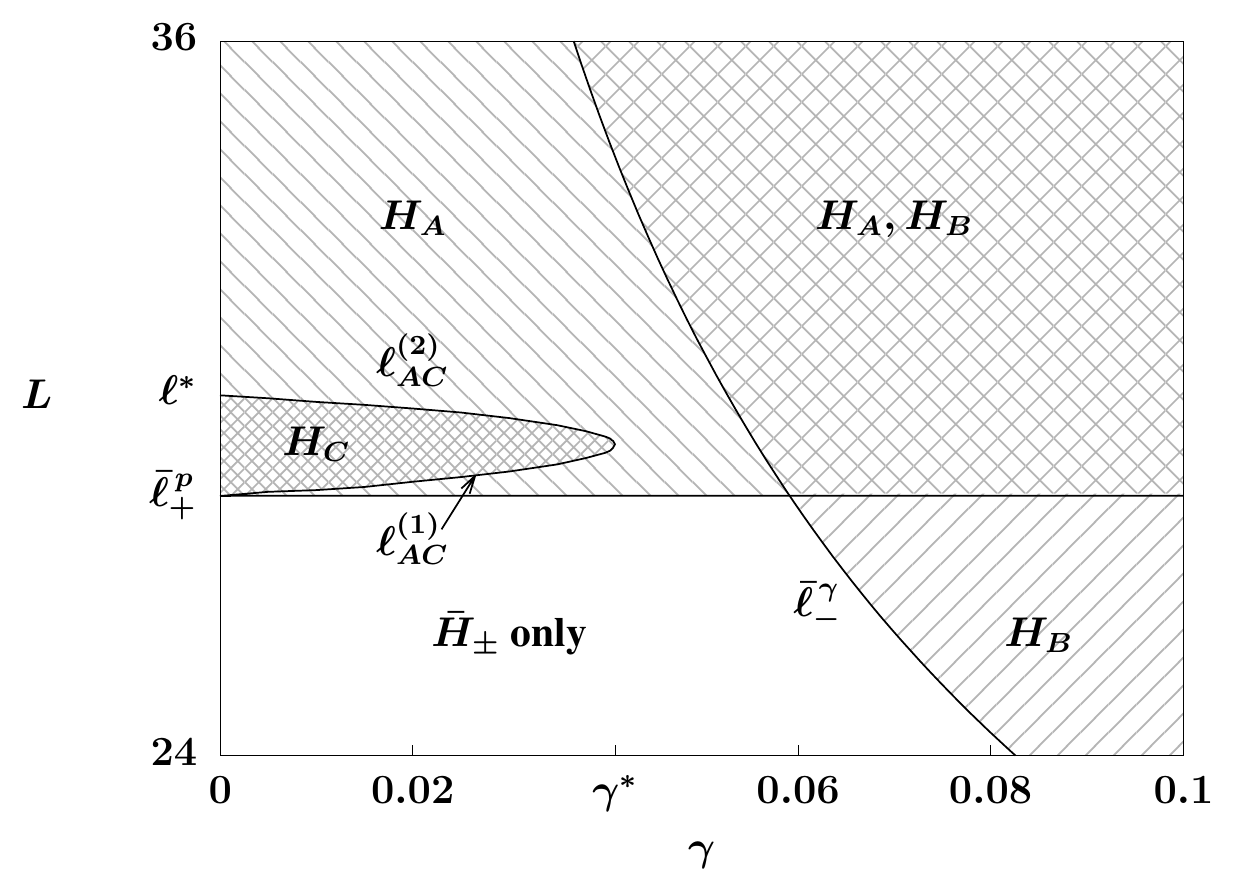}
\caption{Bifurcation diagram of coexisting states:  principal second-order states $\HA(x)$ for $L>\lA$, fourth-order states $\HB(x)$ for $L>\lB(\gamma)$, a finite region for stable limit cycles $\HLC(x,t)$ parametrized by $(\gamma, L)$, and spatially uniform states $\bar{H}_{\pm}$ for all $\gamma,~ L$.}
\lbl{fig:bifurcation_limitcycle}
\end{figure}

\section{Limit cycle dynamics}

\lbl{sec:limitcycle}

\begin{figure}
\centering
\includegraphics[width=6.5cm]{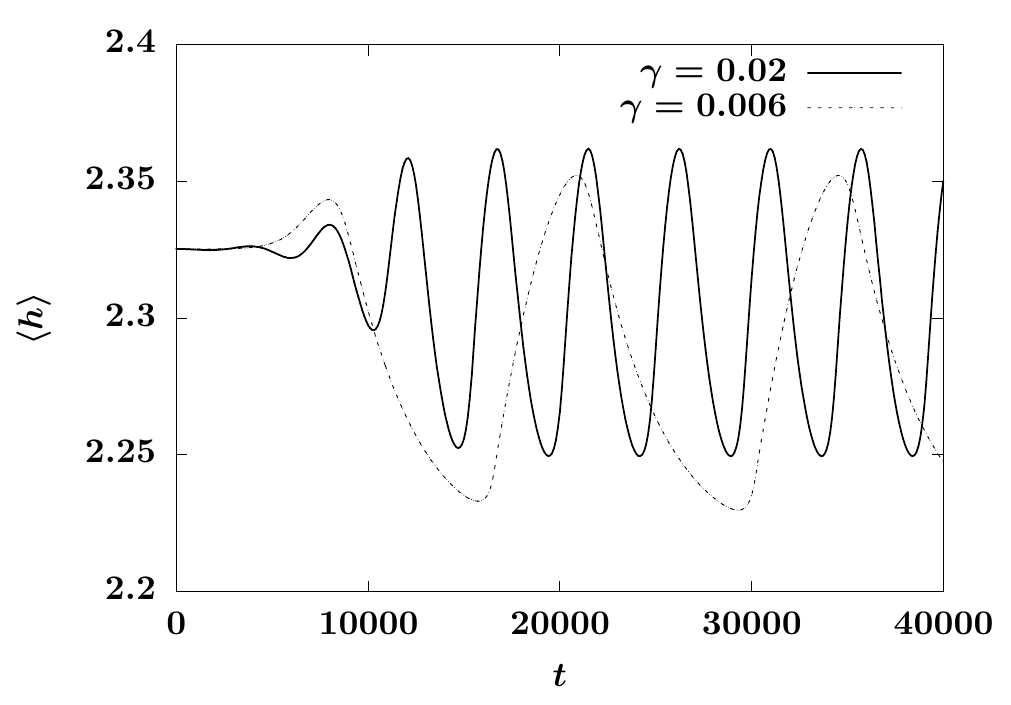}
\includegraphics[width=6.5cm]{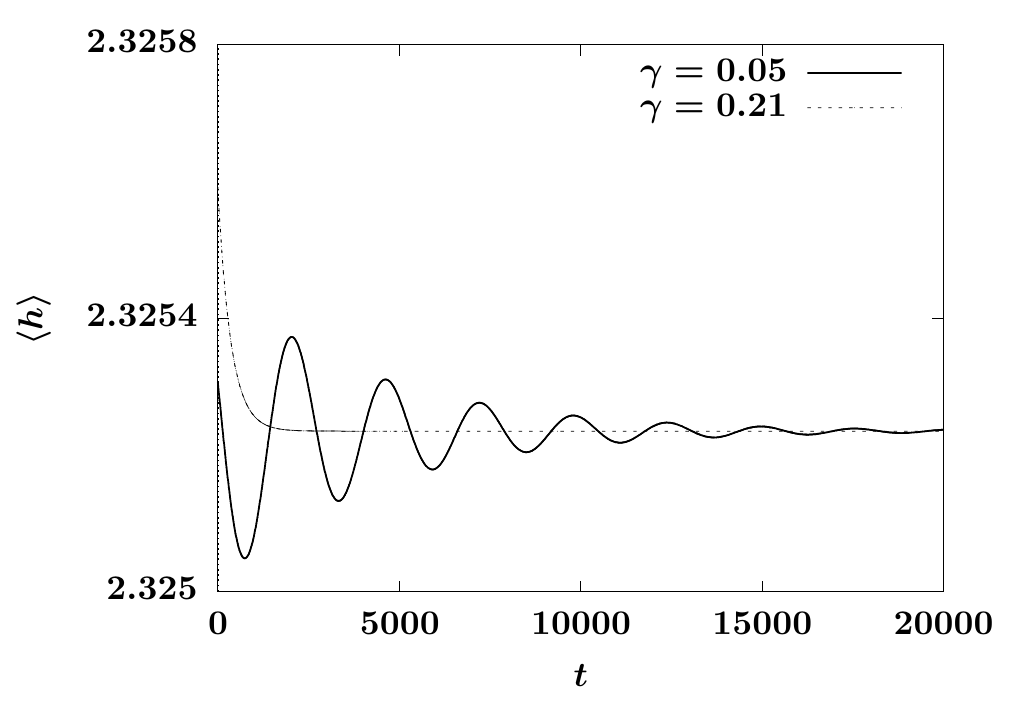}
\caption{Evolution of $\langle h\rangle$ in long-time simulations of \eqref{main} starting from $h_0(x)=\HA(x)+0.01 \Psi(x)$ showing 
(left) limit cycle dynamics for $\gamma=0.02$ and $\gamma=0.006$ (right) oscillatory damping for $\gamma=0.05$ and exponential convergence to $\langle \HA\rangle$ with $\gamma = 0.21$. The domain size is $L=29.2$ with  $\lA<L<\ell^*$.}
\lbl{fig:limitcycle}
\end{figure}
\begin{figure}
\centering
\includegraphics[width=10cm]{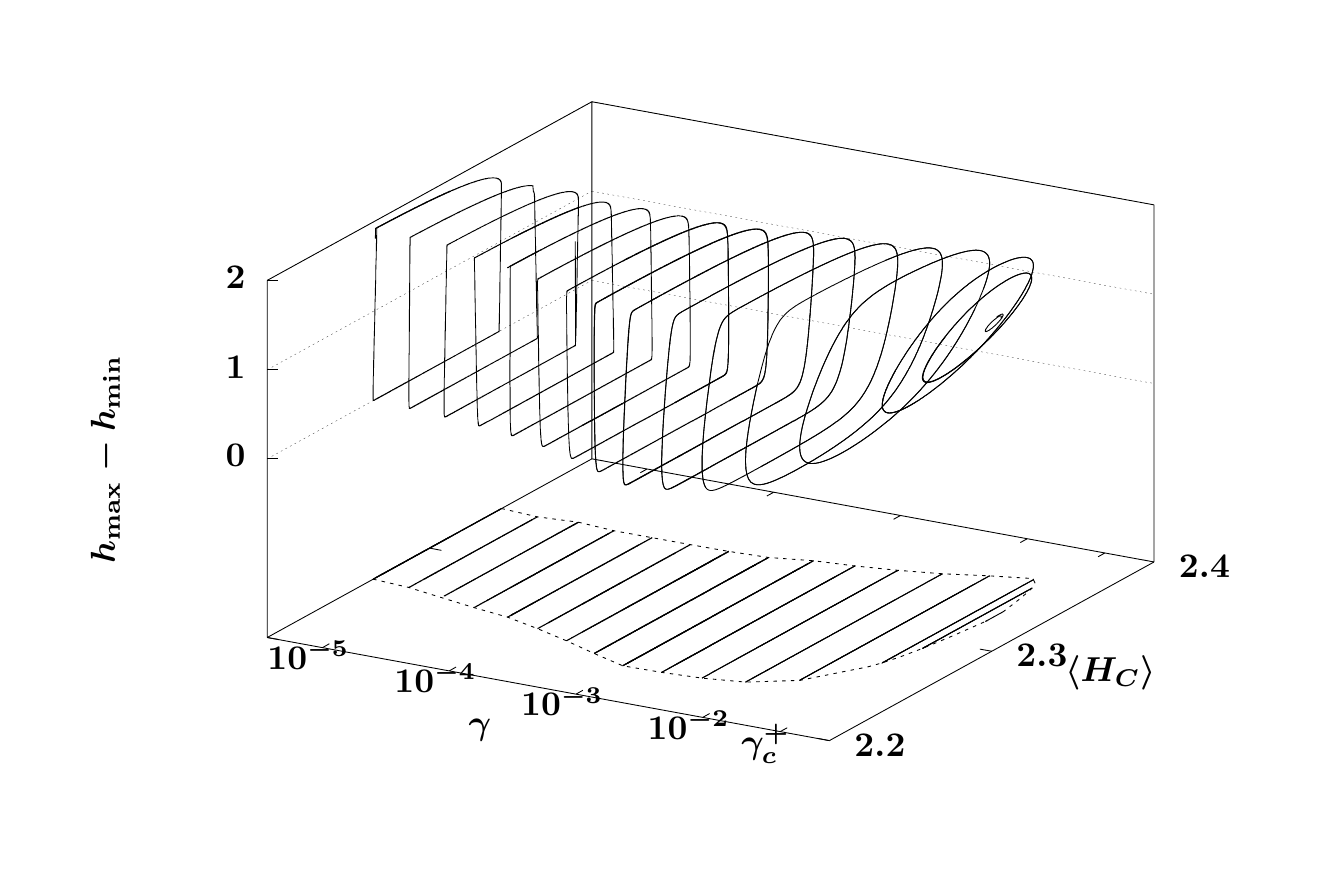}
\caption{Changes of $\langle h \rangle$ and $h_{\max}-h_{\min}$ in direct PDE simulations of stable limit cycles for a range of the parameter $\gamma$, $0<\gamma<\gamma_c^+$, with a fixed domain size $L=29.2$.}
\lbl{fig:gam_ev_LlL*_limitcyclegroup}
\end{figure}
We will show that stable limit cycles occur in the finite region marked as $H_C$ above in Fig.~\ref{fig:bifurcation_limitcycle}.
The dynamics of the observed stable limit cycles can be described by plotting some global properties, such as the average of the solution $\langle h \rangle$ or the extreme values of the solution profile, $h_{\min} = \min_x h$ and $h_{\max} = \max_x h$, as functions of time.
To understand the significance of the stability of the second-order steady states in this regime, we perform PDE simulations starting from initial conditions of the form
\begin{equation}
h_0(x)=H(x)+\delta\Psi(x),
\lbl{idpert}
\end{equation}
where $\Psi(x)$ is an unstable eigenmode of $H = \HA$ associated with an eigenvalue $\lambda$ with $\Re(\lambda) > 0$, and $\|\Psi\|_2 = 1$, $\Psi(0) \ge 0$ and $|\delta| \le 1$. 
At the Hopf bifurcation $\gamma = \gamma_c^+$ a stable limit cycle bifurcates from $\HA$ in the form
\begin{equation}
    \HLC(x,t) \sim \HA(x) +  \cos(\omega t) f(x)\sqrt{\gamma_c^+-\gamma}\qquad 
    \mbox{for $\gamma< \gamma_c^+$,}
\end{equation}
where the period of the limit cycle oscillations is determined by the parameter $\gamma$.
\par 
Fig.~\ref{fig:limitcycle}~(left) shows the evolution of two such limit cycle simulations parametrized by $\langle h \rangle$ with $\gamma = 0.02$ and $\gamma = 0.006$, where both $\gamma$ values satisfy $\gamma < \gamma_c^+$.
It is observed that the period of the limit cycle with $\gamma = 0.02$ is significantly smaller than the one with $\gamma = 0.006$. For $\gamma^+_c<\gamma<\gamma_c^*$, the real part of the complex conjugate pair $\lambda^C, \bar{\lambda}^C$ becomes negative (see Fig.~\ref{fig:beta_ev_LlL*}), and the corresponding simulation with $\gamma = 0.05$ in Fig.~\ref{fig:limitcycle}~(right) shows that the initial spatial perturbations exhibit oscillatory decay over time. For $\gamma > \gamma_c^*$, the dominant eigenvalues of the steady state $\HA(x)$ are real and negative. The corresponding PDE simulation in Fig.~\ref{fig:limitcycle}~(right) with $\gamma = 0.21 > \gamma_c^*$ shows that the dynamic solution converges monotonically to the stable steady state $\HA(x)$.
\par 
Using numerical continuation with respect to the parameter $\gamma$, we also identify a continuous family of stable limit cycles for $0<\gamma< \gamma_c^+$. Fig.~\ref{fig:gam_ev_LlL*_limitcyclegroup} shows the dependence on the parameter $\gamma$ of two properties of the limit cycles, the profile magnitude $h_{\max}-h_{\min}$ and the average of the solution $\langle h \rangle$. As $\gamma\to 0$, the magnitude of the limit cycle oscillations approaches a positive limit magnitude with an increasing period. For $\gamma > \gamma_c^{+}$ the stable limit cycle vanishes, and the dynamic solution converges to the stable steady state $\HA(x)$ (see Fig.~\ref{fig:limitcycle}~(right)).

\begin{figure}
\centering
\includegraphics[width=9cm]{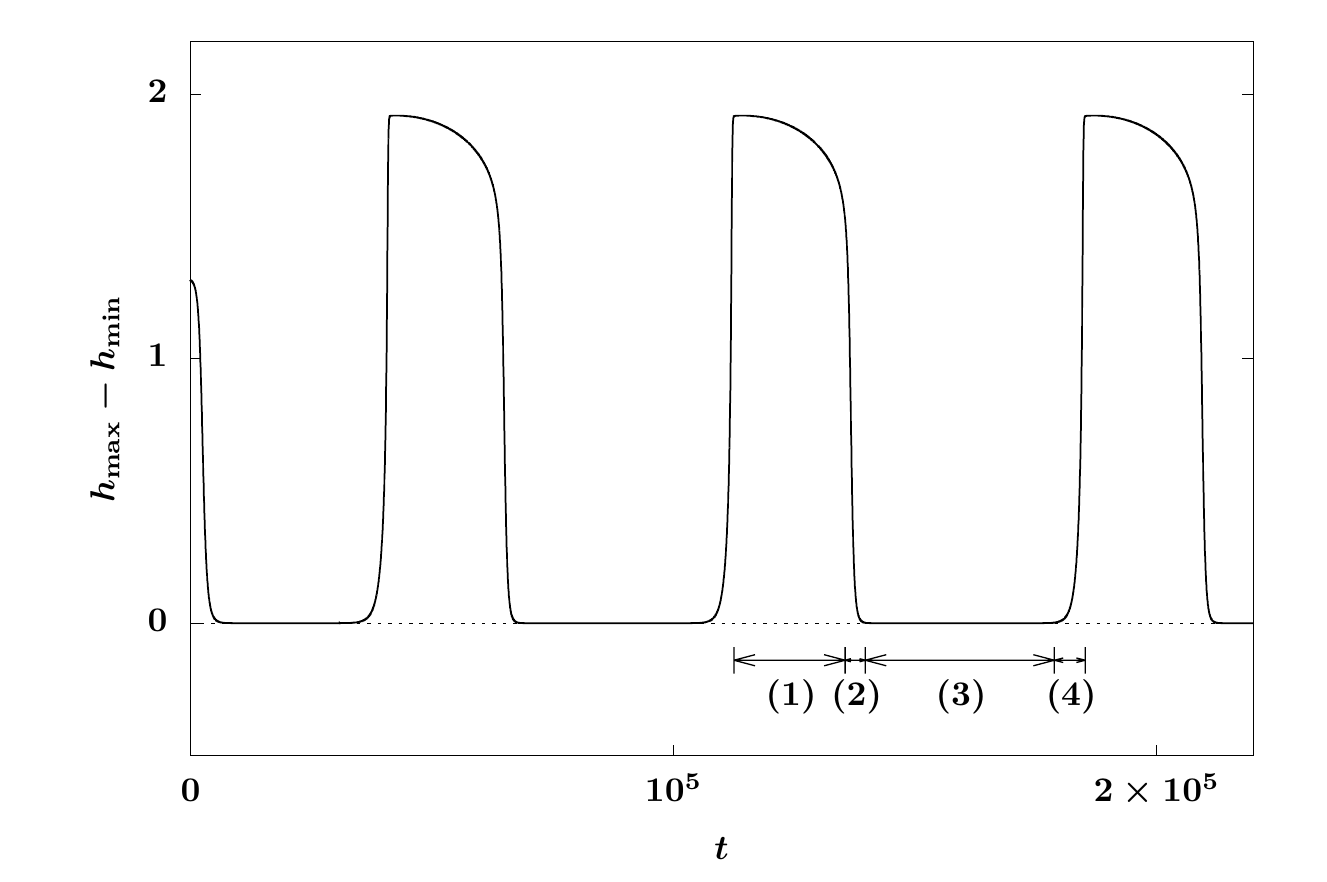}
\caption{Time-profile of $h_{\max}-h_{\min}$ illustrating the four stages in the limit cycle dynamics
starting from the initial data \eqref{idpert} with $H(x) = \HA(x)$ 
and system parameters $L = 29.2$ and $\gamma = 0.001$. }
\lbl{fig:limitcycle_example}
\end{figure}

\begin{figure}
\centering
\includegraphics[width=6.5cm,height=4.8cm]{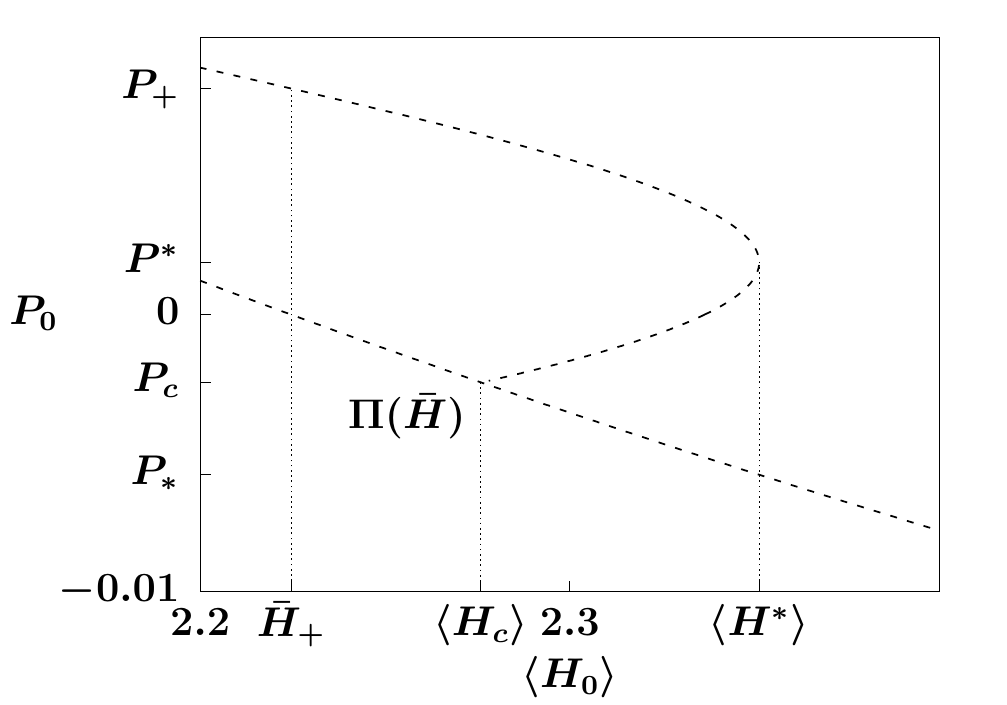}
\includegraphics[width=6.5cm,height=4.8cm]{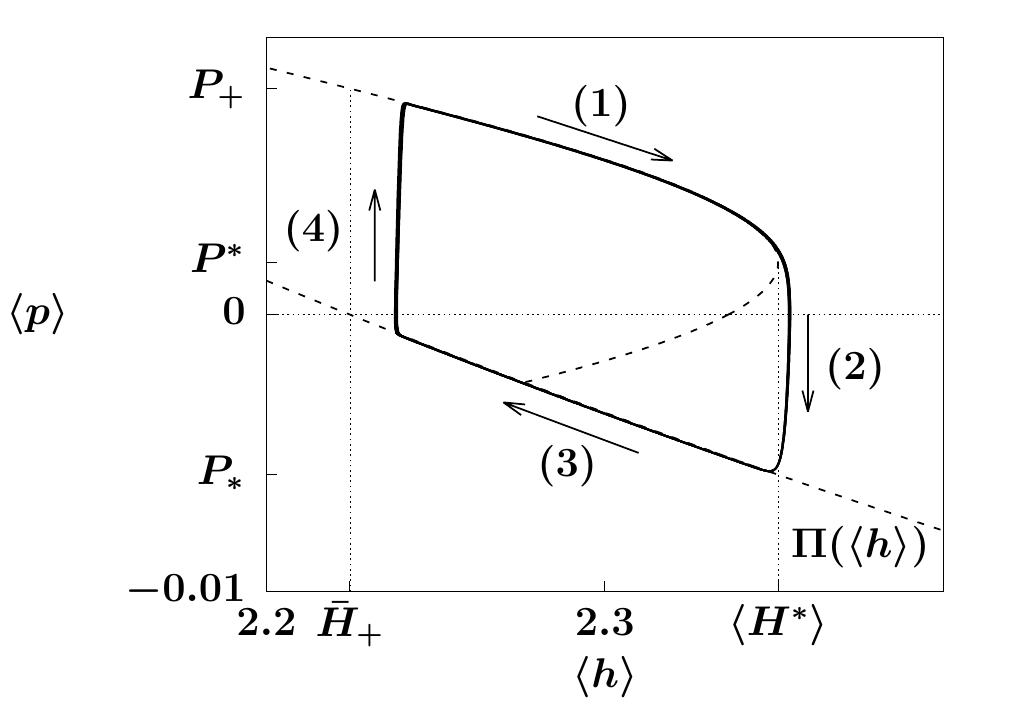}
\caption{(Left) Bifurcation diagram of coexisting solutions to \eqref{eqn:p0ODE} parametrized by $(\langle H_0 \rangle, P_0)$; (Right) The average pressure $\langle p \rangle$ plotted against the average of the solution $\langle h \rangle$ (in solid curves) for the PDE simulation in Fig.~\ref{fig:limitcycle_example}, compared against the slow manifold (in dashed curves) shown in the left plot.
}
\lbl{fig:limitcycle_analysis}
\end{figure}
\par 

Typical relaxation ocillator-type limit cycles in fast/slow dynamical systems \cite{murdock} consist of four stages (shown in Fig.~\ref{fig:limitcycle_example}), where the slow-time solution evolves on a slow manifold in stages $(1)$ and $(3)$, and alternates with fast dynamics occurring in stages $(2)$ and $(4)$. We begin by determining the slow-time solution. By rescaling the original PDE \eqref{main} using 
\begin{equation}
   T = \gamma t , \qquad h(x,t)=H(x,T), \qquad p(h)=P(H), \nonumber
\end{equation}
we get the full PDE system for $(H, P)$
\begin{equation}
\gamma {\partial H\over \partial T} =\frac{\partial}{\partial x} \left(H^3 \frac{\partial P}{\partial x}\right) + \gamma P, \qquad P = \Pi(H) - \frac{\partial^2H}{\partial x^2}.
\lbl{eqn:slowPDE}
\end{equation}
Using regular perturbation expansions for both $H$ and $P$,
\begin{equation}
    H(x, T) = H_0(x,T) + \gamma H_1(x,T) + O(\gamma^2),\nonumber
\end{equation}
\begin{equation}
    P(x,T) = P_0(x, T) + \gamma P_1(x, T) + O(\gamma^2),\nonumber
\end{equation}
and collecting the $O(1)$ and $O(\gamma)$ terms after substituting the expansions in \eqref{eqn:slowPDE}, we obtain the leading order slow system for $\gamma \to 0$
\begin{equation}
    \frac{\partial}{\partial x}\left(H_0^3\frac{\partial P_0}{\partial x}\right) = 0, \qquad P_0 = \Pi(H_0) - \frac{\partial ^2H_0}{\partial x^2},
    \lbl{eq:slow}
\end{equation}
and the order $O(\gamma)$ system
\begin{equation}
   \frac{\partial H_0}{\partial T} = \frac{\partial}{\partial x}\left(H_0^3\frac{\partial P_1}{\partial x}\right) + P_0, \qquad P_1 = \Pi'(H_0)H_1 - \frac{\partial ^2H_1}{\partial x^2}.
    \lbl{eq:slowgamma}
\end{equation}
The leading order system \eqref{eq:slow} gives a quasi-steady solution, and imposing periodic boundary conditions leads to a spatially-uniform leading order pressure term $P_0 \equiv P_0(T)$.
With different values of $P_0$ there are coexisting solutions to the second-order ODE,
\begin{equation}
P_0 = \Pi(H_0) - \frac{d^2H_0}{dx^2},
\lbl{eqn:p0ODE}
\end{equation}
which have been studied extensively in \cite{ji2018instability}.
In addition to the constant solution $\bar{H}$ that satisfies $P_0 =\Pi(\bar{H})$, a family of nonuniform solutions $H_0(x)$ also exists for a range of $P_0$. Fig.~\ref{fig:limitcycle_analysis}~(left) shows a bifurcation diagram of these solutions parametrized by their average $\langle H\rangle$ with the domain size $L = 29.2$. The branch of nonuniform solutions $H_0(x)$ bifurcates from the constant states $\bar{H}$ at a critical $P_0 = P_c$, which corresponds to a critical solution average $P_c = \Pi(\langle H_c\rangle)$. 
Moreover, there exists a critical $P^*>0$ such that for $P_0 > P^*$ the average of the nonuniform solution $\langle H_0\rangle$ monotonically increases as $P_0$ decreases. While for $P_c<P_0<P^*$ the nonuniform solution has a decreasing average $\langle H_0\rangle$ as $P_0$ decreases. As a result, for any average in the range $\langle H_c\rangle<\langle H\rangle < \langle H^*\rangle$, there exist two nonuniform solutions $H_0$ satisfying \eqref{eqn:p0ODE} with different values of $P_0$. With the solution average satisfying $\langle H\rangle < \langle H_c\rangle$, only one nonuniform solution $H_0$ coexists with the constant solution $\bar{H}$.
\par 
Integrating both sides of the $O(\gamma)$ equation \eqref{eq:slowgamma}$_1$ over the periodic domain, one obtains
\begin{equation}
   \frac{d\langle H_0\rangle}{dT} = P_0.
\end{equation}
This allows us to qualitatively determine the dynamics of the average $\langle H_0\rangle$ of the solution on the slow manifold, as $\langle H_0(T)\rangle$ is increasing in time for $P_0 > 0$, and is decreasing for $P_0 < 0$. Restricting the solution to stay on the slow manifold, then starting from any quasi-steady solution $H_0(x)$, with $P_0 < 0$ the solution will be driven to the critical point $(\Hc, 0)$ in Fig.~\ref{fig:limitcycle_analysis}~(left), and for $P_0 > 0$ the solution evolves towards the other critical point $(\langle H^*\rangle, P^*)$.
Fig.~\ref{fig:limitcycle_analysis}~(right) plots the limit cycle dynamics shown in Fig.~\ref{fig:limitcycle_example} parametrized by the average pressure $\langle p\rangle = L^{-1}\int_0^L p(h)~dx$ and the average $\langle h\rangle$ of the dynamic solution. 
It shows that in the slow stages $(1)$ and $(3)$ the PDE solution evolves on the slow manifold in the direction based on the sign of $\langle p\rangle$ as described above.
\par
When the two critical points $(\Hc, 0)$ and $(\langle H^*\rangle, P^*)$ are approached, the dynamic solution jumps off the slow manifold and is governed by the fast-time problem, which can be described using the original PDE \eqref{mainpde} and \eqref{pressure}. For $\gamma \to 0$ we apply the regular perturbation expansions for both $h$ and $p$,
\begin{equation}
    h(x, t) = h_0(x,t) + \gamma h_1(x,t) + O(\gamma^2),\nonumber
\end{equation}
\begin{equation}
    p(x,t) = p_0(x, t) + \gamma p_1(x, t) + O(\gamma^2),\nonumber
\end{equation}
and obtain the leading order fast problem
\begin{equation}
   \frac{\partial h_0}{\partial t}= \frac{\partial}{\partial x}\left(h_0^3\frac{\partial p_0}{\partial x}\right), \qquad p_0 = \Pi(h_0) - \frac{\partial ^2h_0}{\partial x^2},
    \lbl{eq:fast}
\end{equation}
which is identical to the mass-conserving thin film equation \eqref{mainnoevap}. This indicates that for the fast dynamics in the limit cycle, the leading order average of the solution $\langle h_0 \rangle$ is constant in time when $p_0$ is evolving, which explains the motions in the fast stages (2 and 4) in Fig.~\ref{fig:limitcycle_analysis}~(right).

\section{Finite-time singularity formation}
\lbl{sec:singularity}
\par
Next, we illustrate how the linear instability of a fourth-order steady state can lead to another fundamental mode of the dynamics of \eqref{main} -- the formation of a finite-time rupture singularity \cite{eggers2015singularities}.
Fig.~\ref{fig:4th_rup}~(left) depicts a typical numerical simulation of \eqref{main} for $\gamma=0.05$ and $L=50$ with initial data $h_0(x) = \HB(x)+0.03\Psi_B(x)$, where $\Psi_B(x)$ is an unstable eigenmode associated with the complex conjugate pair eigenvalues $\lambda_B \approx 0.0017 \pm \imagi 0.0013$.
After initial transients, the minimum decreases and approaches $h\to 0$ at an isolated point $x_c$ at a
finite critical time, $t_c$. Since $\Pi(h)\to -\infty$ as $h\to 0$, this is a singularity of the dynamics and the solution cannot be continued beyond the critical time. The simulation also suggests that the singularity formation is \textit{localized}, that is, away from $x_c$, the solution remains smoothly and slowly evolving even as properties become singular in a small neighborhood of $x_c$
as $t\to t_c$.
\par
In \cite{ji2017finite} it was shown that first-kind self-similar dynamics lead to this singularity, and we briefly summarize key points here. Let $\tau$ be the time remaining until rupture, $\tau=t_c-t$, we can seek a focusing self-similar solution for rupture 
at $x_c$ of the form
\begin{equation}
    h(x,t)\sim \tau^\alpha  H(\eta)\qquad \eta= {x-x_c\over \tau^\beta}\qquad \mbox{for $\tau\to 0$},
    \lbl{similarity}
\end{equation}
where the positive scaling exponents $\alpha, \beta$ must be determined. Substituting this ansatz into \eqref{main} will not balance all terms in the PDE, but by considering the limit $\tau\to 0$ with $H=O(1)$, we can seek an \textit{asymptotically} self-similar solution that satisfies a leading-order dominant balance of terms. Since the singularity cannot form without the non-conserved term \cite{ji2018instability}, we know the balance must involve part of the $\gamma p$ term balancing the rate of change and other terms as $\tau\to 0$; at leading order this determines the scaling exponents as $\alpha=1/5$ and $\beta=3/10$ and reduces to the ODE for $H(\eta)$,
\begin{equation}
\left({\ts \frac{1}{5}}H-{\ts \frac{3}{10}}\eta {d H\over d\eta}\right)+{d\over d\eta} \left(H^3 {d\over d\eta} \left[{1\over H^4}\right]\right)-\frac{\gamma}{H^4}=0.
\lbl{sim_ODE}
\end{equation}

\begin{figure}
\centering
\includegraphics[width=6.5cm,height = 4.5cm]{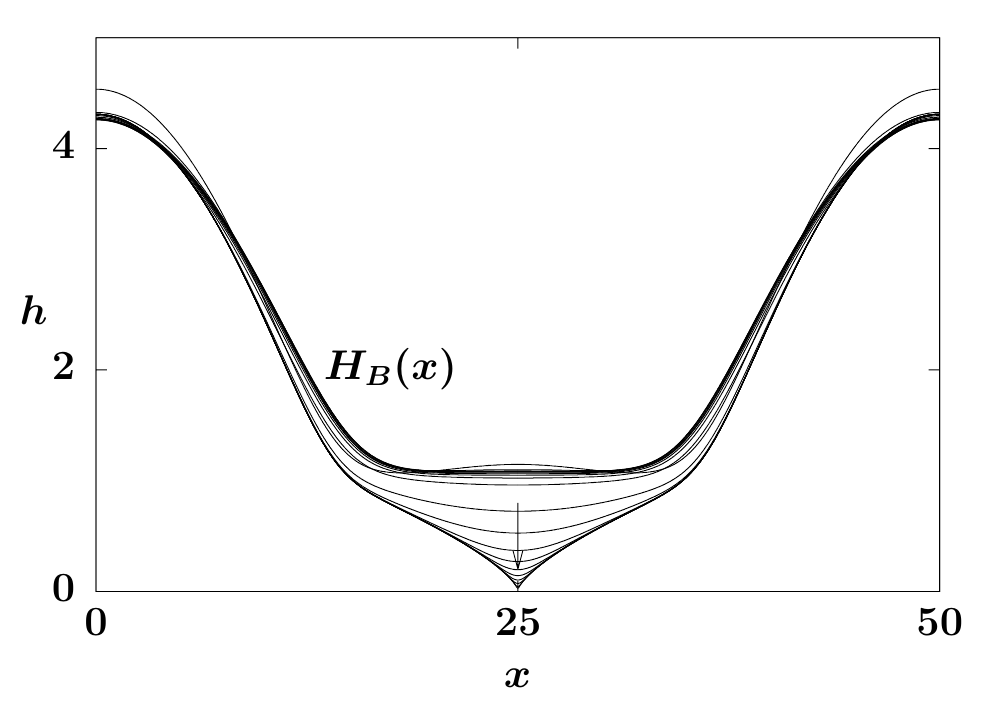}
\includegraphics[width=6.5cm, height = 4.5cm]{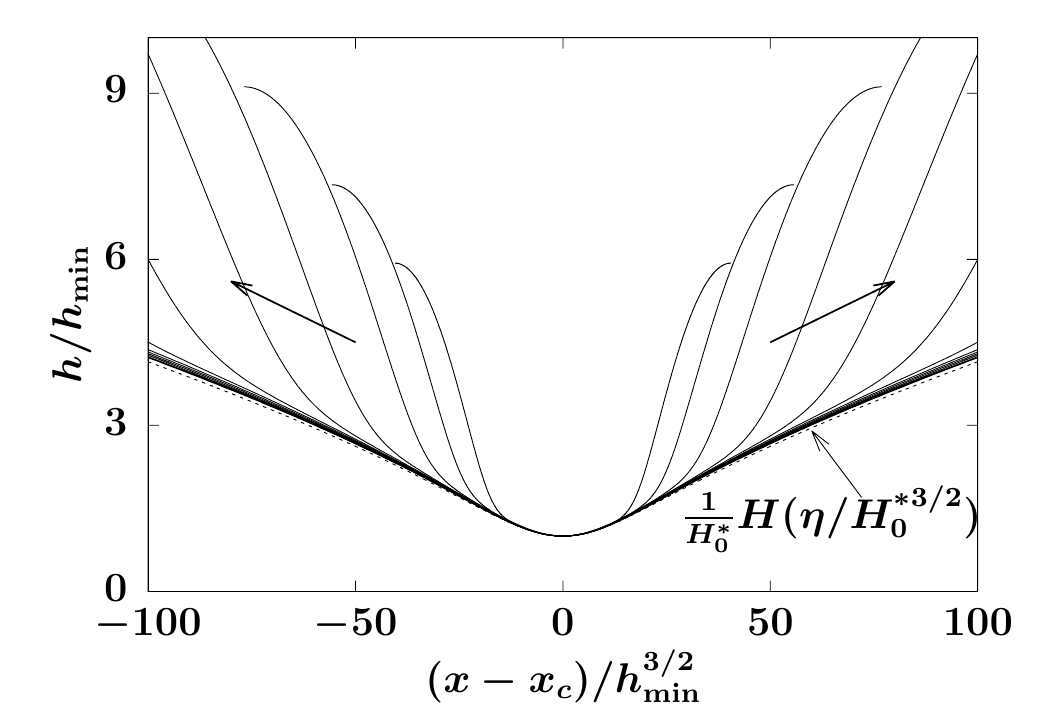}
\caption{(Left) Numerical solution of \eqref{main} with $\p0=\gamma=0.05$ starting from $\HB(x)$ with small perturbations leading to finite-time rupture at $t_c\approx 497.37$; (Right) the numerical profiles rescaled by $h_{\min}(t)$ converging to the similarity solution scaled as $H(\eta/H_0^{*3/2})/H_0^*$ of the leading order ODE.}
\lbl{fig:4th_rup}
\end{figure}
\par 
The localized nature of the self-similar dynamics implies slow evolution away from $x_c$,
namely $h_t=O(1)$ for $|x-x_c|=O(1)$ as $\tau\to 0$ \cite{witelski2000dynamics}. This yields a far-field boundary condition on the similarity solution,
\begin{equation}
{\ts \frac{1}{5}}H-{\ts \frac{3}{10}}\eta H_{\eta} = 0 \qquad \mbox{ as $|\eta| \to \infty.$}
\lbl{sim_bc}
\end{equation}
Using a one-parameter shooting method with $H(0)=H_0$ and $H'(0) = 0$ we  numerically solve the symmetric similarity equation \eqref{sim_ODE} subject to \eqref{sim_bc}. For the unique value of $H(0)=H_0^*\approx 0.711$, we obtained a self-similar solution $H(\eta)$ of this second order boundary value problem, corresponding to the effective
PDE 
$$\partial_t h=\partial_x (h^3\partial_x h)-\gamma/h^4.$$ 
This is in contrast with problems leading to a fourth-order problem for similarity solutions, where multiple solutions can co-exist \cite{witelski2000dynamics}.
\par
A comparison between the PDE profiles scaled using $h_{\min}(t) = \min_x h(x,t)$ and the similarity solution in Fig.~\ref{fig:4th_rup}~(right) shows that the scaled PDE solution converges to $H(\eta)$ as $t \to t_c$. As rupture is approached, the energy \eqref{dissipation} will be dominated by the local contribution near the rupture point,
$\EE\sim {1\over 3} \int h^{-3}\,dx=O(\tau^{-3/10})\to\infty$. The fact that the energy is increasing (and diverging) as $t\to t_c$ again shows that influence of the non-conserved term in breaking the gradient flow structure of \eqref{generalModel} can be dominant in the dynamics. In \cite{ji2017finite}, a fixed value for $\p0$ was used  
in the range $\p0>\peps $ so no steady states would exist and rupture would occur from all generic initial conditions. In this paper, we have $\p0$ in the critical range, and in the next section we will discuss how different dynamics and transitions occur depending on the initial data.

\section{Numerical simulations of dynamic transitions}
\lbl{sec:transition}
We have shown that the existence and stability of the coexisting steady states crucially depend on the domain size and the system parameters.
The complex bifurcation structures also gives the potential for interesting dynamic transitions among these states.  For instance, Fig.~\ref{fig:bifurcation}~(top) depicts a sequence of bifurcation points on the domain size $L$, $\lA < \lB < \ell_{AB}^{\,(1)} < 2\lA < \cdots$, dividing the system into different regimes. With the domain size in each $L$-intervals, equation \eqref{main} has distinctive coexisting steady states and dynamic behaviors. 
\par 
Here we focus on the parameters used in Fig.~\ref{fig:bifurcation} and illustrate typical dynamics of the equation \eqref{main} for two choices of $L$, $L=40$ and $L=50$. In both cases, the principal second-order and fourth-order steady states $\HA(x)$ and $\HB(x)$ coexist with the two constant equilibria $\Hm$ and $\Hc$. However,  we will show that the stability of these states change, leading to qualitatively different dynamics. 

\par
To study the effects of the dominant eigenmodes of a steady state $H(x)$, where $H(x) = \HA(x)$ or $H(x) = \HB(x)$, for $|\delta| \ll 1$ we perform PDE simulations of \eqref{main} starting from perturbed initial conditions \eqref{idpert},
where $\Psi(x)$ is a normalized eigenmode associated with an unstable eigenvalue $\lambda$ and satisfies $\|\Psi\|_2=1$ and $\Psi(0)>0$.
Consequently, the linearized dynamics takes
the form $h(x,t) = H(x)+\delta \Psi(x)e^{\lambda t}$ when close to the steady state $H(x)$. The growth or decay of the perturbation in time can be quantified by the $L^2$-distance of the dynamic solution $h(x,t)$ and the corresponding steady state $H(x)$,
\begin{equation}
\|h(x,t)-H(x)\|_2 = |\delta|\, \|\Psi\|_2\, e^{\Re(\lambda) t}, \mbox{ where } \lambda = \mbox{argmax}_{\lambda \neq 0}\Re(\lambda(H)).
\lbl{growthrate}
\end{equation}

\begin{figure}
\centering
\includegraphics[width=6.5cm]{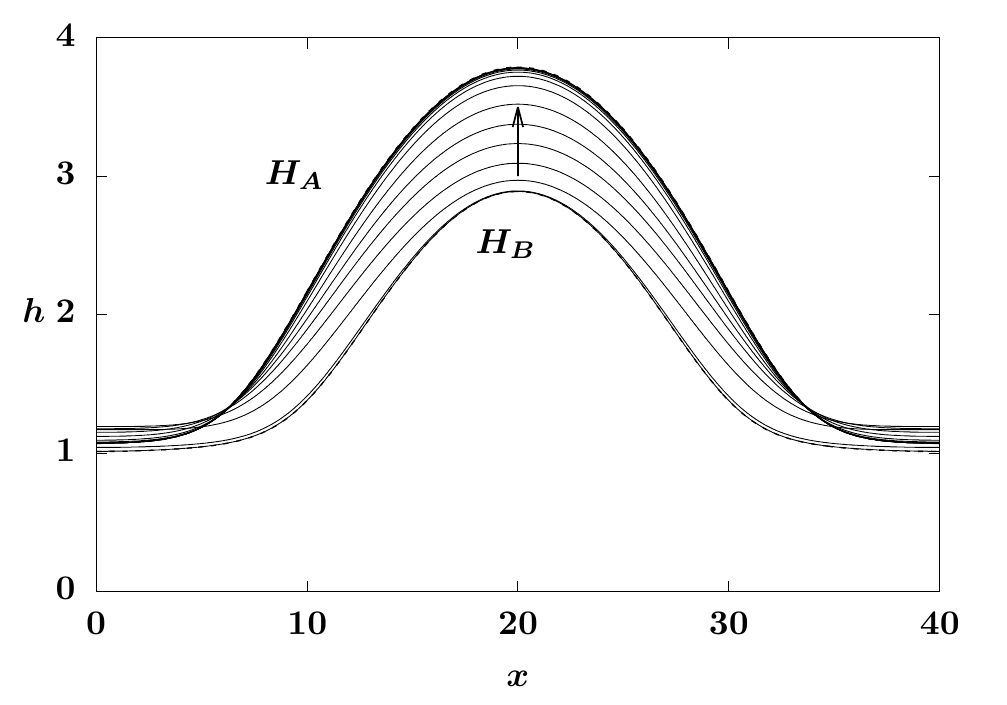}
\includegraphics[width=6.5cm]{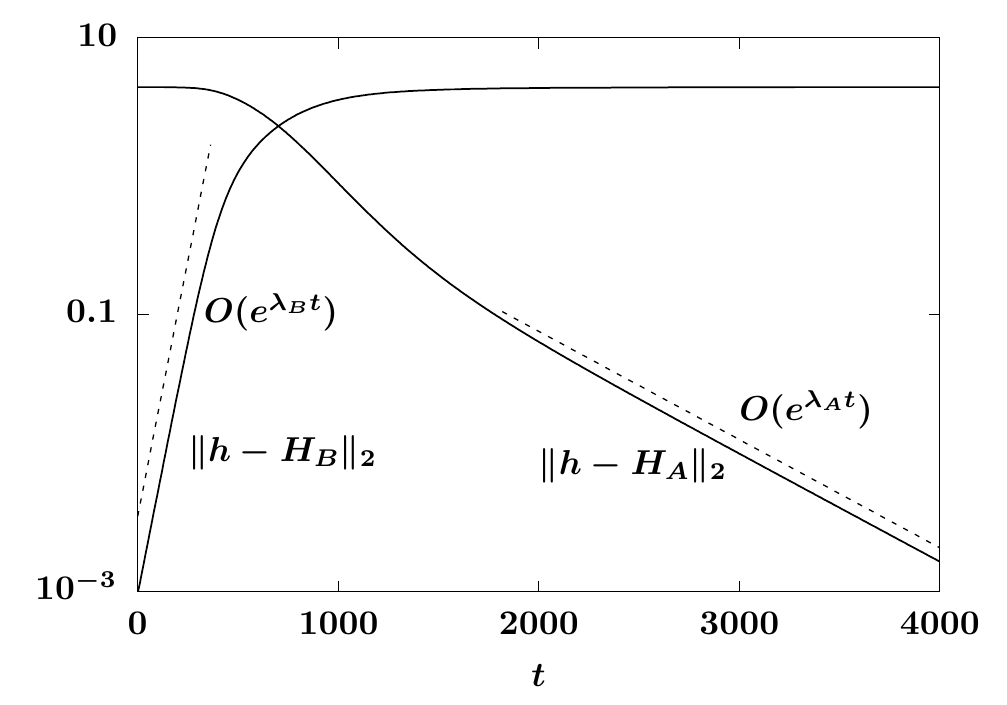}
\caption{(Left) Direct PDE simulation of \eqref{main} starting from the initial condition close to the fourth-order steady state $\HB(x)$, $h_0(x) = \HB(x) - 0.01 \Psi_B(x)$ with the domain size $L = 40$, showing convergence to the second-order steady state $\HA(x)$; (Right) Growth and decay rates of the $L^2$-norm between the PDE solution and the steady states $\HB(x)$ and $\HA(x)$.}
\lbl{fig:h_4th_L40}
\end{figure}
\begin{figure}
\centering
\includegraphics[width=6.5cm]{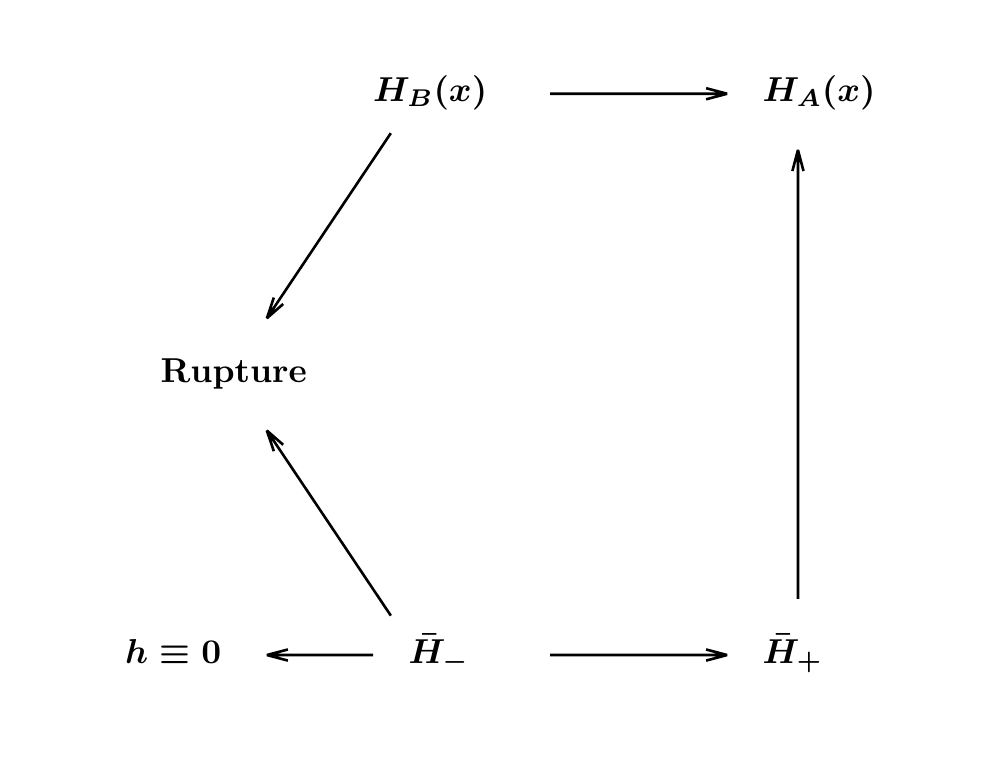}
\includegraphics[width=6.5cm]{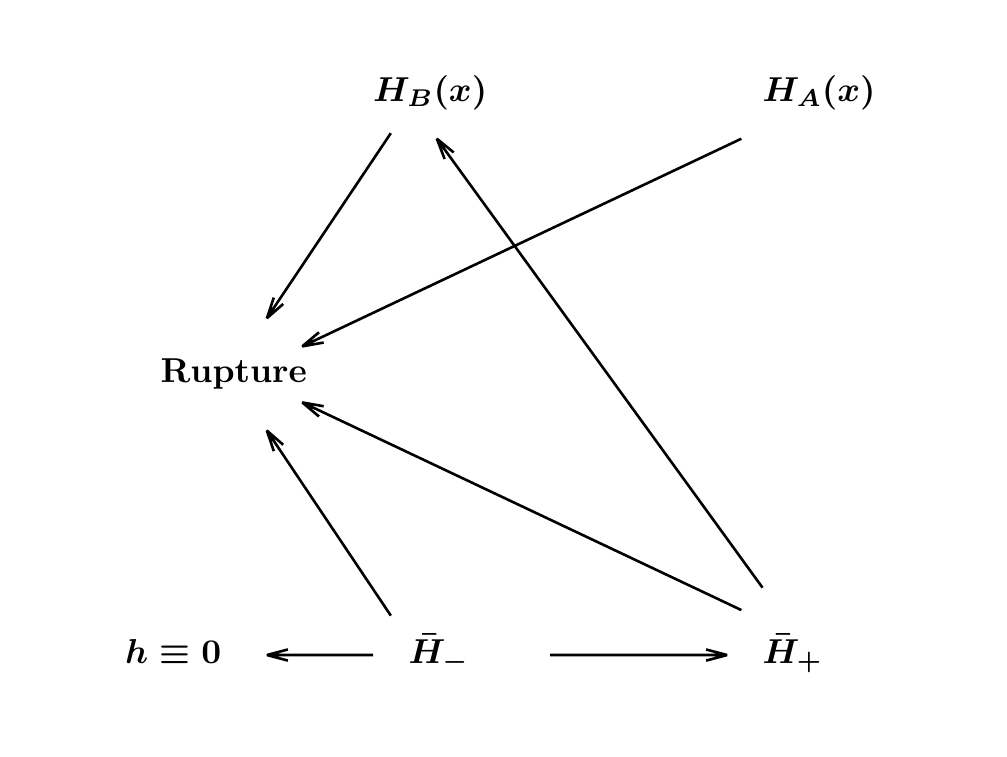}
\caption{Transition-state diagrams with (left) $L = 40$  and (right) $L=50$.}
\lbl{fig:L40-L50phase}
\end{figure}

For the first simulation, we use the domain size $L = 40$
where the second-order steady state $\HA(x)$ is linearly stable up to translations, while the fourth-order steady state $H_{B}(x)$ is linearly unstable. As shown in Fig.~\ref{fig:2nd4thev},
the largest nonzero eigenvalue of $\HA(x)$ is $\lambda_A \approx -0.0018$, and the dominant unstable eigenmode of $\HB(x)$, $\Psi_B(x)$, is associated with $\lambda_B \approx  0.017$. 
Fig.~\ref{fig:h_4th_L40}~(left) shows the PDE simulation of \eqref{main} starting from the initial data \eqref{idpert} with $\Psi = \Psi_B(x)$ and $\delta = -0.01$. In this simulation, the dynamic solution evolves away from $H_{B}(x)$ and approaches the stable steady state $H_{A}(x)$. 
Fig.~\ref{fig:h_4th_L40}~(right) gives the $L^2$-distances between the PDE solution and the two steady states. In the early stage, the distance $\|h(x,t)-\HB(x)\|_2$ grows exponentially in time following \eqref{growthrate} at the growth rate $\lambda = \lambda_B$. As the solution approaches the steady state $\HA(x)$, the distance $\|h(x,t) - \HA(x)\|_2$ decays exponentially at the rate $\lambda = \lambda_A$.
If we perturb the unstable steady state $H_{B}(x)$ in the opposite direction, and specify the initial data using \eqref{idpert} with $\Psi = \Psi_B(x)$ and $\delta = 0.01$, then finite-time rupture phenomenon as shown in Fig.~\ref{fig:4th_rup} will occur.
\par 
Figure \ref{fig:L40-L50phase}~(left) gives a transition-state diagram that summarizes the dynamics among the coexisting states in this system with $L=40$. In addition to the transitions from $\HB(x)$ to $\HA(x)$ and from $\HB(x)$ to finite-time rupture induced by the most unstable eigenmodes $\Psi_B$, it also includes the dynamical transitions starting from spatially-uniform solutions.
Spatially-uniform solutions $\bar{h}(t)$ of the PDE \eqref{main} are governed by the ODE
\begin{equation}
    \frac{d\bar{h}}{dt} = \gamma \Pi(\bar{h}),
    \lbl{constantSoln}
\end{equation}
which is obtained by dropping the spatial derivative terms from \eqref{main}. Since $\gamma > 0$, the evolution of such solutions is determined by the sign of $\Pi(\bar{h})$. Therefore, for $\p0$ in the critical range, spatially uniform perturbations around the uniform steady state $\Hm$ can lead the solution to approach either the uniform state $\Hc$ or the quenched state $h\equiv 0$. We also note that unbounded growth, $\bar{h} \to \infty$, is not possible in \eqref{constantSoln} in the critical range. PDE simulations starting from spatially-perturbed initial data near the $\Hc$ state lead to convergence to the stable steady state $\HA(x)$. 
Finite-time rupture is also numerically observed when the initial data is given by $\Hm$ with spatially nonuniform perturbations. 
\par 
In contrast, for the domain size $L = 50$ 
the four coexisting steady states $\HA(x)$, $\HB(x)$, $\Hm$ and $\Hc$ are all unstable. Fig.~\ref{fig:L40-L50phase}~(right) presents the corresponding transition-state diagram for $L = 50$, showing that all of these steady states with perturbations can lead to finite-time rupture. Again we include arrows from $\Hm$ to $\Hc$ and $h\equiv 0$ that correspond to transitions induced by spatially uniform perturbations governed by \eqref{constantSoln}.

\begin{figure}
\centering
\includegraphics[width=8.5cm]{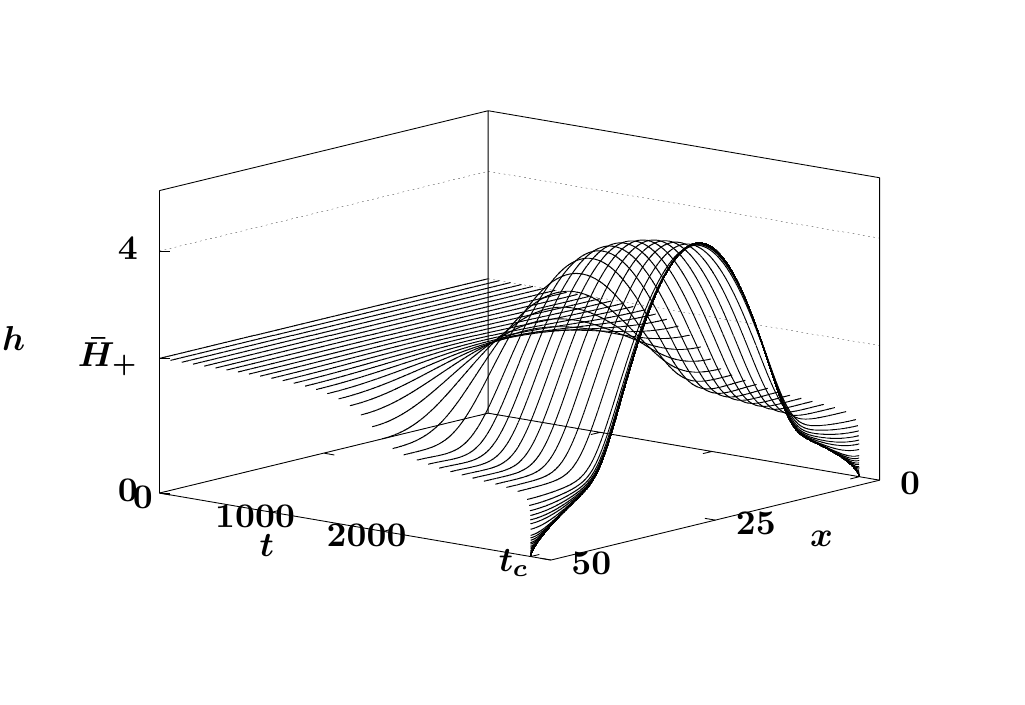}\\
\includegraphics[width=6.5cm]{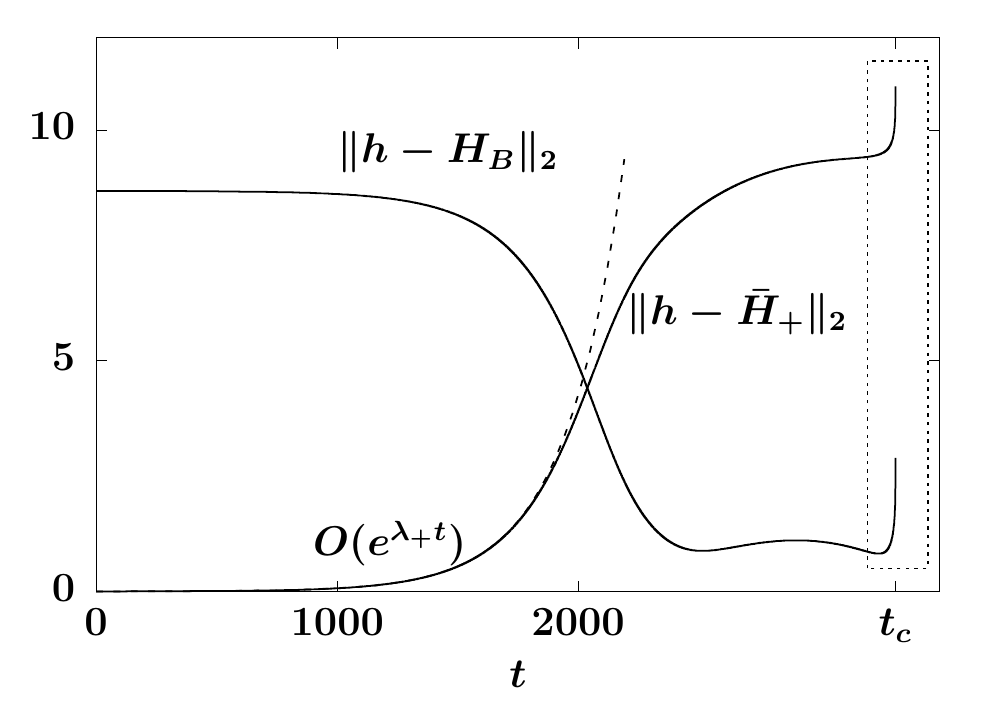}
\includegraphics[width=6.5cm]{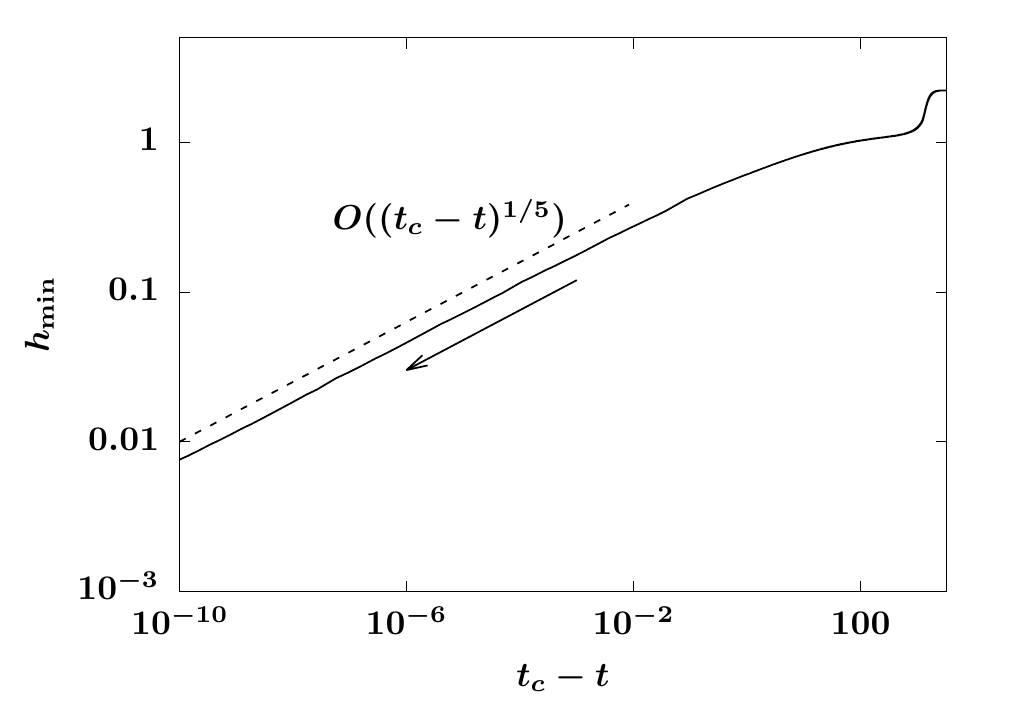}
\caption{(Top) Direct PDE simulation of \eqref{main} starting from $h_0(x) = \Hc + 0.01\sin(2\pi x/L)$ with $L = 50$ showing that the solution evolves towards $H_{B}$ in the early stage, and rupture occurs at a critical time $t_c \approx 3317.1$. (Bottom left) The corresponding evolution of $L^2$-distances between the PDE solution $h(x,t)$ and the steady states $H_{B}(x)$ and $\Hc$ capturing the transient dynamics between the two states. (Bottom right) Plot of $h_{\min}(t)$ as $t\to t_c$ showing agreement with finite-time rupture, $h_{\min} = (t_c-t)^{1/5}H(0)$, corresponding to the behavior highlighted in the dashed box in the bottom left plot.}
\lbl{fig:h_4th_L50}
\end{figure}
Moreover, we include an arrow from $\Hc$ to $\HB(x)$ to represent transient dynamics between the two states.
An example of such dynamics is depicted in Fig.~\ref{fig:h_4th_L50}~(top) where the initial perturbations around $\Hc$ lead the PDE solution to evolve towards the fourth-order state $\HB(x)$, whose instability then yields finite-time singularity at the edges of the domain. The $L^2$-distance between the PDE solution and the uniform steady state $\Hc$ in Fig.~\ref{fig:h_4th_L50}~(bottom left) shows that $\|h(x,t) - \Hc\|_2$ grows monotonically as the $h(x,t)$ evolves away from $\Hc$ with a growth rate $\lambda_+ \approx 0.00412$. This value agrees with the prediction \eqref{dispersion} from the stability analysis for the state $\Hc$.
Meanwhile, the distance $\|h(x,t) - \HB(x)\|_2$ decreases for $0<t < 2515$, indicating that the PDE solution evolves towards $\HB(x)$ in the early stage, but eventually develops a finite time singularity at a critical time $t_c \approx 3317.1$. The self-similar nature of the rupture solution as $t \to t_c$ is captured in Fig.~\ref{fig:h_4th_L50}~(bottom right) which follows $h_{\min}(t) = O\left((t_c-t)^{1/5}\right)$ from \eqref{similarity} with $\alpha = 1/5$.

The eigenvalue plots in Fig.~\ref{fig:2nd4thev} also show that the system can have an unstable $\HA$ state and a coexisting stable $\HB$ state for a range of $L$ (say for $L = 57$), which will yield a transition from $\HA$ to $\HB$ similar to the one shown in Fig.~\ref{fig:h_4th_L40}.
For larger domain sizes, more coexisting steady states are involved in the system and other complicated dynamics like coarsening and symmetry breaking are also expected (see \cite{ji2018instability});  further work is needed on these issues.

\section{Conclusions}
\label{sec:conclusion}
We have studied interesting equilibrium structures and dynamics in the non-conserved thin film-type equation \eqref{main} with small positive $\gamma >0$. Unlike the conserved model \eqref{mainnoevap} (equivalent to \eqref{main} with $\gamma = 0$) and the volatile thin film model (equivalent to \eqref{main} with $\gamma < 0$) where only second-order nonuniform steady states exist in the system, our equation also has fourth-order nonuniform steady state solutions. This is a consequence of not having a gradient flow structure.
\par
In a critical range of the parameter $\p0$ we have investigated various local bifurcation structures that arise from the coexisting uniform and nonuniform steady states.
For instance, the interactions between second-order and fourth-order branches yield different primary-primary and primary-secondary bifurcations. Further work is needed to more fully analyze the structure of the fourth-order steady states \cite{peletier}.
Stability analysis of the steady state solutions also reveals distinctive dynamic behaviors driven by the non-conserved effects, including the appearance of limit cycles, finite-time rupture, and dynamic transitions among steady states.
\par
The weakly non-conserved limit $\gamma\to 0$ is of interest to be studied as a global bifurcation from the mass-conserving equation for $\gamma=0$. While we have shown that some properties persist in this limit, other features of the system change dramatically. For instance, second-order states exist for any value for $p$ in the case $\gamma = 0$, but they are restricted to having only $p=0$ for $\gamma >0$. The appearance of limit cycles, the fourth-order steady states, and finite-time rupture for $\gamma>0$ also motivates further studies from the perspectives of PDE analysis and dynamical systems theory.

A number of interesting questions regarding the model \eqref{main} remain to be solved. First, while this paper focuses on the weak non-conserved effect limit with $\gamma \to 0$, we are also interested in the transient pattern formation under strong non-conserved influences (with large $\gamma$ values). 

Moreover, we have used the Floquet theory to numerically identify the secondary bifurcations from the second-order states. A more careful asymptotic analysis with respect to spatial perturbations is needed to generalize these results to higher-order secondary bifurcations. For instance, we observe that the local structure for $\tilde{H}_3(L)\tilde{H}_4(L)$ for a base state $H_{A,k}$ takes the form of $\tilde{H}_3(L)\tilde{H}_4(L) = O([L-\ell_{AB}^{\,(k)}]^k)$ (see Fig.~\ref{fig:bifurcation}~(bottom)). 

The current analysis is limited to steady states and their stability on relatively small domains. For a larger-scale domain where more non-trivial steady states coexist, more complex bifurcation structures and rich pattern formation are expected and need careful investigation.

\bibliographystyle{abbrv}

\bibliography{SSME_R2.bib}

\newcommand{\noop}[1]{}
\begin{thebibliography}{10}

\bibitem{ajaev2005evolution}
V.~S. Ajaev.
\newblock Evolution of dry patches in evaporating liquid films.
\newblock {\em Physical Review E}, 72(3):031605, 2005.

\bibitem{ajaev2005spreading}
V.~S. Ajaev.
\newblock Spreading of thin volatile liquid droplets on uniformly heated
  surfaces.
\newblock {\em Journal of Fluid Mechanics}, 528:279--296, 2005.

\bibitem{ajaev2001steady}
V.~S. Ajaev and G.~M. Homsy.
\newblock Steady vapor bubbles in rectangular microchannels.
\newblock {\em Journal of Colloid and Interface Science}, 240(1):259--271,
  2001.

\bibitem{bernoff1998axisymmetric}
A.~J. Bernoff, A.~L. Bertozzi, and T.~P. Witelski.
\newblock Axisymmetric surface diffusion: dynamics and stability of
  self-similar pinchoff.
\newblock {\em Journal of Statistical Physics}, 93(3):725--776, 1998.

\bibitem{bertozzi2001dewetting}
A.~L. Bertozzi, G.~Gr{\"u}n, and T.~P. Witelski.
\newblock Dewetting films: bifurcations and concentrations.
\newblock {\em Nonlinearity}, 14(6):1569, 2001.

\bibitem{bertozzi1994lubrication}
A.~L. Bertozzi and M.~C. Pugh.
\newblock The lubrication approximation for thin viscous films: the moving
  contact line with a `porous media' cut-off of van der {W}aals interactions.
\newblock {\em Nonlinearity}, 7(6):1535, 1994.

\bibitem{bestehorn2006regular}
M.~Bestehorn and D.~Merkt.
\newblock Regular surface patterns on {R}ayleigh-{T}aylor unstable evaporating
  films heated from below.
\newblock {\em Physical review letters}, 97(12):127802, 2006.

\bibitem{burelbach1988nonlinear}
J.~P. Burelbach, S.~G. Bankoff, and S.~H. Davis.
\newblock Nonlinear stability of evaporating/condensing liquid films.
\newblock {\em Journal of Fluid Mechanics}, 195:463--494, 1988.

\bibitem{chen2002}
L.-Q. Chen.
\newblock Phase-field models for microstructure evolution.
\newblock {\em Annual Review of Materials Research}, 32:113--140, 2002.

\bibitem{chicone2006ordinary}
C.~Chicone.
\newblock {\em Ordinary differential equations with applications}, volume~34 of
  {\em Texts in Applied Mathematics}.
\newblock Springer, New York, second edition, 2006.

\bibitem{craster2009dynamics}
R.~V. Craster and O.~K. Matar.
\newblock Dynamics and stability of thin liquid films.
\newblock {\em Reviews of Modern Physics}, 81(3):1131, 2009.

\bibitem{cross2009pattern}
M.~Cross and H.~Greenside.
\newblock {\em Pattern formation and dynamics in nonequilibrium systems}.
\newblock Cambridge University Press, 2009.

\bibitem{dai2016}
S.~Dai and Q.~Du.
\newblock Weak solutions for the {C}ahn-{H}illiard equation with degenerate
  mobility.
\newblock {\em Arch. Ration. Mech. Anal.}, 219(3):1161--1184, 2016.

\bibitem{de1985wetting}
P.-G. De~Gennes.
\newblock Wetting: statics and dynamics.
\newblock {\em Reviews of Modern Physics}, 57(3):827, 1985.

\bibitem{eggers2015singularities}
J.~Eggers and M.~A. Fontelos.
\newblock {\em Singularities: Formation, Structure, and Propagation}.
\newblock Cambridge University Press, 2015.

\bibitem{garcke1996}
C.~M. Elliott and H.~Garcke.
\newblock On the {C}ahn-{H}illiard equation with degenerate mobility.
\newblock {\em SIAM J. Math. Anal.}, 27(2):404--423, 1996.

\bibitem{engelnkemper2019continuation}
S.~Engelnkemper, S.~V. Gurevich, H.~Uecker, D.~Wetzel, and U.~Thiele.
\newblock Continuation for thin film hydrodynamics and related scalar problems.
\newblock In {\em Computational Modelling of Bifurcations and Instabilities in
  Fluid Dynamics}, pages 459--501. Springer, 2019.

\bibitem{fife2000}
P.~Fife.
\newblock Models for phase separation and their mathematics.
\newblock {\em Electronic Journal of Differential Equations}, 2000(48):1--26,
  2000.

\bibitem{VSS}
V.~A. Galaktionov.
\newblock {Very Singular Solutions for Thin Film Equations with Absorption}.
\newblock {\em Studies in Applied Mathematics}, {124}({1}):{39--63}, {2010}.

\bibitem{harwin2009}
V.~A. Galaktionov and P.~J. Harwin.
\newblock On center subspace behavior in thin film equations.
\newblock {\em SIAM J. Appl. Math.}, 69(5):1334--1358, 2009.

\bibitem{garcke2017well}
H.~Garcke and K.~F. Lam.
\newblock Well-posedness of a {Cahn--Hilliard} system modelling tumour growth
  with chemotaxis and active transport.
\newblock {\em European Journal of Applied Mathematics}, 28(2):284--316, 2017.

\bibitem{glasner2003pf}
K.~Glasner.
\newblock Spreading of droplets under the influence of intermolecular forces.
\newblock {\em Physics of Fluids}, 15(7):1837--1842, 2003.

\bibitem{glasner2003coarsening}
K.~B. Glasner and T.~P. Witelski.
\newblock Coarsening dynamics of dewetting films.
\newblock {\em Physical Review E}, 67(1):016302, 2003.

\bibitem{grant1993spinodal}
C.~P. Grant.
\newblock Spinodal decomposition for the {Cahn-Hilliard} equation.
\newblock {\em Communications in Partial Differential Equations},
  18(3-4):453--490, 1993.

\bibitem{iakubovich1975linear}
V.~A. Iakubovich and V.~M. Starzhinski{\u\i}.
\newblock {\em Linear differential equations with periodic coefficients},
  volume~2.
\newblock Wiley, 1975.

\bibitem{israel2013wella}
H.~Israel.
\newblock Well-posedness and long time behavior of an {A}llen--{C}ahn type
  equation.
\newblock {\em Communications on Pure \& Applied Analysis}, 12(6), 2013.

\bibitem{israel2013well}
H.~Israel, A.~Miranville, and M.~Petcu.
\newblock Well-posedness and long time behavior of a perturbed
  {C}ahn--{H}illiard system with regular potentials.
\newblock {\em Asymptotic Analysis}, 84(3-4):147--179, 2013.

\bibitem{ji2017thin}
H.~Ji.
\newblock {\em Thin films with non-conservative effects}.
\newblock PhD thesis, Duke University, 2017.

\bibitem{ji2017finite}
H.~Ji and T.~P. Witelski.
\newblock Finite-time thin film rupture driven by modified evaporative loss.
\newblock {\em Physica D}, 342:1--15, 2017.

\bibitem{ji2018instability}
H.~Ji and T.~P. Witelski.
\newblock Instability and dynamics of volatile thin films.
\newblock {\em Physical Review Fluids}, 3(2):024001, 2018.

\bibitem{karali2007role}
G.~Karali and M.~A. Katsoulakis.
\newblock The role of multiple microscopic mechanisms in cluster interface
  evolution.
\newblock {\em Journal of Differential Equations}, 235(2):418--438, 2007.

\bibitem{karali2010}
G.~Karali and T.~Ricciardi.
\newblock On the convergence of a fourth order evolution equation to the
  {A}llen-{C}ahn equation.
\newblock {\em Nonlinear Anal.}, 72(11):4271--4281, 2010.

\bibitem{karali2012existence}
G.~D. Karali and Y.~Nagase.
\newblock On the existence of solution for a {C}ahn-{H}illiard/{A}llen--{C}ahn
  equation.
\newblock {\em Discrete and Continuous Dynamical Systems-Series S}, 2012.

\bibitem{munch2008}
M.~D. Korzec, P.~L. Evans, A.~M\"{u}nch, and B.~Wagner.
\newblock Stationary solutions of driven fourth- and sixth-order
  {C}ahn-{H}illiard-type equations.
\newblock {\em SIAM J. Appl. Math.}, 69(2):348--374, 2008.

\bibitem{kuznetsov}
Y.~A. Kuznetsov.
\newblock {\em Elements of applied bifurcation theory}, volume 112 of {\em
  Applied Mathematical Sciences}.
\newblock Springer-Verlag, New York, third edition, 2004.

\bibitem{LP}
R.~S. Laugesen and M.~C. Pugh.
\newblock Linear stability of steady states for thin film and {Cahn-Hilliard}
  type equations.
\newblock {\em Archive for rational mechanics and analysis}, 154(1):3--51,
  2000.

\bibitem{laugesen2000properties}
R.~S. Laugesen and M.~C. Pugh.
\newblock Properties of steady states for thin film equations.
\newblock {\em European Journal of Applied Mathematics}, 11(03):293--351, 2000.

\bibitem{laugesen2002energy}
R.~S. Laugesen and M.~C. Pugh.
\newblock Energy levels of steady states for thin-film-type equations.
\newblock {\em Journal of Differential Equations}, 182(2):377--415, 2002.

\bibitem{lyushnin2002fingering}
A.~V. Lyushnin, A.~A. Golovin, and L.~M. Pismen.
\newblock Fingering instability of thin evaporating liquid films.
\newblock {\em Physical Review E}, 65(2):021602, 2002.

\bibitem{miranville2017cahn}
A.~Miranville.
\newblock The {Cahn--Hilliard} equation and some of its variants.
\newblock {\em AIMS Mathematics}, 2(3):479--544, 2017.

\bibitem{Mitlin}
V.~S. Mitlin.
\newblock Dewetting of solid surface: Analogy with spinodal decomposition.
\newblock {\em J. Colloid and Interface Science}, 156:491--497, 1993.

\bibitem{moosman1980evaporating}
S.~Moosman and G.~M. Homsy.
\newblock Evaporating menisci of wetting fluids.
\newblock {\em Journal of Colloid and Interface Science}, 73(1):212--223, 1980.

\bibitem{munch2005dewetting}
A.~M{\"u}nch.
\newblock Dewetting rates of thin liquid films.
\newblock {\em Journal of Physics: Condensed Matter}, 17(9):S309, 2005.

\bibitem{munch2005lubrication}
A.~M{\"u}nch, B.~A. Wagner, and T.~P. Witelski.
\newblock Lubrication models with small to large slip lengths.
\newblock {\em Journal of Engineering Mathematics}, 53(3-4):359--383, 2005.

\bibitem{murdock}
J.~A. Murdock.
\newblock {\em Perturbations: Theory and methods}, volume~27 of {\em Classics
  in Applied Mathematics}.
\newblock Society for Industrial and Applied Mathematics (SIAM), Philadelphia,
  PA, 1999.

\bibitem{myers}
T.~G. Myers.
\newblock Thin films with high surface tension.
\newblock {\em {SIAM} {R}eview}, 40(3):441--462, 1998.

\bibitem{ANC1988}
A.~Novick-Cohen.
\newblock Energy methods for the {C}ahn-{H}illiard equation.
\newblock {\em Quart. Appl. Math.}, 46(4):681--690, 1988.

\bibitem{ANC1984}
A.~Novick-Cohen and L.~A. Segel.
\newblock Nonlinear aspects of the {C}ahn-{H}illiard equation.
\newblock {\em Physica D}, 10(3):277--298, 1984.

\bibitem{novick2010thin}
A.~Novick-Cohen and A.~Shishkov.
\newblock The thin film equation with backwards second order diffusion.
\newblock {\em Interfaces and Free Boundaries}, 12:463--496, 2010.

\bibitem{OckendonOckendon}
J.~R. Ockendon and H.~Ockendon.
\newblock {\em Viscous flow}.
\newblock Cambridge University, Cambridge, 1995.

\bibitem{oron1999dewetting}
A.~Oron and S.~G. Bankoff.
\newblock Dewetting of a heated surface by an evaporating liquid film under
  conjoining/disjoining pressures.
\newblock {\em Journal of Colloid and Interface Science}, 218(1):152--166,
  1999.

\bibitem{oron2001dynamics}
A.~Oron and S.~G. Bankoff.
\newblock Dynamics of a condensing liquid film under conjoining/disjoining
  pressures.
\newblock {\em Physics of Fluids}, 13(5):1107--1117, 2001.

\bibitem{oron1997long}
A.~Oron, S.~H. Davis, and S.~G. Bankoff.
\newblock Long-scale evolution of thin liquid films.
\newblock {\em Reviews of Modern Physics}, 69(3):931, 1997.

\bibitem{peletier}
L.~A. Peletier and W.~C. Troy.
\newblock {\em Spatial patterns: Higher order models in physics and mechanics},
  volume~45 of {\em Progress in Nonlinear Differential Equations and their
  Applications}.
\newblock Birkh\"{a}user Boston, Inc., Boston, MA, 2001.

\bibitem{perazzo2017analytical}
C.~A. Perazzo, J.~R. Mac~I., and J.~M. Gomba.
\newblock Analytical solutions for the profile of two-dimensional droplets with
  finite-length precursor films.
\newblock {\em Physical Review E}, 96(6):063109, 2017.

\bibitem{taranets2014unstable}
R.~M. Taranets and J.~R. King.
\newblock On an unstable thin-film equation in multi-dimensional domains.
\newblock {\em Nonlinear Differential Equations and Applications},
  21(1):105--128, 2014.

\bibitem{Taranets2003}
R.~M. Taranets and A.~E. Shishkov.
\newblock The effect of time delay of support propagation in equations of thin
  films.
\newblock {\em Ukra\"{i}n. Mat. Zh.}, 55(7):935--952, 2003.

\bibitem{teletzke1988wetting}
G.~F. Teletzke, H.~T. Davis, and L.~E. Scriven.
\newblock Wetting hydrodynamics.
\newblock {\em Revue de Physique Appliquee}, 23(6):989--1007, 1988.

\bibitem{teschl2012ordinary}
G.~Teschl.
\newblock {\em Ordinary differential equations and dynamical systems}.
\newblock American Mathematical Soc., 2012.

\bibitem{thiele2010thin}
U.~Thiele.
\newblock Thin film evolution equations from (evaporating) dewetting liquid
  layers to epitaxial growth.
\newblock {\em Journal of Physics: Condensed Matter}, 22(8):084019, 2010.

\bibitem{Thiele2014}
U.~Thiele.
\newblock Patterned deposition at moving contact lines.
\newblock {\em Advances in Colloid and Interface Science}, 206:399--413, 2014.

\bibitem{Thiele2018}
U.~Thiele.
\newblock Recent advances in and future challenges for mesoscopic hydrodynamic
  modelling of complex wetting.
\newblock {\em Colloids and Surfaces A}, 553:487--495, 2018.

\bibitem{thiele2016}
U.~Thiele, A.~J. Archer, and L.~M. Pismen.
\newblock Gradient dynamics models for liquid films with soluble surfactant.
\newblock {\em Physical Review Fluids}, 1(083903):1--31, 2016.

\bibitem{todorova2012relation}
D.~Todorova, U.~Thiele, and L.~M. Pismen.
\newblock The relation of steady evaporating drops fed by an influx and freely
  evaporating drops.
\newblock {\em Journal of Engineering Mathematics}, 73(1):17--30, 2012.

\bibitem{witelski2000dynamics}
T.~P. Witelski and A.~J. Bernoff.
\newblock Dynamics of three-dimensional thin film rupture.
\newblock {\em Physica D}, 147(1):155--176, 2000.

\bibitem{zhang2016existence}
X.~Zhang and C.~Liu.
\newblock Existence of solutions to the {C}ahn-{H}illiard/{A}llen--{C}ahn
  equation with degenerate mobility.
\newblock {\em Electron. J. Differential Equations}, 2016(329):1--22, 2016.

\end{thebibliography}
\end{document}